\documentclass[11pt]{amsart}
\usepackage{amsmath, amsfonts, amsthm, amssymb}
\usepackage[all,cmtip,line]{xy}\usepackage{mathtools}
\usepackage{hyperref} 
\usepackage{array}
\usepackage{enumerate}
\usepackage{srcltx} 

\usepackage{todonotes}
\newtheorem{thm}{Theorem}[section]
\newtheorem{lem}{Lemma}[section]
\newtheorem{prop}[lem]{Proposition}
\newtheorem{cor}[lem]{Corollary}

\newtheorem{defn}[lem]{Definition}

\numberwithin{equation}{section}

\newcommand{\bR}{ \mathbb{R}} 
\newcommand{\bN}{ \mathbb{N}} 
\newcommand{\bC}{ \mathbb{C}} 

\newcommand{\diam}{ \mbox{diam}}

\newcommand{\parallelsum}{\mathbin{\!/\mkern-5mu/\!}}

\newcommand \om{\omega}

\newlength{\originalbase}
\title{The uniqueness of Poincar\'e type extremal K\"ahler metric}

\author{Yulun Xu}  
\address{Mathematics Department and the Institute for Mathematical Sciences, Stony Brook University, Stony Brook NY, 11794-3651, USA}
\email{yulun.xu@stonybrook.edu}

\date{}

\begin{document}

\maketitle

\begin{abstract}

Let $D$ be a smooth divisor on a closed K\"ahler manifold $X$. Suppose that $Aut_0(D)=\{Id\}$. We prove that the Poincar\'e type extremal K\"ahler metric with a cusp singularity at $D$ is unique up to a holomorphic transformation on $X$ that preserves $D$. This generalizes Berman-Berndtson's work \cite{BB} on the uniqueness of extremal K\"ahler metrics from closed manifolds to some complete and noncompact manifolds.
\end{abstract}

\tableofcontents

\section{Introduction}
 Let $(V,\omega_X)$ be a compact K\"ahler manifold, then we can define the space of K\"ahler potentials in a K\"ahler class $[\omega_X]$ as:
\begin{equation*}
    \mathcal{H}=\{\varphi \in C^{\infty}(V):\omega_{\varphi}=\omega_X+\sqrt{-1}\partial \bar{\partial}\varphi >0 \text{ on  }V\}.
\end{equation*}
 Locally in a holomorphic coordinate chart, the K\"ahler form $\omega_{\varphi}$ can be written as $$
\omega_\varphi = g_{\varphi,\alpha \bar{\beta}}\sqrt{-1} d\,z^\alpha \wedge \overline{d\, z^\beta} =  \left(g_{\alpha \bar{\beta}} + \frac{\partial^2 \varphi}{\partial z_\alpha \bar \partial z_\beta}\right) \sqrt{-1} d\,z^\alpha \wedge \overline{d\, z^\beta}. $$
Then, its scalar curvature $R_\varphi$  is defined as:
\[
R_\varphi = - g_\varphi^{\alpha\bar \beta} \frac{\partial^2}{\partial z_\alpha \partial \bar{z}_\beta}\log \det(g_{\varphi,i\bar j}).
\]

The central problem in K\"ahler geometry which goes back to Calabi's program \cite{C2} \cite{C3} is to find a K\"ahler metric as a canonical representative in a given K\"ahler class. A candidate for such representative is extremal K\"ahler metrics, which are critical points of the Calabi functional which is defined by:
\begin{equation*}
\Phi(\omega_{\varphi})= \int_X R_{\varphi}^2 \omega_{\varphi}^n.
\end{equation*}
Another characterization of the extremal K\"ahler metric is that the scalar curvature $R$ of the extremal K\"ahler metric $\omega$ satisfies $R^{,\alpha \beta}=0$, where the covariant derivatives are taken with respect to $\omega$. An extremal K\"ahler metric is called a cscK metric if its scalar curvature is a constant function.

According to the Yau-Tian-Donaldson conjecture, the existence of an extremal metric is expected to be equivalent to K-stability (\cite{Do},\cite{T}, \cite{Y}). It was proved by Chen-Cheng \cite{CC} that the properness of K-energy implies the existence of a cscK metric. An interesting question is what will happen in the unstable case. Donaldson predicted in \cite{Do}, \cite{Do1}, \cite{Do2} there exists a divisor $D$ such that one can find a complete extremal metric on $X\setminus D$. Therefore, it is important to study complete extremal metrics on $X\setminus D$ and the goal of this paper is to generalize the uniqueness of extremal K\"ahler metrics from closed manifolds to some complete and noncompact manifolds. Our setting is as follows: On $\mathbb{C}^n$, we can write down the standard local model for the Poincar\'e type K\"ahler metric (cusp metric): 
\begin{equation}\label{standard cusp}
    \omega_0=\frac{\sqrt{-1} dz^n \wedge d\bar z^n}{2|z_n|^2 log^2(|z_n|)} +\Sigma_{i=1}^{n-1}\sqrt{-1}dz^i \wedge d \bar z^i.
\end{equation}
The above metric is an ideal model for constant scalar curvature K\"ahler metrics because its scalar curvature is constant. The above metric can also be seen as a limit of conical metric as the cone angle goes to zero, c.f. \cite{G}, \cite{Ao}.

Then we can define the Poincar\'e type K\"ahler metric. Let $(X,\omega_X)$ be a closed K\"ahler manifold, let $D$ be a smooth divisor on $X$. Let $D=\Sigma_{j=1}^N D_j$ be the decomposition of $D$ into smooth irreducible components. We can define the Poincar\'e type K\"ahler metric as follows:
\begin{defn}
    We say that $\omega$ is a Poincar\'e type K\"ahler metric, of class $\Omega=[\omega_X]_{dR}$, if for any point $p\in D$, and any holomorphic coordinate $U$ of $X$ around $p$ such that in the coordinate $\{D=0\}=\{z_n=0\}$(We call this kind of coordinate cusp coordinate from now on), $\omega$ satisfies:
    \begin{enumerate}
        \item There exists a constant $C$ such that $\frac{1}{C}\omega_0 \le \omega \le C \omega_0$ holds in $U$\\
        \item There exists a function $\varphi$ such that $\omega=\omega_X +dd^c \varphi$. There exist  constants $C(k)$ such that in $U$, $|\nabla^k_{\omega_0}\varphi|_{\omega_0}\le C(k)$ for any $k\ge 1$. Moreover, $\varphi=O(log(-log|z_n|))$.
        \item $\omega$ is a smooth K\"ahler metric on $X\setminus D$.
    \end{enumerate}
\end{defn}
The interesting point of the Poincar\'e type K\"ahler metric is that $(X\setminus D, \omega)$ is complete but non-compact. Compared with the closed manifolds,  many interesting new phenomena appear as a result.  A lot of progress has been made in this case. Auvray proves in \cite{A} the existence of the Poincar\'e type $C^{1,1}$ geodesic. Under the assumption that $K_X[D]$ is ample, he proved that the Poincar\'e type cscK metric is unique. He also discovered a topological constraint for the Poincar\'e type cscK metrics in \cite{A2}, asymptotic properties of Poincar\'e type extremal K\"ahler metrics in \cite{A3} and the Poincar\'e type Futaki characters in \cite{A4}. Sektnan \cite{S} and Feng \cite{F} used gluing arguments to construct examples of Poincar\'e type extremal K\"ahler metrics. Aoi \cite{Ao} proved that Poincar\'e type cscK metrics can be approximated by conical K\"ahler metrics under assumptions on holomorphic vector fields on $D$ and $X$.

Denote $Aut_0^D(X)$ as the identity component of $\{g\in Aut_0(X): g(D)=D\}$.  For the Poincar\'e type extremal K\"ahler metric, we can prove the following Theorem:
\begin{thm}\label{unique extk}
Suppose that  $D$ is a smooth divisor on $X$ and  $Aut_0(D)=\{Id\}$.  Given two Poincar\'e type extremal K\"ahler metrics $\omega_1$ and $\omega_2$ in a given cohomology class $[\omega]$. Then there exists an element $g\in Aut_0^D(X)$ such that $g^* \omega_1 = \omega_2$.
\end{thm}
In our previous paper \cite{XZ}, we proved the uniqueness of Poincar\'e type cscK metrics under the same assumption as the above Theorem.  

The above Theorem is implied by the following openness of a continuity path:
\begin{thm}\label{main theorem}
Suppose that  $D$ is a smooth divisor on $X$ and $Aut_0(D)=\{Id\}$.  Given two Poincar\'e type extremal K\"ahler metrics $\omega_1$ and $\omega_2$ in a given cohomology class $[\omega]$. Then, there exists $\epsilon>0$ with a smooth function $\phi: (1-\epsilon,1]\times (X\setminus D) \rightarrow \bR$ such that $\varphi_{t_1} \triangleq \phi(t_1,\cdot) \in \widetilde{\mathcal{PM}_{[\omega]}}$ for $t_1\in (1-\epsilon,1]$  and
\begin{equation*}
\nabla_{\varphi_t}^{1,0}(R_{\varphi_t}-(1-t)tr_{\varphi_t}\omega_2) \in \mathbf{h}_{\parallelsum}^D.
\end{equation*}
Moreover, there exists $g\in Aut_0^D(X)$ such that $g^* (\omega +dd^c \varphi_1)=\omega_1$.
\end{thm}
In the above, we denote $\mathcal{PM}_{[\omega]}$ as the space of Poincar\'e type K\"ahler metrics in the cohomology class $[\omega]$. Denote $\widetilde{\mathcal{PM}_{[\omega]}}$ as the space of K\"ahler potentials for $\mathcal{PM}_{[\omega]}$ with respect to a given background metric $\omega$. 
We also denote $\mathbf{h}^D_{\parallelsum}$ as the set of holomorphic vector fields on $X$ that are parallel to $D$.

For smooth metrics on a closed K\"ahler manifold $X$, Berman-Berndtsson \cite{BB} and Chen-Paun-Zeng \cite{CPZ} proved the uniqueness of cscK metrics and extremal K\"ahler metrics. The idea of proving the uniqueness of cscK metric is as follows: First, we prove the convexity of the K-energy along $C^{1,1}$ geodesics in the space of K\"ahler potentials whose existence was given by Chen in \cite{C5}.  Suppose that we have that the K-energy is strictly convex, we can get that the cscK metric is unique as its critical point. However, this is not always the case. Instead, we perturb the K-energy by another functional, which is strictly convex. We call the new functional twisted K-energy.  Using a bifurcation argument, which is a version of implicit function theory, we can show that near a given cscK metric, we can get a critical point of the twisted K-energy.  Since twisted K-energy is strictly convex, its critical point is unique. Then, we can take a limit from the critical points of twisted K-energy to the critical points of K-energy to prove the uniqueness of the cscK metric. As for the proof of the uniqueness of extremal K\"ahler metrics, Calabi \cite{C4} proved that the isometry group of an extremal K\"ahler metric is a maximal compact connected subgroup of $Aut_0(X)$. Using this fact, we can show that the extremal K\"ahler vector fields of two extremal K\"ahler metrics are the same after we pull back one of the extremal K\"ahler metrics by an element in $Aut_0(X)$. Then, we can define a modified K-energy by the K-energy and the extremal K\"ahler vector field such that the extremal K\"ahler metric is a critical point of the modified K-energy. The rest of the proof is similar to the uniqueness of cscK metrics. 

One key part in the proof of the Theorem \ref{main theorem} is the solvability of $Re L$, where $L$ is the Lichnerowicz operator. Fix a metric $\omega$, the Lichnerowicz operator is defined by
\begin{equation*}
L u = 2 u^{,\alpha \beta}_{\,\,\,\,\,\,\,\,\, \beta \alpha},
\end{equation*}
where all the covariant derivatives are take with respect to $\omega$. By changing the order of covariant derivatives, we have that
\begin{equation*}
\begin{split}
L u &= 2u^{,\alpha \beta}_{\,\,\,\,\,\,\,\,\, \beta \alpha} = 2u^{,\alpha \beta}_{\,\,\,\,\,\,\,\,\, \alpha \beta} = 2u ^{,\alpha \,\,\,\beta}_{\,\,\,\, \alpha \,\,\,\, \beta} +2(u^{,\alpha} R_{\alpha}^{\,\,\,\, \beta})_{,\beta} \\
& = \frac{1}{2}\Delta^2 u + <dd^c u, Ric_{\omega}> + \frac{1}{2}u^{,\beta} R_{,\beta}.
\end{split}
\end{equation*}
Here $R_{\alpha \bar \beta}$ is the covariant component of the Ricci tensor. Here $dd^c u =-d(Jdu)$ and $\Delta u =\frac{n dd^c u \wedge \omega^{n-1}}{\omega^n}$. The above notations align with the notations in \cite{A3}. If $\omega$ is a cscK metric, then $R$ is a constant which implies  $u^{,\beta} R_{,\beta}=0$. Then $L$ is a real operator which is the linearized operator of the cscK equation:
\begin{equation*}
R= \underline{R}.
\end{equation*}
If $\omega$ is not a cscK metric,  $L$ may not be real. We can consider $Re L$ instead.

\begin{defn}
We say that a Poincar\'e type K\"ahler metric is asymptotic to a product metric, if there exists a K\"ahler metric $\omega_D$ such that there exist  constants $\eta>0$ and $a>0$ such that in any cusp coordinate:
\begin{equation}\label{e 5.1}
\omega= p^* \omega_D + \frac{2 \sqrt{-1} ad z^n \wedge d \bar z^n}{|z_n|^2 \log^2 |z_n|^2} +O(e^{-\eta t}).
\end{equation}
\end{defn}
According to the Lemma \ref{asym} proved in \cite{A3}, Poincar\'e type extremal K\"ahler metrics, including Poincar\'e type cscK metrics, all satisfy (\ref{e 5.1}).

Then we can prove the following Solvability of $Re L$:
\begin{prop}\label{solve l operator1}
Suppose that $\omega$ is a Poincar\'e type K\"ahler metric satisfying (\ref{e 5.1}). Then there exists a constant $0< \delta_1< \frac{1}{2}$. For any $\eta_0 \in (0, \delta_1)$, for any $f\in C_{-\eta_0}^{1,\alpha}$  such that $\int_{M\setminus D}f u \omega^n=0$ for any $u\in \overline{\mathbf{h}^D_{\parallelsum,\mathbb{R}}}$, we can find a function $v\in C_{-\eta_0}^{5,\alpha} \oplus \chi(t) p^* Ker Re L_D$ such that $Re Lv =f$. 
\end{prop}
In the above,  $t$ is the function given in the section 3.7 expressing the normal direction of $D$. We have that
\begin{equation*}
\lim_{ d(x,D) \rightarrow 0} t(x)=\infty.
\end{equation*}
 $\chi(t)$ is a smooth function defined on $[0,+\infty)$ with $\chi(t)=0$ for $t\in [0,\delta_0]$ and $\chi(t)=1$ for $t \ge 2 \delta_0$ with $\delta_0$ to be a small positive constant. $C_{\delta}^{k,\alpha}$ is the weighted H\"older space defined in the section 3.4.  $p$ is the projection map from a neighbourhood of $D$ to $D$ defined in the section 3.7. 
$L_{D}$ is the Lichnerowicz operator of $\omega_D$. 

In our previous paper \cite{XZ}, we proved a weaker result in the sense that the solution $v$ we get lies in $C_{-\eta_1}^{5,\alpha} \oplus \chi(t) Ker Re L_D$ for some $0< \eta_1 < \eta_2$. The proposition \ref{solve l operator1} is sharp because after modding out functions in a finite dimensional space $\chi(t) Ker Re L_D$, the solution has the same decay rate as the right-hand side of the equation. This enables us to use the implicit function theory in many situations. For example, this is used in our proof of the Theorem \ref{main theorem},  Feng's work on the gluing construction of Poincar\'e type cscK metric \cite{F}, Aoi's work in the approximation of Poincar\'e type cscK metric with conical cscK metrics\cite{A} and Sektnan's work on the blow up of Poincar\'e type extremal K\"ahler metrics\cite{S}. Note that Sektnan claimed the solvability of $L$, but unfortunately there is a gap in his proof of Proposition 4.3 about the Fredholm index of the Lichnerowicz operator. In that place he used the result of Lockhart-McOwen\cite{LM}: Let $M$ be a manifold with a cylindrical end, i.e. $$M= D\times [0,+\infty) \cup M_2,$$ where $D$ is a closed manifold and $M_2$ is a compact manifold with boundary. Then we can study the global Fredholm index using the Fredholm index of the same operator restricted to $D\times [0,\infty)$. 

The gap is that Lockhart-McOwen\cite{LM} study manifolds with cylindrical ends. But the manifolds with Poincar\'e type K\"ahler metrics don't have cylindrical ends. We need to mod out a $S^1$ action near the divisor to get a cylindrical end which is elaborated in the section 3.7. It is unclear that how the Fredholm index changes when we mod out a $S^1$ action. As a result, we use another way to prove the Proposition \ref{solve l operator} without using the Fredholm index of $Re L$ at all.

If $\omega$ is a smooth K\"ahler metric on $X$, then the Lichnerowicz operator $L$ corresponding to $\omega$ satisfies the followsing equation using Fredholm alternative:
\begin{equation}\label{smooth dec}
C^{k,\alpha}= Ker Re L |_{C^{k,\alpha}} \oplus Re L (C^{k,\alpha}).
\end{equation}
For any $\delta \in \bR$, we can define the following space:
\begin{equation*}
\widetilde{C}^{k,\alpha}_{\delta}(X \setminus D) \triangleq C_{\delta}^{k,\alpha} (X \setminus D) \oplus \chi p^* C^{k,\alpha}(D) .
\end{equation*}
Here $C_{\delta}^{k,\alpha}$ is a weighted H\"older space defined in the section 3.  Using the Proposition \ref{solve l operator}, we can prove an equation similar to (\ref{smooth dec}) for Poincar\'e type K\"ahler metrics:
\begin{thm}\label{holder dec 12}
Suppose that $\omega$ is a Poincar\'e type K\"ahler metric satisfying (\ref{e 5.1}) with $\eta < \frac{1}{2}$. Then there exists a constant $\delta_1>0$ such that for any $\eta_0 \in (0, \delta_1)$, we have that:
\begin{equation*}
\widetilde{C}_{-\eta_0}^{1,\alpha}= Ker Re L |_{\widetilde{C}_{-\eta_0}^{5,\alpha}} \oplus Re L(t \chi (p^* Ker Re L_D) ) \oplus Re L(\widetilde{C}_{-\eta_0}^{5,\alpha}).
\end{equation*}
\end{thm}
The above Theorem is a key part when we use the implicit function theorem in the proof of the Theorem \ref{main theorem}.

Denote $Iso_0^D (X,\omega)$ as the identity component of $\{g \in Aut(X): g(D)=D, g^* \omega= \omega\}$.  Another key part in the proof of the Theorem \ref{main theorem} is the following Theorem:

\begin{thm}\label{compact isometry}
Suppose that $D$ is a smooth divisor. Suppose that $Aut_0(D)=\{Id\}$. Let $\omega$ be a Poincar\'e type extremal K\"ahler metric. Then  $ Iso_0^D (X,\omega)$ is a maximal compact connected subgroup in $Aut_0^D(X)$.
\end{thm}

The above Proposition was proved by Lichnerowicz in \cite{L} for smooth cscK metrics and Calabi in \cite{C4} for smooth extremal K\"ahler metrics, both on closed manifolds. Note that in these cases, the compactness of the isometry group is not a problem. The completeness and noncompactness of $(X\setminus D, \omega)$ make it much harder to prove in the Poincar\'e type case.

Firstly, we want to uniformly control the behavior of elements in $ Iso_0^D (X,\omega)$ away from $D$. We prove that for any compact set $K\subset X\setminus D$, there exists a compact set $K' \subset X\setminus D$ such that $g (K) \subset K'$ and $g^{-1}(K) \subset K'$ (see the Proposition \ref{compact to compact}) for any $g\in Iso_0^D(X,\omega)$. This helps us rule out the following situation: there is a point $q\in X \setminus D$ and a sequence $\{g_k\}$ in $ Iso_0^D (X,\omega)$ such that $\lim_{k \rightarrow \infty}g_k(q) =q_0$ for some $q_0 \in D$.  Thus, we can prove that for any sequence $\{g_k\}$ in $Iso_0^D(X,\omega)$, we can take a subsequence of $\{g_k\}$ (still denoted as $\{g_k\}$) converging locally uniformly on $X\setminus D$ to a map $g$ which is a holomorphic transformation of $X\setminus D$. Secondly, in order to show that $g$ can be continuously extended to $D$ such that $g\in Iso_0^D(X,\omega)$ and $g_k$ converging uniformly to $g$ on $X$, we need to uniformly control the behavior of elements in $Iso_0^D (X,\omega)$ near $D$ (see the Proposition \ref{gk converge}). We develop a geodesic technique to achieve this. 

One application of the Theorem \ref{compact isometry} is that we can characterize the asymptotic behaviour of Poincar\'e type extremal K\"ahler metrics:

\begin{thm}\label{a for extremal metric}
Suppose that $Aut_0(D)=\{Id\}$ and $D$ is a smooth divisor. Let $\omega_3= \omega + dd^c \varphi_3, \omega_4= \omega+ dd^c \varphi_4$ be any two Poincar\'e type extremal K\"ahler metrics in the same cohomology class. Then we have that
\begin{equation*}
a_j(\omega_1)= a_j(\omega_2)
\end{equation*}
for any $j \le N$.
\end{thm}
Note that the above theorem was proved by Auvray in \cite{A3} for Poincar\'e type cscK metrics, see the Lemma \ref{asym}. The constants $a_j$ in the above theorem are defined in the Lemma \ref{asym}. $a_j$ basically characterize the behaviour of a Poincar\'e type extremal K\"ahler metric in the direction perpendicular to $D_j$.

In the section 3, we introduce some background knowledge about Poincar\'e type K\"ahler metrics and clarify some notations. In the section 4, we prove the Proposition \ref{solve l operator1} and the Theorem \ref{holder dec 12}. In the section 5, we prove the compactness of isometry group. In the section 6, we prove the Theorem \ref{compact isometry}. In the section 7, we prove that the isometry group can determine extremal K\"ahler vector fields. In the section 8, we prove the Theorem \ref{main theorem}, the Theorem \ref{unique extk} and the Theorem \ref{a for extremal metric}.
\section{Acknowledgement}
This project was suggested by Prof. Xiuxiong Chen.  The author thanks his advisors Prof. Xiuxiong Chen and Prof. Jingrui Cheng for their suggestions on this paper and Yueqing Feng for her feedback on the original version of the paper. The author also thanks Junbang Liu for providing a proof of the Proposition \ref{compact to compact} which is different from the author's original proof. This helps remove the additional assumption that $D$ is connected. This research is partially funded by the Simons Foundation.
\section{Preliminaries}

\subsection{Background metric of Poincar\'e type}
First, we can construct a Poincar\'e type K\"ahler metric and use it as a background metric.   We take a holomorphic defining section $\sigma \in (\mathcal{O([D])},|\cdot|)$ for $D$. Then we define 
$$\rho \triangleq -\log(|\sigma|^2)\ge 1$$ 
out of $D$, equivalently, $|\sigma|^2\leq e^{-1}$. Let $\lambda$ be a nonnegative real constant to be determined. Then we set 
$$\mathbf{u}\triangleq \log(\lambda+\rho).$$ We denote 
$$\omega\triangleq\omega_X-Ai \partial \bar \partial \mathbf{u}$$ which is used as a background metric. 

Auvray shows in \cite[Lemma 1.1]{A} that for any $A>0$ and for sufficiently large $\lambda$ depending on $A$ and $\omega_X$, the $(1,1)$-form $\omega_X-Ai \partial \bar \partial \mathbf{u}$ is a Poincar\'e type K\"ahler metric. 

\subsection{Asymptotic behaviour of Poincar\'e type extremal K\"ahler metrics}
Define 
\begin{equation*}
\underline{R} =-4\pi n \frac{c_1(K_X[D])\cdot [\omega]^{n-1}}{[\omega]^n} \text{ and }\underline{R}_{D_j}=-4\pi n \frac{c_1(D_j)\cdot c_1(K_X[D])\cdot [\omega]^{n-2}}{c_1(D_j)\cdot [\omega_X]^{n-1}}.
\end{equation*}
Auvray proved the asymptotic behaviors of Poincar\'e type extremal (constant scalar curvature) K\"ahler metrics in \cite{A3}.
\begin{lem}\label{asym}
Assume that $\omega$ is a Poincar\'e type extremal (constant scalar curvature)  K\"ahler metric of class $[\omega]$ on the complement of a (smooth) divisor $D=\Sigma_{j=1}^N D_j$ with disjoint components in a compact K\"ahler manifold $(X,\omega).$ Then for all $j$ there exist constants $a_j, \eta>0$, and an extremal (constant scalar curvature) K\"ahler metric $\omega_j \in [\omega|_D]$ such that on any open subset $U$ of coordinates $(z^1,z^2,...,z^m)$ such that $U \cap D_j =\{z^n=0\}$, then $\omega=\frac{a_j \sqrt{-1} dz^n \wedge d\bar z^n}{2 |z^n|^2 log^2(|z^n|)}+p^* \omega_j +O(|log(|z^n|)|^{-\eta})$ as $z^n \rightarrow 0$. Moreover, if $\omega$ is a Poincar\'e type cscK metric, then $a_j=\frac{2}{\underline{R}_{D_j}-\underline{R}}$.
\end{lem}

\subsection{Quasi coordinates}
Next, the quasi coordinates, see \cite{TY}, is used in \cite{A} to define function spaces using Poincar\'e type K\"ahler metrics.
Let $\Delta$ be a unit disc and let $\Delta^*$ be a punctured unit disc. For any $\delta \in (0,1)$, we can set $$\varphi_{\delta}: \frac{3}{4}\Delta \rightarrow \Delta^*,\quad\xi \mapsto exp(-\frac{1+\delta}{1-\delta}\frac{1+\xi}{1-\xi}).$$ For any $\delta \in (0,1)$ and any Poincar\'e type K\"ahler metric $\omega$, $\varphi_{\delta}^* \omega$ is quasi-isometric to the Euclidean metric. Then we can take \begin{align*}
\Phi_{\delta}: \mathcal{P} \triangleq \Delta^{n-1} \times (\frac{3}{4}\Delta) &\rightarrow \Delta^{n-1} \times \Delta^*,\quad \delta \in (0,1),\\
 (z_{1},...,z_{n-1},\xi)&\mapsto (z_1,...,z_{n-1},\phi_{\delta}(\xi)).
 \end{align*} 
 We say that a holomorphic coordinate of $X$ is a cusp coordinate if in this coordinate we have $D=\{z_n=0\}$. Let us prove a lemma using the quasi coordinate.
\begin{lem}\label{lem 2.1}
    Let $\omega_X$ be a smooth K\"ahler metric on $X$. Then in any cusp coordinate and for any $k\ge 1$, we have that $|\nabla_{\omega_0}^k \omega_X|_{\omega_0}\le C(k)$ for some constant $C(k)$. Here $\omega_0$ is the standard local Poincar\'e type K\"ahler metric given by (\ref{standard cusp}).
\end{lem}
\begin{proof}
    Using direct calculation, we have that 
    \begin{equation*}
        \Phi_{\delta}^* \omega_{0}= \frac{\sqrt{-1}d \xi \wedge d \bar \xi}{(1-|\xi|^2)^2} +\Sigma_{i=1}^{n-1} \sqrt{-1} dz^i \wedge d \bar z^i.
    \end{equation*}
    This is $C^{\infty}$ quasi-isometric to the Euclidean metric on $\frac{3}{4}\Delta$. In a holomorphic coordinate of $X$, we can write $\omega_X$ as $\omega_X = \Sigma_{i,j} a_{ij} \sqrt{-1} dz^i \wedge d\bar z^j$. Then we have that:
    \begin{equation*}
    \begin{split}
         \Phi_{\delta}^* \omega_X &= \Sigma_{\alpha,\beta} a_{\alpha \beta}(\Phi_{\delta} (z',\xi)) \sqrt{-1} dz^{\alpha} \wedge d \bar z^{\beta} \\
         &+ \Sigma_{\alpha} a_{\alpha n}(\Phi_{\delta} (z',\xi)) \sqrt{-1} dz^{\alpha} \wedge \overline{exp(-\frac{1+\delta}{1-\delta}\frac{1+\xi}{1-\xi})}(-\frac{1+\delta}{1-\delta}\frac{2}{(1-\bar \xi)^2}) d\bar \xi \\
         & + \Sigma_{\beta} a_{n \beta} (\Phi_{\delta} (z',\xi)) \sqrt{-1}exp(-\frac{1+\delta}{1-\delta}\frac{1+\xi}{1-\xi}) (-\frac{1+\delta}{1-\delta}\frac{2}{(1-\xi)^2}) d \xi \wedge d\bar z^{\beta}\\
         & + a_{nn}(\Phi_{\delta} (z',\xi)) \sqrt{-1}exp(-2\frac{1+\delta}{1-\delta}Re\frac{1+\xi}{1-\xi}) (\frac{(1+\delta)^2}{(1-\delta)^2}\frac{4}{|1-\xi|^4}) d \xi \wedge d \bar \xi.
    \end{split}
    \end{equation*}
   Here $\alpha,\beta=1,...,n-1$.  Since $\delta\in (0,1)$ and $\xi \in \frac{3}{4}\Delta$, we have that $Re(-\frac{1+\delta}{1-\delta}\frac{1+\xi}{1-\xi})<0$. As a result, 
    \begin{equation*}
        |exp(-\frac{1+\delta}{1-\delta}\frac{1+\xi}{1-\xi}) (-\frac{1+\delta}{1-\delta}\frac{2}{(1-\xi)^2})|
    \end{equation*}
    is uniformly bounded independent of $\delta$. 
    
    Similarly, the derivatives of $\Phi_{\delta}^* \omega_X$ of any order are bounded with respect to the Euclidean metric. Recall that we have shown that $\Phi_{\delta}^* \omega_0$ is $C^{\infty}$ quasi-isometric to the Euclidean metric. We have that the derivatives of $\Phi_{\delta}^* \omega_X$ of any order are bounded with respect to $\Phi_{\delta}^* \omega_0$. This shows that:
    \begin{equation*}
        |\nabla^k_{\Phi_{\delta}^* \omega_0}\Phi_{\delta}^* \omega_X|_{\Phi_{\delta}^* \omega_0}\le C(k)
    \end{equation*}
    for any $k \ge 0$. Then we have that:
$
        |\nabla^k_{\omega_0}\omega_X|_{ \omega_0}\le C(k).
$
\end{proof}

\subsection{Function spaces}
\begin{defn}
If $U$ is a polydisc neighborhood of $D$ with $U \cap D$ given by $\{z_n=0\}$, we define for $f\in C^{p,\alpha}_{loc}(U \setminus D), (p,\alpha)\in \mathbb{N}\times [0,1),$
\begin{equation*}
||f||_{C^{p,\alpha}(U \setminus D)} \triangleq \sup_{\delta \in (0,1)}||\Phi_{\delta}^* f||_{C^{p,\alpha}(\mathcal{P})},
\end{equation*}
assuming that $U \subset \Delta^{n-1} \times (c\Delta).$ 

Then given a finite number of such open sets $U\in \mathcal{U}$, covering $D$ and an open set $V \subset \subset X \setminus D$ such that $X=V \cup \bigcup_{U \in \mathcal{U}} U$ and a partition of unity $\{\chi_V\} \cup \{\chi_U: U \in \mathcal{U}\},$ we can define the H\"older space 
\begin{equation*}
C^{p,\alpha}(M) \triangleq \{f\in C_{loc}^{p,\alpha}(M) : ||\chi_V f||_{C^{p,\alpha}(V)}+\max_{U\in 
 \mathcal{U}} ||\chi_U f||_{C^{p,\alpha}(U \setminus D)} < \infty \}.
\end{equation*}
\end{defn}
\begin{defn}
We can define the weighted H\"older norm:
\begin{equation*}
C_{\eta}^{k,\alpha} \triangleq \{f\in C^{k,\alpha}_{loc}(M) : ||\chi_V f||_{C^{p,\alpha}(V)}+\sup_{U\in \mathcal{U}} \sup_{\delta \in (0,1)}||(1-\delta)^{\eta}\Phi_{\delta}^* (\chi_U f)||_{C^{k,\alpha}(\mathcal{P})} < \infty\}.
\end{equation*}
Define
\begin{equation*}
C_{\eta,\mathbb{C}}^{k,\alpha}=\{v=f+ \sqrt{-1}g: f,g \in C_{\eta}^{k,\alpha}\}.
\end{equation*}
Since  $\frac{1}{C(1-\delta)} \le \Phi_{\delta}^* \rho \le \frac{C}{1-\delta}$ for some constant $C$, $||(1-\delta)^{\eta}\Phi_{\delta}^* (\chi_U f)||_{C^{k,\alpha}(\mathcal{P})} $ is equivalent to $||\Phi_{\delta}^* (\rho^{-\eta} \chi_U f)||_{C^{k,\alpha}(\mathcal{P})} $.
Heuristically, $f\in C_{\eta}^{k,\alpha}$  implies that $f =O(\rho^{\eta})$. We can also define:
\begin{equation*}
    C_{\eta}^{\infty}= \cap_{k=0}^{\infty} C_{\eta}^{k,\alpha}.
\end{equation*}
\end{defn}
\begin{defn}
 We can also define the weighted Sobolev space:
\begin{equation*}
W_{\eta}^{k,2} \triangleq \{v \in W^{k,2}_{loc}(M): \int_{M} \Sigma_{i=0}^k |\nabla_i v|^2 \rho^{-2\eta}\omega^n < \infty \}.
\end{equation*}
Define
\begin{equation*}
W_{\eta,\mathbb{C}}^{k,2} \triangleq \{v=f+ \sqrt{-1}g: f,g \in W_{\eta}^{k,2} \}.
\end{equation*}
\end{defn}

Clearly, $W_{\eta}^{k,2} \subset W_{\eta'}^{k,2} $, when $\eta\leq \eta'$.

\subsection{Poincar\'e type $C^{1,1}$ geodesic}
 Next, we talk about the setting for the Poincar\'e type $C^{1,1}$ geodesic. Consider the space $\mathfrak X=X\times R$, where $R$ is a cylinder $S^1 \times [0,1]$. Let $\pi$ be the projection from $\mathfrak X$ to $X$. Then the background metric on $\mathfrak X$ can be taken as 
 \begin{align*}
 \omega^\ast \triangleq \pi^* \omega +\sqrt{-1}dz^{n+1} \wedge d\bar z^{n+1},\quad   \omega^\ast_X \triangleq \pi^* \omega_X +\sqrt{-1}dz^{n+1} \wedge d\bar z^{n+1}.
 \end{align*} Here $(z^{n+1})$ is the standard coordinate of the cylinder and we write $$z^{n+1}=t+\sqrt{-1}s.$$  
 Clearly,  we have
  \begin{align*}
 \omega^\ast =   \omega^\ast_X -Ai \partial \bar \partial \pi^*\mathbf{u},\quad \pi^*\mathbf{u}
 = \log[\lambda-\log(|\sigma|^2)].
 \end{align*} Here, $\sigma$ is a section of $\mathfrak D=D\times R$.
 
 S. Semmes \cite{Se} observed that the geodesic can be seen as a $S^1$ invariant function on $\mathfrak X$. We will use this perspective. 
 We denote $\Psi=\varphi-|z^{n+1}|^2$.
 The geodesic connecting $\varphi_0,\varphi_1$ satisfies a degenerate Monge-Amp\'ere equation with Poincar\'e singularity
 \begin{align*}
 (\omega^\ast+dd^c\Psi)^{n+1}=\frac{n+1}{4}(\ddot{\varphi}-|\partial \dot{\varphi}|^2_{\omega_{\varphi}})\cdot\omega_{\varphi}^n \wedge \sqrt{-1}d z^{n+1} \wedge d \bar z^{n+1}=0\text{ in }\mathfrak M=M\times R
 \end{align*}
 with the boundary condition $\Psi=\Psi_0$ on $X\times\partial R$, where we define
 \begin{align*}
 \Psi_0=\varphi_0-s^2\text{ on }X\times \{0\}\times S^1,\quad
\Psi_0=\varphi_1-1-s^2\text{ on }X\times \{1\}\times S^1,
\end{align*} 
where $d,\partial$ and $\bar \partial$ are those of $M$ and the dot $\dot{}$ stands for $\partial_t$. 
 We also set $$\tilde \Psi_0\triangleq (1-t)\varphi_0+t\varphi_1$$ and define $\Psi_1$ to be $\tilde\Psi_0$ plus a sufficiently convex function on $z^{n+1}$, which vanishes on $X \times \partial R$.

 Auvray proved in the Theorem 2.1 and the Corollary 2.2 of \cite{A} the existence of the Poincar\'e type $\epsilon$-geodesic:
\begin{lem}\label{eps geodesic}
    For any $\varphi_0,\varphi_1\in \widetilde{\mathcal{PM}}_{\Omega}$ and any small enough $\epsilon>0$, there exists a path $\varphi^{\epsilon}$, denoted as $\epsilon$-geodesic, from $\varphi_0$ to $\varphi_1$, satisfying the equation of $\Psi^\epsilon=\varphi^\epsilon-|z^{n+1}|^2$
    \begin{equation*}
        (\om^\ast+dd^c\Psi^\epsilon)^{n+1}=\frac{n+1}{4}(\ddot{\varphi}^{\epsilon}-|\partial \dot{\varphi}^{\epsilon}|^2_{\omega_{\varphi^{\epsilon}}})\cdot\omega_{\varphi^{\epsilon}}^n \wedge  \sqrt{-1} dz^{n+1} \wedge d \bar z^{n+1}
        =\epsilon \cdot (\om^\ast+dd^c\Psi_1)^{n+1}.
    \end{equation*}
    There exists $C>0$ such that for all $\epsilon$,
    \begin{equation*}
        |\varphi^{\epsilon}-\Psi_0|,\quad |d \varphi^{\epsilon}|_{\omega}, \quad|\ddot{\varphi}^{\epsilon}|, \quad|d\dot{\varphi}^{\epsilon}|_{\omega},\quad |i\partial \bar \partial \varphi^{\epsilon}|_{\omega}\le C.
    \end{equation*}
     Moreover, we have that:
    \begin{equation*}
        \varphi^{\epsilon}- \Psi_0 \in C^{\infty}.
    \end{equation*}
\end{lem}
Then the Poincar\'e type $C^{1,1}$ geodesic is the limit of  $\epsilon$-geodesics:
\begin{lem}\label{C11 geodesic}
    For any $\varphi_0,\varphi_1\in \widetilde{\mathcal{PM}}_{\Omega}$, there exists a geodesic $\varphi$ such that there exists a constant $C>0$ such that:
    \begin{equation*}
        |\varphi- \Psi_0|,\quad  |d\varphi|_{\omega}, \quad|\ddot{\varphi}|,\quad |d \dot{\varphi}|_{\omega},\quad |i\partial \bar \partial \varphi|_{\omega}\le C.
    \end{equation*}
    and for any compact set $K\subset \subset M\times(0,1)$ and any constant $\alpha\in (0,1)$, we have that 
    \begin{equation*}
        \lim_{\epsilon \rightarrow 0} |\varphi^{\epsilon}-\varphi|_{C^{1,\alpha}(K)}=0.
    \end{equation*}
\end{lem}

\subsection{Energy functionals}
Next we define several functionals defined on $\widetilde{\mathcal{PM}}_{\Omega}$:
\begin{equation}\label{mathcal E}
    \mathcal{E}(\varphi)\triangleq \int_X \varphi  \Sigma_{j=0}^n \omega_\varphi^{n-j}\wedge \omega^j.
\end{equation}
Given a closed $(1,1)$-form (or current) $T$ bounded by a Poincar\'e type K\"ahler metric of any order, we set 
\begin{equation*}
    \mathcal{E}^T (\varphi )\triangleq \int_X \varphi  \Sigma_{j=0}^{n-1}\omega_\varphi ^{n-j-1}\wedge \omega^j \wedge T.
\end{equation*}
 Denote $\mu_0=\omega^n$. For any measure $\mu$ which is absolutely continuous with respect to $\mu_0$, we can also define the entropy term:
\begin{equation*}
    H_{\mu_0}(\mu)\triangleq \int_X log(\frac{d\mu}{d\mu_0})d\mu
\end{equation*}
The $K$-energy can be expressed as
\begin{equation}\label{decom kenergy}
    \mathcal{M}(\varphi )\triangleq \frac{\bar R}{n+1}\mathcal{E}(\varphi )-\mathcal{E}^{Ric_{\omega}}(\varphi ) +H_{\mu_0}(\omega_\varphi ^n).
\end{equation}
We can define the $J_{\chi}$ functional as follows:
\begin{equation*}
J_{\chi}(\varphi)=\frac{1}{n!} \int_X \varphi \Sigma_{k=0}^{n-1} \chi \wedge \omega_0^k \wedge \omega_{\varphi}^{n-1-k}- \frac{1}{(n+1)!} \int_M \underline{\chi} \varphi \Sigma_{k=0}^n \omega_0^k \wedge \omega_{\varphi}^{n-k}.
\end{equation*}
Here $$\underline{\chi}= \frac{\int_X \chi \wedge \frac{\omega_0^{n-1}}{(n-1)!}}{\int_X \frac{\omega_0^n}{n!}}.$$
The following Gaffney's Stokes theorem \cite{MR0062490}
 is pivotal when we compute the derivative of the functionals we defined above:
\begin{lem}\label{G-S}
Let $(X,g)$ be a complete n-dimensional Riemannian manifold where $g$ is a $C^2$ metric tensor. Let $\Theta$ be a $C^1$ $(n-1)$ form on $M$  such that both $|\Theta|_g$ and $|d\Theta|_g$ are in $L^1(X,g)$. Then we have that $\int_X d\Theta = 0$.
\end{lem}
With this Lemma, we can get the following lemma. The details was shown in our previous paper \cite{XZ}
\begin{lem}\label{functional first order}
    Suppose that $\varphi \in \widetilde{\mathcal{PM}}_{\Omega}$ and $v= O(\mathbf{u})$ with the derivatives of any order bounded with respect to a Poincar\'e type K\"ahler metric. Let $\mathcal{E}$, $\mathcal{E}^T$ and $H_{\mu_0}(\omega_\varphi^n)$ be defined as before. Then we have that:
    \begin{equation*}
        d \mathcal{E}|_\varphi (v)=(n+1)\int_X v \omega_\varphi^n,\quad 
        d \mathcal{E}^T |_\varphi (v) = n \int_X v \omega_\varphi^{n-1}\wedge T.
    \end{equation*}
    and
    \begin{equation*}
        d H_{\omega^n}(\omega_\varphi^n)(v)= \int_X v(-R_\varphi +tr_{\omega_\varphi} Ric_{\omega}) \omega_\varphi^n.
    \end{equation*}
    Here $R_\varphi$ is the scalar curvature of $\omega_\varphi$.
\end{lem}

We can directly compute that 
\begin{equation*}
\frac{d}{dt} J_{\omega} (\varphi_t)= \int_M (tr_{\varphi_t}\omega - n) \dot{\varphi}_t \frac{\omega_{\varphi_t}^n}{n!},
\end{equation*}
and
\begin{equation*}
\frac{d^2}{dt^2} J_{\omega} (\varphi_t) = \int (\ddot{\varphi}- |\nabla \dot{\varphi}|_{\varphi_t}^2)(tr_{\varphi_t}\omega -n) \omega_{\varphi_t}^n +\int \dot{\varphi}_{,\alpha} \dot{\varphi}_{,\bar \beta}\omega_{\bar \alpha \beta} \omega_{\varphi_t}^n >0.
\end{equation*}
This implies that the functional $J_{\omega}$ is strictly convex along smooth Poincar\'e type geodesic. By approximating Poincar\'e type $C^{1,1}$ geodesics with Poincar\'e type $\epsilon-$geodesic as in \cite{XZ}, we can see that $J_{\omega}$ is also strictly convex along Poincar\'e type $C^{1,1}$ geodesics.

\subsection{Fiber bundle structure of a neighbourhood of D}

According to the Section 3 of  \cite{A2}, a neighbourhood of $D$, denoted as $\mathcal{N}_A$, can be seen as a $S^1$ bundle over $[A,\infty) \times D$. This fiber bundle can be written as $$q: \mathcal{N}_A \setminus D \xrightarrow{q=(t,p)} [A,\infty) \times D.$$   The function $t$ is defined in \cite{A2}. We have that $t=\mathbf{u}$ up to a perturbation which is a $O(e^{-t})$, that is, a $O(\frac{1}{|log |\sigma||})$, as well as its derivatives of any order with respect to Poincar\'e type K\"ahler metrics. Denote $p$ as the projection from $\mathcal{N}_A \setminus D$ to $D$. We can also define a connection $\widetilde{\eta}$ in $\mathcal{N}_A \setminus D$ which can be seen as a volume form on each $S^1$ fibre such that $$J dt =2e^{-t} \widetilde{\eta} +O(e^{-t}).$$
In a cusp coordinate $(z_1,...,z_n=r e^{i\theta})$, one has 
\begin{equation}\label{eta dtheta}
\widetilde{\eta} =d \theta +O(1)
\end{equation}
in the sense that $\widetilde{\eta}- d\theta$ and all the derivatives of it of any order with respect to $\omega$ is bounded.  Then we can express Poincar\'e type K\"ahler metrics using $t$, $\widetilde{\eta}$ and $\theta$ as follows, according to the Proposition 1.2 of \cite{A}:
\begin{equation}\label{t coordinate}
p^* g_D + \frac{|dz^n|^2}{2|z^n|^2 \log^2 (|z^n|)} = p^* g_D + dt^2 + 4e^{-2t} \widetilde{\eta}^2 +O(e^{-t}) = p^* g_D + dt^2 + 4e^{-2t} d \theta^2 + O(e^{-t}),
\end{equation}
and
\begin{equation}\label{t coordinate 2}
p^* \omega_D + \frac{\sqrt{-1}dz^n \wedge d \bar z^n}{2|z^n|^2 \log^2 |z^n|}= p^* \omega_D -2e^{-t} dt \wedge d\theta +O(e^{-t}) = p^* \omega_D + dd^c (-t) +O(e^{-t}).
\end{equation}

Given an arbitrary function $f$ supported in a neighbourhood $\mathcal{N}_A$ of $D$, we can decompose $f$ as:
\begin{equation}\label{f decom}
f= f_0(t,p) + f^{\bot},
\end{equation}  
where $$f_0(t,p)=\frac{1}{2\pi}\int_{q^{-1}(t,p)} f \widetilde{\eta}$$ is the $S^1$ invariant part and $f^{\bot}$ is the part that is perpendicular to $S^1$ invariant functions. For any $S^1$ invariant function $u$, we have that
\begin{equation}\label{ddc}
dd^c u = 2(u_t - u_{tt})e^{-t} dt \wedge \widetilde{\eta} -2e^{-t} d_D u_t \wedge \widetilde{\eta} - dt \wedge d_D^c u_t + d d_D^c u +O(e^{-t}),
\end{equation}
where $d_D$ and $d_D^c$ are differential operators on $D$, according to the section 3 of \cite{A2}. Note that the definition of $dd^c$ in our case differs from the definition in \cite{A2} by a sign.
\subsection{Holomorphic vector fields}

\begin{defn}
We define the following things:
\begin{enumerate}
\item Define $\mathbf{h}_{\parallelsum}^D$ as the set of holomorphic vector fields on $X$ that are parallel to the divisor $D$.
\item Define $\mathbf{h}_{\parallelsum,\mathbb{C}}^D=\{V\in \mathbf{h}_{\parallelsum}^D: V=\nabla^{1,0}f \text{ for some complex valued function f}\}$.
\item Define $\mathbf{a}^D_{\parallelsum}(M)$  as the Lie subalgebra of $\mathbf{h}^D_{\parallelsum}$ consisting of the autoparallel, holomorphic vector fields of $M$ in $\mathbf{h}^D_{\parallelsum}$. 
\item Define $\mathbf{h}_{\parallelsum,\mathbb{R}}^D=\{V\in \mathbf{h}_{\parallelsum}^D: V=\nabla^{1,0}f \text{ for some real valued function f}\}$. 
\item Define $\mathbf{h}^D$ as the set of holomorphic vector fields on $D$.
\item Denote $Aut_0^D(X)$ as the identity component of the set of biholomorphisms on $M$ that preserve $D$.
\item Denote $Iso^D(X,\omega)$ as the set of biholomorphisms of $X$ preserving $D$ and preserving $\omega$.
\item Denote $Iso_0^D(X,\omega)$ as the identity component of $Iso^D(X,\omega)$.
\item Denote $Iso(D,\omega_D)= \{g \in Aut(D): g^* \omega_D = \omega_D\}$. 
\item 
Define the Mabuchi distance on $\widetilde{\mathcal{PM}}_{\Omega}$ as follows: for any two K\"ahler potentials $\varphi_0 , \varphi_1 \in \widetilde{\mathcal{PM}}_{\Omega}$, let  $\varphi_t$ be the Poincar\'e type $C^{1,1}$ geodesic connecting them given by the Lemma \ref{C11 geodesic}. Denote $\omega_t =\omega +dd^c \varphi_t$. Denote $b_t$ as the average of $\dot{\varphi_t}$ with respect to $\omega_t^n$. Then the Mabuchi distance is:
\begin{equation*}
d(\omega_1, \omega_0)^2 =\int_0^1 dt \int_X |\dot{\varphi_t}-b_t|^2 \omega_t^n.
\end{equation*}
\item For any $K \subset Aut_0^D(X)$, we can  define $\mathbf{h}^D_{\parallelsum, \mathbb{C}, K}= \{V\in \mathbf{h}^D_{\parallelsum, \mathbb{C}}: \text{ the flow of }ImV \text{ lies in }K\}$.
\item Given a vector field $V\in \mathbf{h}^D_{\parallelsum}$. we can define $\mathcal{PM}_{\Omega,V}\triangleq \{\omega \in \mathcal{PM}_{\Omega}: \omega \text{ is invariant under } Im V\}$.
\end{enumerate}
\end{defn}

\section{Solvability of Lichnerowicz operator}
In our previous paper \cite{XZ}, we solved the Lichnerowicz operator. Now we improve the method used in \cite{XZ}. In fact, in \cite{XZ}, we mainly used Sobolev space to improve the decay rate of a solution to the Lichnerowicz operator. However, we assume that the right-hand side of the equation in a H\"older space instead of a Sobolev space. The transition between H\"older space and a Sobolev space can cause a loss of the decay rate. This is why in \cite{XZ}, the decay rate of the solution we got is slightly weaker than the decay rate of the right-hand side of the equation. 

  Throughout this section, we assume that the Poincar\'e type K\"ahler metric is asymptotic to a product metric in the sense of (\ref{e 5.1}). The main proposition we want to prove in this section is the Proposition \ref{solve l operator1} which is as follows:
\begin{prop}\label{solve l operator}
Suppose that $\omega$ is a Poincar\'e type K\"ahler metric satisfying (\ref{e 5.1}). Then there exists a constant $0< \delta_1< \frac{1}{2}$. For any $\eta_0 \in (0, \delta_1)$, for any $f\in C_{-\eta_0}^{1,\alpha}$  such that $\int_{M\setminus D}f u \omega^n=0$ for any $u\in \overline{\mathbf{h}^D_{\parallelsum,\mathbb{R}}}$, we can find a function $v\in C_{-\eta_0}^{5,\alpha} \oplus \chi(t) p^* Ker Re L_D$ such that $Re Lv =f$. 
\end{prop}

\subsection{Kernel and range}
To begin with, we need the following Lemma which characterizes the image of operators with closed range (See the Theorem 2.19 in \cite{B}).
\begin{lem}\label{lem 4.2}
Let $A: D(A) \subset E \rightarrow F$ be an unbounded linear operator that is densely defined and closed. The following properties are equivalent:
\begin{enumerate}
\item  $Im(A)$ is closed,
\item $Im(A^*)$ is closed,
\item $Im(A)=Ker(A^*)^{\bot}$,
\item $Im(A^*)=Ker(A)^{\bot}$.
\end{enumerate}
\end{lem}
In the above Lemma, we denote the range of $A$ as $Im(A)$ and we denote the kernal of $A$ as $Ker(A)$. 
In our case, we set $E=F=L_{\eta}^2$ and $A=Re L$. We define $Re L^*$ as:
if $u\in D(Re L^*)$, then for any $v\in D(Re L)$, $\int_M v Re L^* u \omega^n=\int_M Re L v  u \omega^n$. Since $L$ is self-adjoint, so is $\bar L$. Thus $Re L$ is self-adjoint.
Then we have that $D(Re L)=W_{\eta}^{4,2}=\{u: \Sigma_{k=0}^4 (\int_X |\nabla^k u|^2 e^{-2\eta t} \omega^n)^{\frac{1}{2}}< \infty\}$. $Re L^*=Re L|_{W_{-\eta}^{4,2}}$. 

Now we want to show that $Im(Re L)$ is closed. We need to use the Lemma below:
\begin{lem}\label{lem 4.3}
Suppose that $Re L$ satisfies the following formula for any $v \in W_{\delta}^{m,2}$ :
\begin{equation}\label{estimate for closed range}
||v||_{W_{\delta}^{m,2}(X \setminus D)}\le C(||Re L v||_{W_{\delta}^{m-4,2}(X\setminus D)}+||v||_{L^2(K)})
\end{equation}
for some compact set $K \subset \subset X \setminus D$ and some constant $\delta$ and $m \ge 4$. Then we have that $dim Ker(Re L|_{W_{\delta}^{m,2}})< \infty$ and  $Im(Re L|_{W_{\delta}^{m,2}})$ is closed.
\end{lem}
\begin{proof}
The proof of this Lemma is similar to the proof of the Lemma 5.3 in our previous paper \cite{XZ} with $L$ replaced by $Re L$ in our case.
\end{proof}

Using (\ref{e 5.1}), we have that:
\begin{equation}\label{e 5.3}
Ric_{\omega} = \frac{- \sqrt{-1} dz^n \wedge  d \bar z^n}{2|z^n|^2 \log^2 (|z^n|)} + p^* Ric_{\omega_D} +O(e^{-\eta t}).
\end{equation}
Then we have that:
\begin{equation}\label{R dec}
\begin{split}
R_{\omega}&=2n \frac{Ric_{\omega}\wedge \omega^{n-1}}{\omega^n}\\
&=2n\frac{(n-1)p^* Ric_{\omega_D}\wedge p^* \omega_D^{n-2}\wedge (\frac{\sqrt{-1} adz^n \wedge d \bar z^n}{2 |z^n|^2 \log^2 (|z^n|)})+(\frac{- \sqrt{-1} dz^n \wedge d \bar z^n}{2|z^n|^2 \log^2 (|z^n|)})\wedge p^* \omega_D^{n-1}+O(e^{-\eta t})}{n p^* \omega_D^{n-1}\wedge(\frac{\sqrt{-1} a dz^n \wedge  d \bar z^n}{2 |z^n|^2 \log^2 (|z^n|)})+O(e^{-\eta t})}\\
& =p^* R_{\widetilde{\omega}_j}-\frac{2}{a} +O(e^{-\eta t}),
\end{split}
\end{equation}
which gives that 
\begin{equation}\label{asymptotic S}
<\uparrow \bar \partial \varphi, \partial S>_{\omega}=<\uparrow \bar \partial \varphi,\partial S_{\omega_D}>_{\omega_D}+O(e^{-\eta t}).
\end{equation}

Next, we restrict $Re L$ on the space of $S^1$ invariant functions and consider $\frac{1}{2}\Pi_0 \circ (L+ \bar L) \circ q^* $. 

Recall that $$q: \mathcal{N}_{A} \setminus D \xrightarrow{q=(t,p)} [A,\infty) \times D.$$ So $q^*$ means canonically map a function defined on $[A,\infty) \times D$ to a function defined on $\mathcal{N}_{A} \setminus D$ which is invariant along each $S^1$ fiber. $\Pi_0$ is the map of a function to its $S^1$-invariant part. Using the (\ref{asymptotic S}) and the asymptotic behavior of $\omega$, we can see that  $\frac{1}{2}\Pi_0 \circ (L +\bar L)\circ q^* $ is asymptotic to the following operator (see \cite[Proposition 3.4]{A3}):
\begin{equation*}
Re L^0 \triangleq \frac{1}{2}(\frac{\partial}{\partial t}-\frac{\partial^2}{\partial t^2})^2 +(\frac{\partial}{\partial t}-\frac{\partial^2}{\partial t^2}) + \frac{1}{2}(L_{p^* \omega_D}+\bar L_{p^* \omega_D})+\Delta_{\omega_D}\circ (\frac{\partial}{\partial t}-\frac{\partial^2}{\partial t^2}).
\end{equation*}
In fact, we can prove the following Lemma:

\begin{lem}\label{AA0}
Suppose that $\omega$ is a Poincar\'e type K\"ahler metric satisfying (\ref{e 5.1}). Then for any $\delta \in \mathbb{R}$, there exists a constant $C$ such that for any $S^1-$invariant function $v$, we have that:
\begin{equation*}
||Re Lv-Re L^0 v||_{C_{\delta-\eta}^{k,\alpha}}  \le C ||v||_{C^{k+4,\alpha}_{\delta}}
\end{equation*}
and
\begin{equation*}
||(Re L-Re L^0) v||_{W_{\delta -\eta}^{k,2}}\le C ||v||_{W_{\delta}^{k+4,2}}.
\end{equation*}
\end{lem}
\begin{proof}
Recall that (\ref{ddc}) implies that
\begin{equation*}
dd^c v = 2(v_t - v_{tt})e^{-t} dt \wedge \widetilde{\eta} -2e^{-t} d_D v_t \wedge \widetilde{\eta} - dt \wedge d_D^c v_t + d d_D^c v +O(e^{-t}),
\end{equation*}
Using (\ref{e 5.1}), we have that
\begin{equation*}
\Delta_{\omega} v =(\partial_t -\partial_t^2)v +p^* \Delta_{\omega_D} v +O(e^{-\eta t}(|\nabla^2 v|+|\nabla v|)).
\end{equation*} 
and
\begin{equation}\label{delta2 asy}
\Delta_{\omega}^2 v =(\partial_t -\partial_t^2)^2 v +(p^*\Delta_{\omega_D})^2 v +(\partial_t -\partial_t^2)p^* \Delta_{\omega_D} v +p^* \Delta_{\omega_D}(\partial_t -\partial_t^2)v +O(e^{-\eta t}(|\nabla^2 v|+|\nabla^3 v|+|\nabla^4 v|)).
\end{equation}
Using (\ref{e 5.1}) again, we have that the Ricci form of $\omega$ has an asymptotic behavior:
\begin{equation}\label{ric asym}
Ric_{\omega}=dt \wedge 2e^{-t}\eta +p^* Ric_{\omega_D} +O(e^{-\eta t})
\end{equation}
Note that the real part of Lichnerowicz operator can be expressed as:
\begin{equation}\label{L dec}
Re Lv=\frac{1}{2} \Delta_{\omega}^2 v +<Ric_{\omega},dd^c v>_{\omega} +1/2(v^{\alpha}R_{\alpha} + v_{\alpha} R^{\alpha}).
\end{equation}
We define an operator:
\begin{equation}
Re L^0 =\frac{1}{2}(\partial_t-\partial_t^2)^2 +(\partial_t -\partial_t^2)+ ReL_D +\Delta_{\omega_D}\circ (\partial_t -\partial_t^2).
\end{equation}
Using the Formulae (\ref{delta2 asy}), (\ref{ric asym}), (\ref{L dec}) and (\ref{asymptotic S}), we have that:
\begin{equation}\label{AA0 holder difference}
||Re Lv-Re L^0 v||_{C_{\delta-\eta}^{k,\alpha}}  \le C ||v||_{C^{k+4,\alpha}_{\delta}}
\end{equation}
and
\begin{equation*}
||(Re L-Re L^0) v||_{W_{\delta -\eta}^{k,2}}\le C ||v||_{W_{\delta}^{4,2}}.
\end{equation*}
\end{proof}

Then we have that:
\begin{prop}\label{fredholm}
\begin{equation*}
Im(Re L |_{W_0^{4,2}})=Ker (Re L|_{W_0^{4,2}})^{\bot}.
\end{equation*}
\end{prop}
\begin{proof}
Note that (\ref{estimate for closed range}) is proved by the Proposition 3.2 of \cite{S} for any $\delta$ which is not an indicial root for $Re L$. Using the Lemma \ref{l0 iso}, we have that $\delta=0$ is not an indicial root. Then, we can use the Lemma \ref{lem 4.2} and the Lemma \ref{lem 4.3} to conclude the proof of the proposition.
\end{proof}
\subsection{Kernel and holomorphic vector fields.}
 Note that $\overline{\mathbf{h}^D_{\parallelsum, \mathbb{R}}}=\{f\in C_{\mathbb{R}}^{\infty}(M \setminus D): \nabla^{1,0}f \in \mathbf{h}_{\parallelsum}^D\}$.  We record the following Lemma:
\begin{lem}\label{ker a}
Suppose that $\omega$ is a Poincar\'e type K\"ahler metric, then
\begin{equation*}
Ker(Re L|_{W_{0}^{k,2}})=\overline{h^D_{\parallelsum,\mathbb{R}}},
\end{equation*}
for $k \ge 4$.
\end{lem}
\begin{proof}
The proof of this Lemma is similar to the proof of the Lemma 5.6 in \cite{XZ} with $L$ replaced by $Re L$.
 \end{proof}
We also record the following Lemma which we will use the a section below.
 \begin{lem}\label{ker l}
Suppose that $\omega$ is a Poincar\'e type K\"ahler metric, then
\begin{equation*}
Ker(L |_{W_{0,\mathbb{C}}^{k,2}})=\overline{h^D_{\parallelsum,\mathbb{C}}},
\end{equation*}
for $k \ge 4$.
 \end{lem}
 \begin{proof}
 The proof is similar to the proof of the last lemma.  The formula $Ker(L |_{W_{0,\mathbb{C}}^{k,2}}) \subset \overline{h^D_{\parallelsum,\mathbb{C}}}$ can be shown by using the local Taylor expansion of holomorphic functions near the divisor.  Indeed, for any $u\in Ker(L|_{W_{0,\mathbb{C}}^{k,2}})$, we have that 
\begin{equation*}
0=\int_M L u \bar u \omega^n = \int_M \mathcal{D}^* \mathcal{D}u \bar u \omega^n =  \int_M |\mathcal{D}u|^2 \omega^n .
\end{equation*}

This implies that $\mathcal{D}u=0$ which means that $V \triangleq \nabla^{1,0}u$ is a holomorphic vector field on $X\setminus D$. We should be careful that we don't know if $V$ is a holomorphic vector field on $X$ or not. We will prove that $V$ can be extended to $D$. Since $u\in W_{0,\mathbb{C}}^{k,2}$, we can get that $|V| \in L^2(\omega^n)$. In an arbitrary cusp coordinate domain $U$, we denote $V=v^i \frac{\partial}{\partial z^i}$. Since $\omega$ is equivalent to the standard cusp metric (\ref{standard cusp}), we have that:
\begin{equation*}
\int_U |V|_{\omega_0}^2 \omega_0^n \le C \int_M |V|_{\omega}^2 \omega^n < +\infty
\end{equation*}
Then we have that:
\begin{equation*}
\begin{split}
\int_U |V|_{\omega_0}^2 \omega_0^n &\ge C \int_U |v^n|^2 (\omega_0)_{n \bar n} \omega_0^n =\int_U |v^n|^2 \frac{1}{|z_n|^2 log^2 (|z_n|^2)}\frac{n!}{|z_n|^2 log^2 (|z_n|^2)} dVol_E \\
&= \int_U |v^n|^2 \frac{n!}{|z_n|^4 log^4(|z_n|^2)}dVol_E.
\end{split}
\end{equation*}
Here $dVol_E=(\sqrt{-1})^ndz^1 \wedge d \bar z^1 \wedge... \wedge dz^n \wedge d\bar z ^n$. Note that we have a Laurent series of $v^n$(c.f. the proposition 1.4 of \cite{R})
 \begin{equation*}
 v^n =\Sigma_{\mu \in \mathbb{N}^{n-1}} \Sigma_{k\in \mathbb{Z}}C_{\mu k} z'^{\mu}z_n^k.
 \end{equation*}
Here $z'=(z_1,...,z_{n-1}).$  Let $\epsilon>0$ be a constant such that $U(\epsilon)\triangleq \{z: |z_i|\le \epsilon \text{ for any }i\} \subset U$. We also denote $U'(\epsilon)\triangleq \{z: |z_i|\le \epsilon \text{ for any }i\le n-1\}$. Then we have that:
\begin{equation*}
\begin{split}
&\int_{U_{\epsilon}}|v^n|^2 \frac{n!}{|z_n|^4 log^4 (|z_n|^2)} dVol_E \\
&=\int_{U'_{\epsilon}} d Vol_E (z')\int \Sigma_{\mu \in \mathbb{N}^{n-1}}\Sigma_{k\in \mathbb{Z}}|C_{\mu k}|^2 |z'|^{2\mu} \frac{|z_n|^{2k-4}}{log^4(|z_n|^2)}  \sqrt{-1} dz^n \wedge d\bar z^n.
\end{split}
\end{equation*}
Combining the above Formulae, we have that $C_{\mu k}=0$ for any $k\le 0$. This proves that  $v^n$ can be extended holomorphically to $D$ and vanishes on $D$. Similarly, we can show that $v^i$ can be extended for any $i\le n-1$. This concludes the proof of $Ker(L |_{W_{0,\mathbb{C}}^{k,2}}) \subset \overline{h^D_{\parallelsum,\mathbb{C}}}.$\\

The formula $\overline{h^D_{\parallelsum,\mathbb{C}}} \subset Ker(L|_{W_{0,\mathbb{C}}^{k,2}}) $ can be shown as follows:
For any $f\in \overline{h^D_{\parallelsum,\mathbb{C}}}$, we have that $V \triangleq \nabla^{1,0}_{\omega}f$
 is a holomorphic vector field on $X$.  First, we claim that:
 \begin{equation}\label{nabla diff metric}
 V =\nabla_{\omega_X}^{1,0}(f-V(u)),
 \end{equation}
 where $u$ is the potential such that $\omega =\omega_X +dd^c u$. Indeed, we can calculate that:
 \begin{equation*}
v^i ((g_0)_{i\bar l}+u_{i\bar l})= v^i g_{i\bar l}=g^{i\bar j}f_{\bar j}g_{i\bar l}=f_{\bar l}.
 \end{equation*}
 Multiply $g_0^{k \bar l}$ on the both sides of the above formula and take the sum with respect to $l$. We get:
 \begin{equation*}
 g^{k\bar j} f_{\bar j} \frac{\partial}{\partial_z^k} +\nabla_{\omega_X}^{1,0} (V(u))=\nabla_{\omega_X}^{1,0}f.
 \end{equation*}
 This concludes the proof of the claim.  By the definition of the Poincar\'e type metric, we have that the derivative of $u$ are bounded with respect to a given poincar\'e type metric. As a result, we have that $V(u)\in C_{0,\bC}^{k,\alpha}$ for any $k$.  Using (\ref{nabla diff metric}) and the fact that $\omega_X$ is a smooth K\"ahler metric, we get that  $f-V(u)$ is a smooth function on $M$. Then we have that $f\in C_{0,\bC}^{k,\alpha}\subset W_{0,\mathbb{C}}^{4,2}$. This concludes the proof of the Lemma. 
 \end{proof}

 \subsection{\texorpdfstring{$u_0$}{u0} and \texorpdfstring{$u^{\bot}$}{u}} Recall that we defined $u_0$ and $u^{\bot}$ in (\ref{f decom}). We have the following technical lemmas:
\begin{lem}\label{u dec}
For any $\delta,$ there exists a uniform constant $C$ such that $||u^{\bot}||_{C_{\delta}^{k,\alpha}}\le C ||u||_{C^{k,\alpha}_{\delta}}$ and $||u_0||_{C_{\delta}^{k,\alpha}}\le C ||u||_{C^{k,\alpha}_{\delta}}$ for any $u$ such that $||u||_{C^{k,\alpha}_{\delta}}< \infty$.
\end{lem}
\begin{lem}\label{u dec sobolev}
For any $\delta,$ there exists a uniform constant $C$ such that $||u^{\bot}||_{W_{\delta}^{k,2}}\le C ||u||_{W^{k,2}_{\delta}}$ and $||u_0||_{W_{\delta}^{k,2}}\le C ||u||_{W^{k,2}_{\delta}}$ for any $u$ such that $||u||_{W^{k,2}_{\delta}}< \infty$.
\end{lem}
Recall that $t$ is a function defined in section 3.7. Since the integral of $u^{\bot}$ on each $S^1$ fiber is zero and the length of the $S^1$ fiber exponentially decay to zero as $t$ goes to $\infty$, we have the following Lemma basically saying that the decay rate of $u^{\bot}$ con be improved if we have control on its higher order derivatives. See the section 3 of \cite{A2} and the Formula (3.6) of \cite{S}.

\begin{lem}\label{improve decay}
For any $\delta \in \mathbb{R}$ and $k\in \mathbb{N}$, there exists a constant $C$ such that:
\begin{equation*}
||u^{\bot}||_{W_{\delta}^{k,2}}\le C||u^{\bot}||_{W_{\delta+1}^{k+1,2}}. 
\end{equation*}
holds for any $k$ and $u$ such that $||u^{\bot}||_{W_{\delta+1}^{k+1,2}}< \infty$.
\end{lem}

\begin{lem}\label{improve decay holder}
For any $\delta \in \mathbb{R}$ and $k\in \mathbb{N}$, there exists a constant $C$ such that:
\begin{equation*}
||u^{\bot}||_{C_{\delta}^{k,\alpha}}\le C||u^{\bot}||_{C_{\delta+1}^{k+1,\alpha}}. 
\end{equation*}
holds for any $k$ and $u$ such that $||u^{\bot}||_{C_{\delta+1}^{k+1,\alpha}}< \infty$.
\end{lem}

\subsection{Operator $Re L^0$}
Denote 
\begin{equation*}
W_{0,\delta}^{k,2}=\{u: u|_{t=0}=0, u_t |_{t=0}=0, \int_{D\times [0,\infty)} \Sigma_{l=0}^k |\nabla^l u|^2 e^{(-2\delta-1)t} dt dVol_D <\infty \}.
\end{equation*}
 Denote 
 \begin{equation*}
 C_{0,\delta}^{k,\alpha}=\{u: u|_{t=0}=0, u_t|_{t=0}=0, \sup_{(x,t)\in D \times[0,+\infty)} e^{-\delta t}|\nabla_i u|\le +\infty \text{ for any }i \le k \}.
 \end{equation*}
 The following results were proved by Auvray (see the Lemma 3.8 in \cite{A}):
\begin{lem}\label{l0 iso}
\begin{enumerate}
\item For any $\delta \in(-\frac{1}{2},\frac{1}{2})$, 
\begin{equation*}
Re L^0 : W_{0,\delta}^{4,2} ([0,\infty)\times D) \rightarrow L_{\delta}^2 ([0,\infty)\times D)
\end{equation*}
is an isomorphism.
\item There exists $\delta_0>0$ such that 
\begin{equation*}
Re L^0 : W_{0,\delta}^{4,2} ([0,\infty)\times D) \oplus \chi p^* ker Re L_D \rightarrow L_{\delta}^2 ([0,\infty)\times D)
\end{equation*}
is an isomorphism for $\delta \in (-\frac{1}{2}-\delta_0, -\frac{1}{2}).$
\end{enumerate}
\end{lem}

\begin{lem}\label{l0 iso holder}
\begin{enumerate}
\item For any $\delta \in (0,1)$, we have that:
\begin{equation*}
Re L^0 : C_{0,\delta}^{k+4,\alpha} ([0,\infty) \times D)  \rightarrow C_{\delta}^{k,\alpha}([0,\infty) \times D)
\end{equation*}
is an isomorphism.
\item There exists $\delta_0>0$ such that  for all $\delta \in (-\delta_0,0)$ 
\begin{equation*}
Re L^0 : C_{0,\delta}^{k+4,\alpha} ([0,\infty) \times D) \oplus \chi p^* ker Re L_D  \rightarrow C_{\delta}^{k,\alpha}([0,\infty) \times D)
\end{equation*}
is an isomorphism.
\end{enumerate}
\end{lem}

 We can prove the following Lemma:

\begin{lem}\label{l0 adjoint}
For any $u \in W_{0,0}^{4,2}$, we have that:
\begin{equation*}
 \int_{D\times [0,\infty)} Re L^0 u  u e^{-t} dt dvol_D= \int |u_{tt}|^2 e^{-t}+\int |u_t|^2 e^{-t}+\int e^{-t} |\mathcal{D}u|_D^2 +\int |\nabla_D \dot{u}|^2 e^{-t}.
\end{equation*}
\end{lem}
\begin{proof}
The proof of this Lemma is similar to the proof of the Lemma 5.10 in \cite{XZ} with $L^0$ replaced by $Re L^0$.
\end{proof}

\subsection{Regularity results}
Sektnan proved the following regularity result in \cite{S}:
\begin{lem}\label{l regularity}
Suppose $u \in W_{\delta-\frac{1}{2}}^{2,0}$ and suppose that $Re L u\in C_{\delta}^{k-4,\alpha}$ in the sense of distributions for a weight $\delta$. Then $u\in C_{\delta}^{k,\alpha}.$ Moreover, there is a $C>0$ such that:
\begin{equation*}
||u||_{C_{\delta}^{k+4,\alpha}} \le C(||Re L u||_{C_{\delta}^{k,\alpha}}+||u||_{W_{\delta-\frac{1}{2}}^{2,0}}).
\end{equation*}
\end{lem}

We also need the following regularity lemma:
\begin{lem}\label{lem 4.12}
Suppose that $u\in W_{\delta}^{0,2}$ and $Re L u \in W_{\delta}^{k,2}$. Then we have that $u\in W_{\delta}^{k+4,2}$ and
\begin{equation*}
||u||_{W^{k+4,2}_{\delta}}\le C(||Re Lu||_{W^{k,2}_{\delta}} +||u||_{W^{0,2}_{\delta}}).
\end{equation*}
\end{lem}
\begin{proof}
This lemma can be proved in the same way as the Lemma 1.12 of \cite{A}. We just sketch the proof here. We can use a covering of $X\setminus D$ using the quasi-conformal coordinate mentioned in section 3. In each coordinate, the Poincar\'e type K\"ahler metric is quasi-isometric to the Euclidean metric. Then, we can use the standard $L^p$ estimate for $Re L$ in each quasi-conformal coordinate since $Re L$ is a fourth-order elliptic operator and the coefficient of it is uniformly bounded in each quasi-conformal coordinate. Then, we patch them together to prove the lemma. 
\end{proof}

\subsection{Improve decay rate}
We want to prove some lemmas that help us improve the decay rate of $S^1$ invariant functions which are pivotal in our proof of the Proposition \ref{solve l operator}.
\begin{lem}\label{improve decay sobolev}
Suppose that $v$ is supported in a small neighbourhood of $D$ which can be seen as a $S^1$ bundle over $[0,+\infty) \times D$ as in the section 3.7. Suppose that $v$ is $S^1$ invariant.  Suppose that $Re Lv \in L_{\delta_1}^{2}$ and $v\in W_{\delta_2}^{4,2}$ with $-\eta \ge \delta_2 -\eta \ge \delta_1 > -\frac{1}{2}$, where $\eta$ is given by (\ref{e 5.1}). Then we have that 
\begin{equation*}
v \in W^{4,2}_{\delta_2 -\eta}.
\end{equation*} 
\end{lem}
\begin{proof}
Using the Lemma \ref{AA0}, we have that 
\begin{equation*}
||(Re L-Re L^0)v||_{W^{0,2}_{\delta_2 -\eta}}\le C ||v||_{W^{4,2}_{\delta_2}}.
\end{equation*}
Since $\delta_2 -\eta \ge  \delta_1$ and $Re Lv \in L_{\delta_1}^{2}$, we have that:
\begin{equation*}
Re L^0 v\in W_{\delta_2 -\eta}^{0,2}
\end{equation*}
Using the Lemma \ref{l0 iso}, we can get a function $h \in W_{0,\delta_2 -\eta}^{4,2}$ such that $Re L^0 h= Re L^0 v$. It suffices to prove that $h=v$. We can apply the Lemma \ref{l0 adjoint} with $u$ replaced by $v-h$ to get that:
\begin{equation*}
\begin{split}
 \int_{D\times [0,\infty)} Re L^0 (v-h)  (v-h) e^{-t} dt dvol_D& = \int |(v-h)_{tt}|^2 e^{-t}+\int |(v-h)_t|^2 e^{-t} \\
 &+\int e^{-t} |\mathcal{D}(v-h)|_D^2 +\int |\nabla_D (v-h)_t|^2 e^{-t}.
\end{split}
\end{equation*}
Since $Re L^0 h= Re L^0 v$, we get that $(v-h)_t=0$. Since $v|_{t=0}= h|_{t=0}$, we can get that $v=h$. This finishes the proof of this Lemma.
\end{proof}

\begin{lem}\label{improve decay sobolev 2}
Suppose that $v$ is supported in a small neighbourhood of $D$ which can be seen as a $S^1$ bundle over $[0,+\infty) \times D$ as in the section 3.7. Suppose that $v$ is $S^1$ invariant. Suppose that $Re Lv \in L_{\delta_1}^{2}$ and $v\in W_{\delta_2}^{4,2}$ with $-\frac{1}{2}+\eta > \delta_2 > -\frac{1}{2} >  \delta_1 \ge \delta_2 -\eta$, where $\eta$ is given by (\ref{e 5.1}). Then we have that 
\begin{equation*}
v \in W^{4,2}_{\delta_1} \oplus \chi(t) p^* Ker Re L_D.
\end{equation*} 
\end{lem}
\begin{proof}
Using the Lemma \ref{AA0}, we have that 
\begin{equation*}
||(Re L-Re L^0)v||_{L^{2}_{\delta_2 -\eta}}\le C ||v||_{W^{4,2}_{\delta_2}}.
\end{equation*}
Since  $-\frac{1}{2}+\eta > \delta_2 > -\frac{1}{2} >  \delta_1 > \delta_2 -\eta$ and $Re Lv \in L_{\delta_1}^{2}$, we have that:
\begin{equation*}
Re L^0 v \in L_{\delta_1}^{2}. 
\end{equation*}
Using the Lemma \ref{l0 iso}, we can get a function $h \in W_{0,\delta_1}^{4,2} \oplus \chi(t) p^* Ker Re L_D$, such that $Re L^0 h= Re L^0 v$. It suffices to prove that $h=v$. We can apply the Lemma \ref{l0 adjoint} with $u$ replaced by $h-v$ to get that:
\begin{equation*}
\begin{split}
 \int_{D\times [0,\infty)} Re L^0 (v-h)  (v-h) e^{-t} dt dvol_D& = \int |(v-h)_{tt}|^2 e^{-t}+\int |(v-h)_t|^2 e^{-t} \\
 &+\int e^{-t} |\mathcal{D}(v-h)|_D^2 +\int |\nabla_D (v-h)_t|^2 e^{-t}.
\end{split}
\end{equation*}
Since $Re L^0 h= Re L^0 v$, we get that $(v-h)_t=0$. Since $v|_{t=0}= h|_{t=0}$, we can get that $v=h$. This finishes the proof of this Lemma.
\end{proof}

\begin{lem}\label{improve decay holder}
Suppose that $v$ is supported in a small neighbourhood of $D$ which can be seen as a $S^1$ bundle over $[0,+\infty) \times D$ as in the section 3.7. Suppose that $v$ is $S^1$ invariant. Suppose that $Re Lv \in C_{\delta_1}^{0,\alpha}$ and $v\in C_{\delta_2}^{0,\alpha}$ with $0 \ge \delta_2  \ge \delta_1 \ge  -\eta$, where $\eta$ is given by (\ref{e 5.1}). Then we have that 
\begin{equation*}
v \in C^{4,\alpha}_{\delta_1} \oplus \chi(t) p^* Ker Re L_D.
\end{equation*} 
\end{lem}
\begin{proof}
Using the Lemma \ref{AA0}, we have that 
\begin{equation*}
||(Re L-Re L^0)v||_{C^{0,\alpha}_{\delta_2 -\eta}}\le C ||v||_{C^{4,\alpha}_{\delta_2}}.
\end{equation*}
Since $0 \ge \delta_2  \ge \delta_1 > -\eta$ and $Re Lv \in C_{\delta_1}^{0,\alpha}$, we have that:
\begin{equation*}
Re L^0 v\in C_{\delta_1}^{0,\alpha}
\end{equation*}
Using the Lemma \ref{l0 iso holder}, we can get a function $h \in C_{\delta_1}^{4,\alpha} \oplus \chi(t) p^* Ker Re L_D$ such that $Re L^0 h= Re L^0 v$. It suffices to prove that $h=v$. We can apply the Lemma \ref{l0 adjoint} with $u$ replaced by $h-v$ to get that:
\begin{equation*}
\begin{split}
 \int_{D\times [0,\infty)} Re L^0 (v-h)  (v-h) e^{-t} dt dvol_D& = \int |(v-h)_{tt}|^2 e^{-t}+\int |(v-h)_t|^2 e^{-t} \\
 &+\int e^{-t} |\mathcal{D}(v-h)|_D^2 +\int |\nabla_D (v-h)_t|^2 e^{-t}.
\end{split}
\end{equation*}
Since $Re L^0 h= Re L^0 v$, we get that $(v-h)_t=0$. Since $v|_{t=0}= h|_{t=0}$, we can get that $v=h$. This finishes the proof of this Lemma.
\end{proof}

\subsection{Proof of the Proposition \ref{solve l operator}}

In this proof we replace $\eta$  and $\eta_0$ by $\min\{\eta_0,\eta\}$ and assume that $\eta_0=\eta$ without loss of generality. Here $\eta_0$ is the constant in the Proposition \ref{solve l operator} and $\eta$ is the constant in (\ref{e 5.1}). We can also assume that $\eta < \delta_0$, where $\delta_0$ is given by the Lemma \ref{l0 iso} and the Lemma \ref{l0 iso holder}. We first sketch the proof of this Proposition. Note that $Ker (Re L |_{W_0^{4,2}})=\overline{\mathbf{h}^D_{\parallelsum,\mathbb{R}}}$.
As a result, for any $f\in C_{-\eta}^{1,\alpha} \cap (\overline{\mathbf{h}^D_{\parallelsum,\mathbb{R}}})^{\bot}$ with some $\eta>0$, we have that $f\in W_0^{1,2} \subset W_0^{0,2}$. Then using the  proposition \ref{fredholm}, we can find $u\in W_0^{4,2}$ such that $Re Lu =f$. Then we can use the Lemma \ref{lem 4.12} to get that $u\in W_0^{5,2}$. 

Then, we will show that the decay rate of $u$ can be improved. The idea is as follows: We can localize the problem in a neighbourhood of $D$ and assume that $u$ is supported in this neighbourhood of $D$.  Then we can decompose $u=u_0+ u^{\bot}$, where $u_0$ is the $S^1$ invariant part and $u^{\bot}$ is perpendicular to $S^1$ invariant functions. We improve the decay rate of $u_1$ using the Lemma \ref{l0 iso} and the Lemma \ref{l0 iso holder}. $u^{\bot}$ has a good decay rate using the Lemma \ref{improve decay} and the Lemma \ref{improve decay holder}.

Next, we prove the above argument rigorously.  Using the argument above, we can find $u \in W_0^{5,2}$ such that $Re Lu=f$. Using the standard local elliptic estimates, we can show that $u\in C_{loc}^{5,\alpha}(X\setminus D)$. Then we can take a small neighbourhood of $D$, denoted as $V_1$ and let $\chi$ be a cut-off function supported in $V_1$ which is equal to $1$ in a smaller neighbourhood of $D$, denoted as $V_2$. Note that in the rest of the proof we only need to use the property of $f$ near $D$. Since $Re L (\chi u) \in C_{loc}^{1,\alpha}(X \setminus D)$ and is equal to $Re L u =f$ in $V_2$, we can replace $u$ by $\chi u$ and assume that $u$ is supported in $V_1$ which is a $S^1$ bundle over $[0,\infty) \times D$ as in the section 3.7. Using the Lemma \ref{u dec sobolev}, we have that:
\begin{equation*}
||u_0||_{W^{5,2}_{0}}\le C ||u||_{W^{5,2}_{0}},\,\,\,||u^{\bot}||_{W^{5,2}_{0}}\le C ||u||_{W^{5,2}_{0}}
\end{equation*} 
Then we can use the Lemma \ref{improve decay} to get that:
\begin{equation}\label{e 5.9}
||u^{\bot}||_{W_{-1}^{4,2}} \le ||u^{\bot}||_{W^{5,2}_{0}}.
\end{equation}
Then, we can get that 
\begin{equation*}
||Re L u^{\bot}||_{W^{0,2}_{-1}} \le ||u^{\bot}||_{W_{-1}^{4,2}}  < \infty.
\end{equation*}
Combining this and the fact that $Re L u \in C_{-\eta }^{0,\alpha} \subset W_{-\eta -\frac{1}{2}+\epsilon}^{0,2}$ for any $\epsilon>0$, we get that 
\begin{equation}\label{e 5.10}
Re L u_0 \in W_{-\eta-\frac{1}{2} +\epsilon}^{0,2}.
\end{equation}
Let $\epsilon$ small such that $\epsilon \le \eta$. Without loss of generality, we can assume that $\eta <\frac{1}{2}$. Then we can apply the Lemma \ref{improve decay sobolev} with $\delta_1=-\eta$ and $\delta_2=0$ to get that:
\begin{equation*}
u_0 \in W_{-\eta}^{4,2}.
\end{equation*}
Without loss of generality, we can assume that there doesn't exist an integer $k$ such that $k \eta = \frac{1}{2}$. Denote $k_0$ as the biggest integer such that $k_0 \eta < \frac{1}{2}$. Then we can repeat the above argument to get that:
\begin{equation*}
u_0 \in W_{-k_0 \eta}^{4,2}.
\end{equation*}
We can let $\epsilon$ be small enough such that $-\eta -\frac{1}{2}+\epsilon < -(k_0+1)\eta$.  Then we can use the Lemma \ref{improve decay sobolev 2} with $\delta_2= -k_0 \eta$ and $\delta_1 = - (k_0+1) \eta$ to get that:
\begin{equation*}
u_0 \in W_{-(k_0+1)\eta}^{4,2} \oplus \chi(t) p^* Ker L_D.
\end{equation*}
Then we can write:
\begin{equation*}
u_0 = \widetilde{u} + \Sigma_{i=1}^N p^* u_i \chi(t),
\end{equation*}
where $\widetilde{u} \in W_{-(k_0+1)\eta}^{4,2} $ and $u_i \in Ker L_D$. Next, we want to use H\"older space instead of Sobolev space.
Using the Lemma \ref{AA0}, we have that: 
\begin{equation}\label{e 5.11}
(Re L-Re L^0) \Sigma_{i=1}^N p^*u_i \chi(t) \in C_{-\eta}^{0,\alpha}.
\end{equation}
We also calculate that:
\begin{equation*}
Re L^0 \Sigma_{i=1}^N p^*u_i \chi(t) =\Sigma_{i=1}^N Re L_D u_i \chi +\Sigma_{i=1}^N \Delta u_i (\partial_t -\partial_t^2)\chi +\Sigma_{i=1}^N u_i [\frac{1}{2}(\partial_t-\partial_t^2)^2 +(\partial_t -\partial_t^2)]\chi .
\end{equation*}
Since $\chi=1$ in a neighbourhood of $D$, the second term and the third term on the right-hand side of the above equation is zero in a neighbourhood of $D$. Since $u_i\in Ker(Re L_D)$, the first term on the right hand side of the above equation is zero. As a result, $Re L^0 \Sigma_{i=1}^N p^*u_i \chi(t) \in C_{-\eta}^{0,\alpha}$. Combining this with (\ref{e 5.11}), we have that $Re L \Sigma_{i=1}^N p^*u_i \chi(t) \in C_{-\eta}^{0,\alpha}$. Then we have that $Re L( u^{\bot}+ \widetilde{u})= Re Lu- Re L \Sigma_{i=1}^N p^*u_i \chi(t) \in C^{0,\alpha}_{-\eta}$. Since we have that $u^{\bot} \in W^{4,2}_{-1}$ and $\widetilde{u}\in W_{-(k_0+1)\eta}^{4,2}$, we have that $u^{\bot}+ \widetilde{u}\in W_{-(k_0+1)\eta}^{4,2}$. Then
we apply the lemma \ref{l regularity} to $u^{\bot}+ \widetilde{u}$ to get that 
\begin{equation}\label{e 5.12}
u^{\bot}+ \widetilde{u}\in C_{-(k_0+1)\eta +\frac{1}{2}}^{4,\alpha},
\end{equation}
where $-(k_0+1)\eta +\frac{1}{2}<0$. 
Now, we want to improve the regularity of $u^{\bot}+ \widetilde{u}$ from $C^{4,\alpha}_{-(k_0+1)\eta +\frac{1}{2}}$ to $C^{4,\alpha}_{-\eta}$. We can apply the Lemma \ref{improve decay holder} to $u^{\bot}+ \widetilde{u}$ with $\delta_2= -(k_0+1)\eta +\frac{1}{2}$ and $\delta_1 = -\eta$ to get that
\begin{equation*}
u^{\bot}+ \widetilde{u} \in C_{-\eta}^{4,\alpha} \oplus \chi(t) p^* Ker Re L_D. 
\end{equation*}
Using the (\ref{e 5.12}), we have that $u^{\bot}+ \widetilde{u}$ goes to zero near $D$. As a result, $u^{\bot}+ \widetilde{u}$ doesn't have a nonzero component in $\chi(t) p^* Ker Re L_D$.  Then we have that:
\begin{equation*}
u^{\bot}+ \widetilde{u} \in C_{-\eta}^{4,\alpha}.
\end{equation*}
This concludes the proof of the Proposition \ref{solve l operator}.

\subsection{Decomposition of a modified weighted H\"older space using Lichnerowicz operator}

For any $\delta \in \bR$, we can define the following modified weighted H\"older space:
\begin{equation*}
\widetilde{C}^{k,\alpha}_{\delta}(X \setminus D) \triangleq C_{\delta}^{k,\alpha} (X \setminus D) \oplus \chi p^* C^{k,\alpha}(D) .
\end{equation*}
Here $\chi$ is a cut-off function supported in a small neighborhood of $D$ and is equal to $1$ in a smaller neighborhood of $D$.

In this subsection, we want to prove the following Theorem which is the Theorem \ref{holder dec 12} in the Introduction section:
\begin{prop}\label{holder dec}
Suppose that $\omega$ is a Poincar\'e type K\"ahler metric satisfying (\ref{e 5.1}) with $\eta < \frac{1}{2}$. Then there exists a constant $\delta_1>0$ such that for any $\eta_0 \in (0, \delta_1)$, we have that:
\begin{equation*}
\widetilde{C}_{-\eta_0}^{1,\alpha}= Ker Re L |_{\widetilde{C}_{-\eta_0}^{5,\alpha}} \oplus Re L(t \chi (p^* Ker Re L_{D}) ) \oplus Re L(\widetilde{C}_{-\eta_0}^{5,\alpha}).
\end{equation*}
\end{prop}

Before proving the above Proposition, first we need the following Lemma:
\begin{lem}\label{cusp g estimate}
Suppose that $\omega$ is a Poincar\'e type K\"ahler metric satisfying (\ref{e 5.1}) with $0< \eta \le 1$. Then in any cusp coordinate, there exists a constant $C$ such that
\begin{equation*}
|g_{z^n \bar z^{\alpha} }| \le C\frac{1}{|z_n||\log |z_n||^{1+\eta}}
\end{equation*}
for $\alpha \le n-1$ and
\begin{equation*}
|g_{z^n \bar z^n}|\le C\frac{1}{|z_n|^2 \log^2 |z_n|}.
\end{equation*}
\end{lem}
\begin{proof}
We want to use the Quasi coordinates to prove this Lemma. Let $\varphi_{\delta}$, $\delta \in (0,1)$, be the map define in the section 3.3. Then we have that
\begin{equation*}
\partial_{\xi}= \partial_{z^n} \frac{\partial \varphi_{\delta}}{\partial \xi} =  exp(-\frac{1+\delta}{1-\delta}\frac{1+\xi}{1-\xi}) (\frac{-(1+\delta)}{1-\delta}) \frac{2}{(1-\xi)^2}\partial_{z^n}.
\end{equation*}
As a result, we have that
\begin{equation}\label{gznzalpha}
|g_{z^n \bar z^\alpha }|= |g_{\xi \bar z^\alpha }| |exp(\frac{1+\delta}{1-\delta}\frac{1+\xi}{1-\xi}) \frac{(1-\xi)^2(1-\delta)}{2(1+\delta)}| \le  C |g_{\xi \bar z^\alpha }| \frac{1}{|z^n| |\log |z^n||}.
\end{equation}
Using (\ref{e 5.1}) and   
\begin{equation*}
        \Phi_{\delta}^* \omega_{0}= \frac{\sqrt{-1}d \xi \wedge d \bar \xi}{(1-|\xi|^2)^2} +\Sigma_{i=1}^{n-1} \sqrt{-1} dz^i \wedge d \bar z^i,
    \end{equation*}
we have that  $|g_{\xi \bar z^\alpha }| \le C \frac{1}{|\log |z^n||^{\eta}}$. Combining this with (\ref{gznzalpha}), we can concludes the proof of the first part of this Lemma. The second part of this Lemma follows from the fact that $\omega$ is quasi-isometric to the standard Poincar\'e type metric $\omega_0$.
\end{proof}

Then we can prove the following Lemma:
\begin{lem}\label{lem 7.21}
Suppose that $\omega$ is a Poincar\'e type K\"ahler metric satisfying (\ref{e 5.1}) with $0< \eta \le 1$. Then for any $\eta_0 \in [0,\eta]$, we have that $Ker ReL |_{\widetilde{C}^{5,\alpha}_{-\eta_0}} = Ker ReL |_{W^{5,2}_{0}}$.
\end{lem}
\begin{proof}
For any $h\in Ker ReL |_{W^{5,2}_{0}}$, we have that $\nabla_{\omega}^{1,0}h =V \in \mathbf{h}^D_{\parallelsum}$, according to the Lemma \ref{ker a}. So $V|_D$ is a holomorphic vector field parallel to $D$. Take a cusp coordinate $(z)$. Denote $V=v^i \partial_{z^i}$. Then we can get that $|v^n|\le C|z_n|$ and $|v^{\alpha}|\le C$ for $\alpha \le n-1$ because $v\in \mathbf{h}^D_{\parallelsum}$. This implies that
\begin{equation}\label{e 7.14}
|h_n|= |v^{\bar k} g_{\bar k n}|\le \Sigma_{\alpha=1}^{n-1} |v^{\bar \alpha} g_{\bar \alpha n}| + |v^{\bar n} g_{n \bar n}|\le \frac{1}{|z_n||\log |z_n||^{1+\eta}}.
\end{equation}
Here we use the Lemma \ref{cusp g estimate}. 
Since 
\begin{equation}\label{e 7.15}
\int_0^s \frac{1}{\lambda |\log \lambda|^{1+\eta}} d\lambda = \frac{1}{|\log s|^{\eta}} 
\end{equation}
 which goes to zero as $s$ goes to zero, we have that $h$ can be extended continuously to $D$.  Again we can use $\nabla^{1,0} h |_D = V|_D $ to get that $h|_D$ is smooth on $D$. Combining this with (\ref{e 7.15}), we have that
 \begin{equation*}
  h-p^* h|_D=O(e^{-\eta t}).
 \end{equation*}
 Since $L h=0$ and $L \chi p^* (h|_D) \in C_{-\eta}^{1,\alpha}$, we can use some standard elliptic estimates in quasi coordinates to get that 
 \begin{equation*}
 h-p^* h|_D \in C_{-\eta}^{5,\alpha}.
 \end{equation*}
 This implies that $h \in \widetilde{C}_{-\eta}^{5,\alpha} \subset \widetilde{C}^{5,\alpha}_{-\eta_0}$ which concludes the proof of this Lemma.
\end{proof}

We also need the following Lemma proved by Sektnan in \cite{S}:
\begin{lem}\label{ker image}
Let $\omega$ be a Poincar\'e type metric on $X\setminus D$ satisfying  (\ref{e 5.1}) with $\eta>0$. Then there exists $\eta_0>0$ such that for all $\widetilde{f}\in Ker Re L_D$ there exists $\sigma \in C^{0,\alpha}_{-\eta_0}$, $\phi \in C^{4,\alpha}(D)$ and $f\in Ker ReL_D$ such that
\begin{equation*}
ReL_D (\chi p^* \phi + t \chi p^* f)= \chi p^* \widetilde{f}+ \sigma.
\end{equation*}
Moreover, $f$ is unique and $\varphi$ is unique up to an element of $Ker Re L_D$. Finally, if $\widetilde{f}=1$, we can take $f=1$ and $\phi=0$.
\end{lem}

Now we are ready to prove the Proposition \ref{holder dec}:
\begin{proof}
(of the Proposition \ref{holder dec})
.Let $\{v_i\}_{i=1}^N$ be an orthogonal unit basis of $Ker ReL |_{\widetilde{C}_{-\eta_0}^{5,\alpha}}$ with respect to the $L^2$ norm defined by $\omega$. For any $f\in \widetilde{C}_{-\eta_0}^{1,\alpha}$, we have that $f- \Sigma_{i=1}^N <f,v_i> v_i \in \widetilde{C}_{-\eta_0}^{1,\alpha}$ is perpendicular to $Ker Re L |_{\widetilde{C}_{-\eta_0}^{5,\alpha}}$. Let $v$ be the function defined on $D$ such that
\begin{equation}\label{holder dec 1}
f- \Sigma_{i=1}^N <f,v_i> v_i |_D =v.
\end{equation}
Using the Fredholm alternative, we have that 
\begin{equation*}
C^{1,\alpha}(D)=  ReL_D |_{C^{5,\alpha}(D)} \oplus Ker ReL_D |_{C^{5,\alpha}(D)}.
\end{equation*}
Then we can decompose $v$ as
\begin{equation}\label{holder dec 2}
v=u_1 + ReL_D(u_2),
\end{equation}
where $u_1,u_2 \in C^{5,\alpha}_D$ and $u_1 \in Ker ReL_D |_{C^{5,\alpha}_D}$. 
Using the Lemma \ref{ker image} below, we can find $u_3\in C^{5,\alpha}(D)$ and $u_4 \in Ker ReL_D$ such that
\begin{equation}\label{holder dec 3}
\chi p^* u_1 - Re L(\chi p^* u_3 + t \chi p^* u_4) \in C^{1,\alpha}_{-\eta_0}.
\end{equation}

According to the asymptotic behaviour of $Re L$, i.e. (\ref{AA0 holder difference}), we have that
\begin{equation}\label{holder dec 4}
Re L (\chi p^* u_2)- \chi p^* Re L_D u_2 \in C_{-\eta_0}^{1,\alpha}.
\end{equation}

Combining (\ref{holder dec 1}), (\ref{holder dec 2}), (\ref{holder dec 3}) and (\ref{holder dec 4}), we have that
\begin{equation*}
f- \Sigma_{i=1}^N <f,v_i>  v_i - ReL (\chi p^* u_2) - Re L(\chi p^* u_3 + t \chi p^* u_4) \in C_{-\eta_0}^{1,\alpha}
\end{equation*}

Using the expression of Poincar\'e type K\"ahler metrics using $t$ variable, i.e. (\ref{t coordinate}), we can get that $t \chi p^* u_4 \in W_0^{5,2}$. Thus we have that $Re L(\chi p^* u_2+ \chi p^* u_3 + t \chi p^* u_4)$ is perpendicular to $Ker ReL |_{\widetilde{C}_{-\eta_0}^{5,\alpha}}$. As a result, $f- \Sigma_{i=1}^N <f,v_i> v_i  - Re L (\chi p^* u_2) - Re L(\chi p^* u_3 + t \chi p^* u_4) $ is also perpendicular to $Ker ReL |_{\widetilde{C}_{-\eta_0}^{5,\alpha}}$. Then we can apply the Proposition  \ref{solve l operator} to get a function $u\in \widetilde{C}_{-\eta_0}^{5,\alpha}$ such that 
\begin{equation}\label{Au decom}
ReLu = f- \Sigma_{i=1}^N <f,v_i>  v_i - ReL (\chi p^* u_2) - ReL(\chi p^* u_3 + t \chi p^* u_4)  .
\end{equation}
(\ref{Au decom}) implies that 
\begin{equation*}
\widetilde{C}_{-\eta_0}^{1,\alpha} \subset Ker ReL |_{\widetilde{C}_{-\eta_0}^{5,\alpha}} +ReL(t \chi (p^* Ker Re L_D)) + ReL(\widetilde{C}_{-\eta_0}^{5,\alpha}).
\end{equation*}
We want to show that
\begin{equation}\label{subset 2}
Ker ReL |_{\widetilde{C}_{-\eta_0}^{5,\alpha}} +ReL(t \chi (p^* Ker Re L_D))   + ReL(\widetilde{C}_{-\eta_0}^{5,\alpha}) \subset \widetilde{C}_{-\eta_0}^{1,\alpha}.
\end{equation}
In fact, we just need to show that $ReL(t \chi (p^* Ker Re L_D))  \subset  \widetilde{C}_{-\eta_0}^{1,\alpha}$. First, we can calculate that for any $v\in Ker Re L_D$,
\begin{equation*}
Re L^0 (t \chi p^* v)=  p^* \Delta_{\omega_D} v +  t p^* Re L_D v +\widetilde{v}= \chi p^* \Delta_{\omega_D}  v  +\widetilde{v},
\end{equation*}
where $\widetilde{v}$ is a function which is zero in a neighbourhood of $D$. Here we use that $\chi$ is equal to 1 in a neighbourhood of $D$ and $v\in Ker Re L_D$.
Thus we have that 
\begin{equation}\label{l0 v}
Re L^0 (t \chi p^* v) \in \widetilde{C}^{\infty}_{-\eta_0}.
\end{equation}
Using the Lemma \ref{AA0}, we can get that
\begin{equation}\label{l0l v}
||Re L(t \chi p^* v)-Re L^0 (t \chi p^* v)||_{C_{-\eta_0}^{1,\alpha}}  \le C ||t \chi p^* v||_{C^{5,\alpha}_{\eta -\eta_0}} \le C  ||t \chi p^* v||_{C^{5,\alpha}_{\epsilon}} 
\end{equation}
In the second inequality above, we assume that $\eta_0$ is small enough such that $\eta- \eta_0 \ge \epsilon>0$ for some small constant $\epsilon$ without loss of generality. Using (\ref{e 5.1}) and (\ref{t coordinate}), we can see that 
\begin{equation}\label{v15 bdd}
||t \chi p^* v||_{C^{5,\alpha}_{\epsilon}} < + \infty.
\end{equation}
Then we can combine (\ref{l0 v}), (\ref{l0l v}) and (\ref{v15 bdd}) to get that
 \begin{equation*}
 ||Re L(t \chi p^* v)||_{C_{-\eta_0}^{1,\alpha}} < + \infty.
 \end{equation*}
This finishes the proof of (\ref{subset 2}). Then we have that
\begin{equation*}
\widetilde{C}_{-\eta_0}^{1,\alpha} =Ker ReL |_{\widetilde{C}_{-\eta_0}^{5,\alpha}} + ReL(t \chi (p^* Ker Re L_D))   + ReL(\widetilde{C}_{-\eta_0}^{5,\alpha}).
\end{equation*}
In order to show that the $+$ above is in fact $\oplus$, we need to show that if there exists $u \in Ker ReL |_{\widetilde{C}_{-\eta_0}^{5,\alpha}} $, $\rho \in Ker ReL_D$, $v\in \widetilde{C}_{-\eta_0}^{5,\alpha}$ such that 
\begin{equation}\label{uA0}
u + ReL(p^* \rho t \chi + v)=0,
\end{equation}
 then we have that $u=ReL(v)=\rho=0$. Since $u \in Ker ReL|_{W_0^{5,2}}$ and $ReL(p^* \rho t \chi + v) \in Im ReL|_{W_0^{5,2}}$, (\ref{uA0}) implies that $u=0$ and $ReL(p^* \rho t \chi + v)=0$. Then we use the Lemma \ref{lem 7.21} to get that $p^* \rho t \chi + v \in \widetilde{C}_{-\eta_0}^{5,\alpha}$. This implies that $p^* \rho t \chi =0$. Thus, we have that $ReLv=0$. This concludes the proof of this Proposition.
\end{proof}

\section{Compactness of isometry group}
In this section, we want to prove the following Theorem:
\begin{thm}\label{compactness of isometry group}
Suppose that $D$ is a smooth divisor. Suppose that $Aut_0(D)=\{Id\}$. Let $\omega$ be a Poincar\'e type extremal K\"ahler metric. Then the isometry group $Iso_0^D (X,\omega)$ is a compact set in $Aut_0^D(X)$.
\end{thm}
In order to prove the compactness of $Iso_0^D (X,\omega)$, we need to get uniform control on elements in $Iso_0^D (X,\omega)$ both in the interior of $X \setminus D$ and near $D$. 

For the control on elements in $Iso_0^D (X,\omega)$ in the interior of $X \setminus D$, we will prove the following Proposition:
\begin{prop}\label{compact to compact}
Let $\omega$ be a Poincar\'e type K\"ahler metric. Then the following holds:
 \begin{enumerate}
 \item For any compact set $K\subset X\setminus D$, there exists a compact set $K' \subset X\setminus D$ such that $g (K) \subset K'$ and $g^{-1}(K) \subset K'$ for any $g\in Iso_0^D(X, \omega)$.
 \item  Let $g_k$ be a sequence of biholomorphisms in $Iso_0^D(X,\omega)$. Then there exists a biholomorphism $g$ from $X\setminus D$ to itself such that after taking a subsequence, $g_k$ converge to $g$ locally compactly on $X\setminus D$. 
 \end{enumerate}
\end{prop}

For the control on elements in $Iso_0^D (X,\omega)$ near $D$, we first need to make some definitions. Let $\{U_i\}_{i\in \mathcal{A}}$ be a finite cover of $D$, where $U_i$ are coordinate balls on $D$. Then, the normal bundle $N_D$ over each $U_i$ is complex trivial. Fix a smooth K\"ahler metric $\omega_X$ on $X$. We identify $N_D$ as a subbundle of $TX_D$ consisting of vectors that are perpendicular to $TD$ with respect to $\omega_X$. Then, we can take a section $\sigma_i$ of $N_D(U_i)$ such that $|\sigma_i|_{\omega_X}=1$. We define a map from $U_i \times \{z\in \mathbb{C}: |z|\le \delta\}$ to $X$ as 
\begin{equation*}
\Phi_i (w',w_n) \triangleq exp_{w'}(w_n \sigma_i).
\end{equation*} 
Assuming that $\delta$ is small enough, this map is a diffeomorphism onto its image. We can also define
\begin{equation*}
\Phi_{i,w_n}(w') \triangleq \Phi_i (w',w_n).
\end{equation*}
For any automorphism $g$ of $X\setminus D$ which preserve $\omega$, we can define $g_{U_i,w_n}: U_i \rightarrow D$ by:
\begin{equation*}
g_{U_i,w_n} \triangleq p \circ g\circ \Phi_{i,w_n}.
\end{equation*}

Recall that $D$ can be written as $D= \Sigma_{i=1}^N D_i$ where $D_i$ are smooth connected divisors. Then we will prove the following Proposition which basically shows that $g$ sends a neighborhood of $D_i$ to be a neighborhood of $D_i$:
\begin{prop}\label{g vi}
Suppose that $D$ is a smooth divisor. For any $i\le N$ and any open neighbourhood $V_i$ of $D_i$, for any $j \in \mathcal{A}$ such $U_j \in D_i$, there exists a constant $\delta_0>0$ such that for any $|w_n|\le \delta_0$  and $g \in Iso_0^D(X,\omega)$, we have that 
\begin{equation*}
g \circ \Phi_{j,w_n} (U_j)\subset V_i. 
\end{equation*} 
Moreover, we have that $g_{U_j,w_n} (U_j) \subset D_i$.
\end{prop}

Then we can prove the following Propositions:
 
\begin{prop}\label{gk converge}
Assume that $Aut_0(D)=\{ Id \}$. Assume that $D$ is  smooth. Then for any $\epsilon>0$, there exists $\delta>0$ such that for any $i$, $|w_n|\le \delta$ and any $g \in Iso_0^D(X,\omega)$, we have that:
\begin{equation*}
|g_{U_i,w_n}- Id|\le \epsilon.
\end{equation*}
\end{prop}

$Id$ above is the identity map on $D$. According to the Proposition \ref{g vi}, the image of $g_{U_i, w_n}$ lie in the same divisor $D_j$ as $U_i$ for some $i$. As a result, $|g_{U_i,w_n}- Id|$ is well defined using a distance function on $D_j$.

We will prove the above propositions in the following subsections. Next, we use the above Propositions to prove the main theorem in this section:
\begin{proof}
(of the Theorem \ref{compactness of isometry group}). 
For any $\epsilon>0$, we can let $\delta$ be small and use the Proposition \ref{gk converge} to get that;
\begin{equation}\label{e 4.2}
| g_{U_i,w_n}- Id|\le \epsilon.
\end{equation}
Note that the Proposition \ref{g vi} controls $g(z)$ in the normal direction of $D$, making $g(z)$ close to $D$, while (\ref{e 4.2}) controls $g(z)$ in the parallel direction of $D$, making $p(g(z))$ close to $p(z)$. Combining the Proposition \ref{g vi}  and (\ref{e 4.2}) together and letting $U$ be close to $D$ depending on $\epsilon$ and making $\delta$ be small, we can get that:
\begin{equation}\label{e 4.3}
d(g(z),z)_{\omega_X} \le C \epsilon
\end{equation}
for any $z\in \Phi_i(U_i \times \{w_n: |w_n|\le \delta\})$ , where $C$ is independent of $g$. Here $d(\cdot, \cdot)_{\omega_X}$ means the distance function defined using a smooth K\"ahler metric $\omega_X$ on $X$.  For any subsequence $\{g_k\}$ in $Iso_0^D (X,\omega)$, we can use the Proposition \ref{compact to compact} to get an automorphism $\widetilde{g}$ of $X\setminus D$ such that  $g_k$ converge to $\widetilde{g}$ locally uniformly on $X \setminus D$ and $\widetilde{g}^* \omega = \omega$.
Take $g=g_k$ in (\ref{e 4.3}) and let $k \rightarrow \infty$. Then we get that 
\begin{equation}\label{e 4.4}
d(\widetilde{g}(z),z)_{\omega_X} \le C \epsilon
\end{equation}
for any $z\in \Phi_i(U_i \times \{w_n: |w_n|\le \delta\})$. This implies that $\widetilde{g}$ can be continuously extended to $D$ and $\widetilde{g}|_D = Id$. Combining this with the fact that $\widetilde{g}$ is holomorphic on $X\setminus D$, we have that $\widetilde{g}$ is holomorphic on $X$. Since $X\setminus (\cup_i \Phi_i(U_i \times \{w_n: |w_n|< \delta\}))$ is a compact set in $X\setminus D$, we have that $g_k$ converge to $\widetilde{g}$ on $X\setminus (\cup_i \Phi_i(U_i \times \{w_n: |w_n|< \delta\}))$. Combining this with (\ref{e 4.3}) and (\ref{e 4.4}), we get that $g_k$ converge to $\widetilde{g}$ uniformly on $X$. This implies that $\widetilde{g} \in Iso_0^D(X, \omega)$ and thus  $Iso_0^D(X, \omega)$ is compact. 
\end{proof}

In this section we assume that the constant $a_j$ in the Lemma \ref{asym} is equal to $\frac{1}{2}$ just for the convenience of calculation.

\subsection{control of elements of isometry group in the interior of $X\setminus D$}
We want to prove the Proposition \ref{compact to compact} in this subsection. 

\begin{proof}
For any compact set $K \subset X\setminus D$, there exists a positive number $C_0>0$ such that for any $q\in K$, we have that
\begin{equation}\label{volume lower bound}
Vol(B_1(q)) \ge C_0.
\end{equation}
Here $Vol$ is the Volume defined using $\omega$ and $B_1(q)$ is the unit geodesic ball with respect to the metric $\omega$. Since $(X\setminus D, \omega)$ is a complete and noncompact manifold with finite volume, there exists an open neighborhood of $D$, denoted as $U$ such that for any $x\in U$, we have that
\begin{equation}\label{volume upper bound}
Vol(B_1 (x)) \le \frac{C_0}{2}.
\end{equation}
Since $g$ is a diffeomorphism and preserve $\omega$, we have that
\begin{equation}\label{volume equal}
Vol(B_1(g(q)))= Vol(g(B_1(q))) = Vol(B_1(q)).
\end{equation}
Combining (\ref{volume lower bound}), (\ref{volume upper bound}) and (\ref{volume equal}), we can get that for any $q\in K$, $g(q)\in X\setminus U$. Denote $K'= X\setminus U$, we finish the proof of (1). 
For any compact set $K  \subset (X \setminus D)$, there exists a compact set $K' \subset (X \setminus D)$ such that $g_k (K) \subset K'$ for any $k$, according to part (1). We can let $\epsilon$ be small enough and find a $\epsilon-$net $\{x_i\}_{i=1}^N$ of $K$ such that:
\begin{equation*}
K \subset  \cup_{i=1}^N B_{\epsilon}(x_i),
\end{equation*}
and
\begin{equation*}
\overline{B_{\epsilon}(x_i)} \subset X \setminus D.
\end{equation*}
We can also assume that each $B_{\epsilon}(x_i)$ is a coordinate ball. Since $g_k(x_i) \subset K'$ and $K'$ is compact, we can take a subsequence of $\{g_k\}$  (still denoted as $\{g_k\}$) such that for any $i$, there exists $y_i \in K'$ such that
\begin{equation*}
\lim_{k \rightarrow \infty}g_k(x_i)=y_i.
\end{equation*}
Since $g_k$ preserves the metric $\omega$, we have that $B_{\epsilon}(g_k(x_i))= g_k (B_{\epsilon}(x_i))$. Then we can let $k$ be big enough such that $g_k (B_{\epsilon}(x_i)) \subset B_{2 \epsilon}(y_i)$. We can assume that $\epsilon$ be small enough such that $B_{2\epsilon}(y_i)$ is contained in a coordinate ball and $\overline{B_{2\epsilon}(y_i)} \cap D = \emptyset$. Then we can use the coordinates of $B_{\epsilon}(x_i)$ and $B_{2\epsilon}(y_i)$ to see $g_k$ as holomorphic maps from a compact set of $\mathbb{C}^n$ to another compact set in $\mathbb{C}^n$, so we can get a subsequence of $\{g_k\}$ such that it converge on each $B_{\epsilon}(x_i)$. Then we can use a standard diagonal argument to get a subsequence of $\{g_k\}$ such that it converges to $g$ locally uniformly on $X \setminus D$ for some holomorphic map $g$ from $X\setminus D$ to $X\setminus D$. Repeat the above procedure on $\{g_k^{-1}\}$, we can find a holomorphic map $g'$ such that $g_k^{-1}$ converge to $g'$ locally uniformly on $X\setminus D$. We can see that $g'=g^{-1}$. So $g$ is an automorphism of $X \setminus D$. Since $g_k$ are holomorphic, we can use standard elliptic regularity result to get that the derivatives of $g_k$ of any order converge to that of $g.$ So we can use the fact that $g_k^* \omega = \omega$ to get that:
\begin{equation*}
g^* \omega =\omega.
\end{equation*}
This concludes the proof of part (2).
\end{proof}

\subsection{control of elements of isometry group near $D$}
 In this subsection, we want to prove the Proposition \ref{gk converge}.
 
\subsubsection{estimate about the exponential map}
  Since we heavily use the exponential map $\Phi_{i}(w)=exp_{w'}(w_n \sigma_i)$, we want to record estimates about this map. Let $(w)$ be a coordinate of $U_i \times \Delta^*$. Let $(z)$ be a cusp coordinate of $X$ containing $U_i$ with $D \cap U_i \subset \{z_n=0\}$. Let $J_0$ be the product complex structure on $U_i \times \Delta^*$. Let $J$ be the complex structure on $X$.  Then we have that:
\begin{equation*}
J dz_k = \sqrt{-1} dz_k
\end{equation*}
and
\begin{equation*}
J_0 dw_k =\sqrt{-1} dw_k.
\end{equation*}
 
The following lemma can be found in \cite{D}.
\begin{lem}\label{l5.13New}
The exponential map on a Hermitian manifold has the Taylor expansion in the following form under local coordinates:
\begin{equation*}
\exp_z(\zeta)_m=g_m(z,\xi)+\sum_{j,k,l}c_{jklm}(\frac{1}{2}\bar{z}_k+\frac{1}{6}\bar{\xi}_k)\xi_j\xi_l+O\big(|\xi|^2(|z|+|\xi|)^2),
\end{equation*}
where
$\zeta = \zeta_i \partial_{z^i}$,	
\begin{equation*}
g_m(z,\xi)=z_m+\xi_m-\sum_{j,l}a_{jlm}z_j\xi_l+\sum_{j,k,l,p}a_{jlp}a_{kpm}z_jz_k\xi_l-\sum_{j,k,l}b_{jklm}(z_jz_k\xi_l+z_k\xi_j\xi_l+\frac{1}{3}\xi_j\xi_k\xi_l),
\end{equation*}
and $\xi$ and $\zeta$ are related through:
\begin{equation*}
\xi_m=\zeta_m+\sum_{j,l}a_{jlm}z_j\zeta_l+\sum_{j,k,l}b_{jklm}z_jz_k\zeta_l.
\end{equation*}
In the above,  $(\exp_z\zeta)_m$ denotes the $m$-th component of the exponential map under local coordinates. 
\end{lem}

\begin{lem}\label{zw coordinates}
Let $A$ be big enough. Let $t$, $p$ be defined in the section 3.7 above. Then for any $q\in \mathcal{N}_A$, there exists $i$ such that $q\in \Phi_i(U_i \times \Delta^*)$ and there exist a coordinate $(w)$ of $U_i \times \Delta^*$ and a cusp coordinate $(z)$ and a constant $C$ depending on $\omega$, $X$ and $D$  such that:
\begin{enumerate}
\item $(\omega_X)_{i, \bar j}(p(q))= \delta_{ij}$.
\item $\sigma_i(p(q))= \partial_{z^n}$.
\item $d w^{\alpha} - d z^{\alpha} =O(e^{-t})$ for any $\alpha \le n-1$ at $q$.
\item $\frac{1}{|w_n|(-\log |w_n|^2)}dw^n - \frac{1}{ |z_n|(-\log |z_n|^2)} dz^n =O(e^{-t})$ at $q$.
\item $\partial_{w^{\alpha}} - \partial_{z^{\alpha}} =O(e^{-t})$, for any $\alpha \le n-1$ at $q$.
\item $ |w_n|(-\log |w_n|^2) \partial_{w^n}- |z_n|(-\log |z_n|^2) \partial_{z^n} =O(e^{-t})$ at $q$.
\end{enumerate} 
Here $O(e^{-\eta t})$ is uniformly bounded independent of $q$.
\end{lem}
\begin{proof}
We can first take a normal coordinate $(\widetilde{z})$ of $\omega_X$
 at $p$ such that $D$ is tangent to the $(n-1)-$plane spanned by $\widetilde{z}^i$ for $i \le n-1$. Then we change the coordinate by replacing $\widetilde{z}^n$ by $\widetilde{z}^n - f(\widetilde{z}')$ where $f$ is a holomorphic function such that $D$ locally is a zero set of $\widetilde{z}^n=f(z')$, and take $z^i =\widetilde{z}^i$ for $i\le n-1$.  Since $f(0)=0$ and $f'(0)=0$, this change of coordinate won't change $\omega_X$ at $p$. As a result, (1) holds. Using $(1)$ and the assumption that $|\sigma_i|=1$ and $\sigma_i$ is perpendicular to $TD$ with respect to $\omega_X$, we can get that $\sigma_i(p(q))= e^{i\theta_0} \partial_{z^n}$ for some constant $\theta_0$.  Then we can change the coordinate of $(z)$ by replacing $z_n$ by $e^{-i\theta_0} z_n$ to make (2) hold. Next we define the coordinate $(w)$ as follows: We define $w^{\alpha}=z^{\alpha}|_{D}$ for $\alpha \le n-1$. Let $w^n$ be the standard coordinate on $\Delta^*$.
 
Using the Lemma \ref{l5.13New}, we can get that:
\begin{equation*}
z_i = w_i +O(|w'||w_n|+|w_n|^3).
\end{equation*}
Then we can calculate that at $p=(0, w_n)$,
\begin{equation}\label{zw coordinates 1}
\begin{split}
&d z^i = dw^i + \Sigma_{\alpha=1}^{n-1} (d w^{\alpha}  O(|w_n|) + d \bar w^{\alpha}  O(|w_n|)) + (d w^n+ d \bar w^n)O(|w_n|^2)\\
&d \bar z^i = d \bar w^i + \Sigma_{\alpha=1}^{n-1} (d w^{\alpha}  O(|w_n|) + d \bar w^{\alpha}  O(|w_n|)) + (d w^n+ d \bar w^n)O(|w_n|^2)
\end{split}
\end{equation}
for any $i\in \{1,2,...,n\}$. Denote 
\begin{equation*}
\beta_{\alpha} \triangleq dz^{\alpha}, \beta_{\bar \alpha} \triangleq d \bar z^{\alpha},  \widetilde{\beta}_{\alpha} \triangleq dw^{\alpha},  \widetilde{\beta}_{ \bar \alpha} \triangleq d \bar w^{\alpha}
\end{equation*}
for $\alpha \le n-1$, and
\begin{equation}\label{zw coordinates 2}
\begin{split}
& \beta_n \triangleq \frac{1}{|z^n| (-\log |z^n|)} dz^n, \beta_{\bar n} \triangleq \frac{1}{|z^n| (-\log |z^n|)} d \bar z^n, \\
& \widetilde{\beta}_n \triangleq \frac{1}{|z^n| (-\log |z^n|)} d w^n, \widetilde{\beta}_{\bar n} \triangleq \frac{1}{|z^n| (-\log |z^n|)} d \bar w^n
\end{split}
\end{equation}
Using (\ref{zw coordinates 1}), we can get
\begin{equation}\label{zw coordinates 3}
(\beta_1, \beta_{\bar 1}, ..., \beta_n, \beta_{\bar n}) = A (\widetilde{\beta}_1, \widetilde{\beta}_{\bar 1},..., \widetilde{\beta}_n, \widetilde{\beta}_{\bar n}),
\end{equation}
where $A$ is a $2n \times 2n$ matrix satisfying $A= Id+ O(e^{- t})$. Then we can get that $A^{-1}= Id +O(e^{-\eta t})$. According to the Lemma \ref{asym} and assuming that the constant $a_j$ for $\omega$ in that lemma is equal to $\frac{1}{2}$, we can get that $\{\beta_i, \beta_{\bar i}\}$ are almost a unit orthogonal basis up to an error $O(e^{-\eta t})$. Then we can use (\ref{zw coordinates 3}) to get that $\{\widetilde{\beta}_i, \widetilde{\beta}_{\bar i}\}$ is also almost a unit orthogonal basis up to an error $O(e^{-\eta t})$. In particular, they are bounded. Combining this with (\ref{zw coordinates 3}) and the fact that $z_n= w_n +O(|w_n|^3)$ at $q$. we can thus finish the proof of (3) and (4) of this Lemma. (5) and (6) of this Lemma can be proved in a similar way.
\end{proof}
 
 \begin{lem}\label{JJ0}
 Let $\eta$ be the constant given in the Lemma \ref{asym}. Then,
 $J$ is assymptotic to $J_0$ in the sense that:
 \begin{equation*}
 J_{\Phi_i} \triangleq \Phi_i^*J=J_0 +O(e^{-\eta t}),
 \end{equation*}
 for any $i$.
 \end{lem}
 \begin{proof}
Let $(z)$ and $(w)$ be coordinates given by the Lemma \ref{zw coordinates}. Denote $v^{\alpha}= \partial_{z^{\alpha}}$ and $\widetilde{v}^{\alpha}= \partial_{w^{\alpha}}$ for $\alpha \le n-1$. Denote $v^n = |z_n|(-\log |z_n|^2) \partial_{z^n}$ and $\widetilde{v}^n = |w_n|(-\log |w_n|^2)\partial_{w^n}$.  Using the Lemma \ref{asym}, we have that 
\begin{equation*}
|<v^i, v^j>_{\omega} - \delta_{ij}| \le C e^{-\eta t}.
\end{equation*}
Then we can use the Lemma \ref{zw coordinates} to get that:
\begin{equation*}
|v^i- \widetilde{v}^i|\le C e^{-\eta t}.
\end{equation*}
Thus we have that
\begin{equation*}
|<\widetilde{v}^i, \widetilde{v}^j> -\delta_{ij}| \le C e^{-\eta t}.
\end{equation*}
As a result, $\{v^i\}$ and $\{\widetilde{v}^i\}$ are both almost unitary bases, and they are close to each other. Note that $J$ and $J_0$ are defined by
\begin{equation*}
J v^i = \sqrt{-1} v^i,
\end{equation*}
\begin{equation*}
J_0 \widetilde{v}^i = \sqrt{-1} \widetilde{v}^i.
\end{equation*}
Then, we can get that:
 \begin{equation*}
 \Phi_i^*J=J_0 +O(e^{-\eta t}).
 \end{equation*}
 \end{proof}

Then we can prove the following Lemma:
\begin{lem}\label{lem 4.15}
Suppose that $g\in Iso_0^D(X,\omega)$. Denote $g_{\Phi_i}= \Phi_i^{-1} \circ g \circ \Phi_i$ as the pull back of $g$ using $\Phi_i$. Then we have that $g_{\Phi_i}$ is almost holomorphic with respect to $J_0$ in the sense that:
\begin{equation*}
    D g_{\Phi_i} \circ J_0 \circ D g_{\Phi_i}^{-1}(x) - J_0 (x)= O(e^{- \eta \min\{t_1,t_2\}})
\end{equation*}
Here $t_1$ is the value of the function $t$ at $x$ and $t_2$ is the value of the function $t$ at $g_{\Phi}^{-1}(x)$. 
\end{lem}
\begin{proof}
Since $g$ is $J-$holomorphic, we have that:
\begin{equation*}
Dg \circ J \circ D g^{-1} =J.
\end{equation*}
Using the pull back with $\Phi_i$, we get that:
\begin{equation*}
D g_{\Phi_i} \circ J_{\Phi_i} \circ D g_{\Phi_i}^{-1} = J_{\Phi_i}.
\end{equation*}
Using the Lemma \ref{JJ0} and the fact that $g$ preserve $\omega$, we can get that:
\begin{equation*}
    D g_{\Phi_i} \circ J_0 \circ D g_{\Phi_i}^{-1}(x) - J_0 (x)= O(e^{- \eta \min\{t_1,t_2\}}).
\end{equation*}
\end{proof}

We can also prove the following Lemma:
\begin{lem}\label{normal bundle metric}
The pullback of the metric $g_{\omega}$ to $N_D$ is asymptotic to the product Poincar\'e type metric on $D \times \Delta^*$ in the sense that
\begin{equation*}
\Phi_i^* g_{\omega} = g_D + \frac{2|dw_n|^2}{|w_n|^2 \log^2(|w_n|^2)}+ O(e^{-\eta t}).
\end{equation*}
\end{lem}
\begin{proof}
According to the proof of the Lemma \ref{zw coordinates}, fix a point $p$,  we can choose appropriate coordinate $(w)$ for $N_D$ and cusp coordinate $(z)$ such that at $p$ we have that 
\begin{equation*}
dz^i = dw^i + \Sigma_{\alpha=1}^{n-1} (h^i_{\alpha} dw^{\alpha} + h^{i}_{\bar \alpha} d \bar w^{\alpha}) +h^i_n dw^n + h^i_{\bar n} d \bar w^n,
\end{equation*}
with $h^i_{\alpha}, h^i_{\bar \alpha} \in O(|w^n|)$, $h^i_n, h^i_{\bar n} \in O(|w^n|^2)$ and
\begin{equation*}
z_i = w_i + O(|w_n||w'|+|w_n|^3).
\end{equation*}
As a result, we have that at $w=(0,w_n)$,
\begin{equation*}
\begin{split}
&\frac{\sqrt{-1} 2 d z^n \wedge d\bar z^n}{|z_n|^2 \log^2 |z_n|^2} + \Sigma_{\gamma=1}^{n-1} \sqrt{-1}dz^{\gamma} \wedge d \bar z^{\gamma} \\
&=\frac{2 \sqrt{-1}(dw^n +h_{\alpha}^n dw^{\alpha} + h_{\bar \alpha}^n d \bar w^{\alpha} +h^i_n dw^n + h^i_{\bar n} d \bar w^n)\wedge (d \bar w^{n} + \overline{h_{\beta}^n}d\bar w^{\beta} + \overline{h_{\bar \beta}^n}dw^{\beta}+\overline{h^i_n} d \bar w^n + \overline{h^i_{\bar n}} d  w^n)}{(1+O(|w_n|^2))|w_n|^2 \log^2(|w_n|^2)}\\
&+ \Sigma_{\gamma=1}^{n-1} \sqrt{-1} (dw^{\gamma}+ h^{\gamma}_{\alpha}dw^{\alpha} +h^{\gamma}_{\bar \alpha}d \bar w^{\alpha} + h^{\gamma}_{n} dw^n + h^{\gamma}_{\bar n} d \bar w^n) \wedge (d \bar w^{\alpha} +\overline{h^{\gamma}_{\beta}}d\bar w^{\beta} + \overline{h^{\gamma}_{\bar \beta}}dw^{\beta}+ \overline{h^{\gamma}_n} d \bar w^n + \overline{h^{\gamma}_{\bar n}} dw^n)\\
& = \frac{2\sqrt{-1} dw^n \wedge d \bar w^n}{|w_n|^2 \log^2 |w_n|^2} + \Sigma_{\gamma=1}^{n-1} \sqrt{-1} dw^{\gamma} \wedge d \bar w^{\gamma} +  O(e^{-t}).
\end{split}
\end{equation*}
Then this Lemma follows from the above formula and the Lemma \ref{asym}.
\end{proof}

\subsubsection{uniform estimate about $g_{U_i,w_n}$}
Denote $g_{U_i,w_n} \triangleq p \circ g \circ \Phi_{i,w_n}$. First we want to prove the following Lemma:
\begin{lem}\label{image goes to d}
For any $A>0$, there exists $\delta_0$ such that for any $|w_n|\le \delta_0$, we have that $t  (g\circ \Phi_{i,w_n}  )\ge A$ for any $i$ and $g \in Iso_0^D(X,\omega)$.
\end{lem}
\begin{proof}
for any open neighbourhood $U$ of $D$ in $X$, we have that $K \triangleq X \setminus U$ is a compact set in $X \setminus D$. Then by the Proposition \ref{compact to compact}, we can find a compact set $K'$ such that for any $g\in Iso_0^D(X, \omega)$, we have that $g(K) \subset K'$. Thus $g^{-1}(K) \subset K'$ since $g^{-1}\in Iso_0^D(X, \omega)$. Then we can take $\delta$ small such that for any $i$, $\Phi_i(U_i \times \{w_n: |w_n|\le \delta\}) \subset X\setminus K'$. Then, for any $z\in \Phi_i(U_i \times \{w_n: |w_n|< \delta\})$, we have that:
\begin{equation}
g(z) \subset g(X\setminus K')= X \setminus g(K') \subset X \setminus K = U.
\end{equation}
So we see that as $w_{n}$ goes to zero, the image of $g(U_i \times \{w_{n}\})$ will go to $D$. This concludes the proof of this Lemma.
\end{proof}

Then we can prove the Proposition \ref{g vi}
\begin{proof}
(of the Proposition \ref{g vi})
Fix small open neighborhoods $V_i$ of $D_i$ for any $i \le N$ such that $V_i$ don't intersect with each other. Using the Lemma \ref{image goes to d}, we can find a constant $\delta_0>0$ such that for any $|w_n|\le \delta_0$ and $g \in Iso_0^D(X,\omega)$, we have that 
\begin{equation}\label{g in nbhd}
g \circ \Phi_{j,w_n} (U_j) \subset \cup_{i=1}^N V_i. 
\end{equation}
Since $g \in Iso_0^D(X,\omega)$, we have that $g |_D =Id$ which implies that $g \circ \Phi_{j,0} (U_j) \subset V_i$. Using the continuity of $g$,  (\ref{g in nbhd}) and the assumption that $V_i$ don't intersect with each other, we have that 
\begin{equation}
g \circ \Phi_{j,w_n} (U_j) \subset V_i. 
\end{equation}
This immediately implies that $g_{U_i,w_n} (U_i) \subset D_i$.
\end{proof}

Then, we can prove the following Lemma:
 \begin{lem}\label{gui estimate}
For any $g \in Iso_0^D(X,\omega)$, for any $i$ and $w_n \in B_{\delta_0}(0)$ with $\delta_0$ depending on $(X,\omega)$, we have that $g_{U_i, w_n}$ is locally diffeomorphic onto its image.
 \end{lem} 
\begin{proof}
For any $q_1\in \Phi_i(U_i \times \Delta^*)$, we can use the Lemma \ref{zw coordinates} to find two coordinates $(w)$ and $(z)$. Denote the value of $w_n$ at $q$ as $w^*$. Then $\{ \partial_{w^{\alpha}}\}_{\alpha=1}^{n-1}$ is a basis of $T_{q}(U_i \times \{w^*\})$.  In order to show that $g_{U_i,w_n}$ is a diffeomorphism, it suffices to show that $\{p_* g_*(\partial_{w^{\alpha}})\}$ are linearly independent. We will use geodesics to show this.
 
Fix a point $q_0 \in X\setminus D$. Let $\delta$ be a small constant to be determined. Let $A$ and $T$ be the constants given by the Lemma \ref{normal vector} below.  Since $(X, \omega)$ is a complete manifold, there exists a constant $A_1$ such that for any point $q_1$ with $t(q_1) \ge A_1$, we have that $d_{\omega}(q_1, \{q: t(q)=A\}) \ge T$. We can also use the fact that $(X, \omega)$ is a complete manifold to find a minimizing unit geodesic $\gamma$ connecting $q_1$ with $q_0$. Let $A$ be big depending on $q_1$ and let $q_1$ be close to $D$ depending on $A$. Then $\gamma$ intersect with $ \{q: t(q)=A\}$ at some point $q_2$. Denote $v= \nabla_s \gamma(q_1)$. Then, we can apply the Lemma \ref{normal vector} to show that $|\nabla_s (p \circ \gamma)|_{\omega}(q_1) \le \delta.$ According to the Lemma \ref{normal bundle metric}, $\Phi_i^* g_{\omega}$ is asymptotic to the standard Poincar\'e metric which is a product metric. As a result, we have that :
\begin{equation*}
|<v, \partial_{w^{\alpha}}>_{\omega}| \le 2\delta, \,\, |<Jv, \partial_{w^{\alpha}}>_{\omega}| \le 2\delta
\end{equation*}
for any $\alpha \le n-1$. 
Denote $v_1= g_* v$. Since $g$ preserve $J$ and $\omega$, we have that:
\begin{equation*}
|<Jv_1, g_* \partial_{w^\alpha}>_{\omega}|=|<g_* Jv, g_* \partial_{w^\alpha}>_{\omega}|=|<Jv, \partial_{w^{\alpha}}>_{\omega}| \le 2 \delta 
\end{equation*}
and
\begin{equation}\label{e v1 w}
|<v_1, g_* \partial_{w^\alpha}>_{\omega}|=|<g_* v, g_* \partial_{w^\alpha}>_{\omega}|=|<v, \partial_{w^{\alpha}}>_{\omega}| \le 2 \delta
\end{equation}
On the other hand,  since $g$ preserves $\omega$, we have that $g(\gamma)$ is a minimizing geodesic connecting $g(q_0)$ with $g(q_1)$. We have that  
\begin{equation*}
v_1=g_* v = \nabla_s (g\circ \gamma)(g(q_1)).
\end{equation*}
We can take coordinates $(\widetilde{w})$ and $(\widetilde{z})$ for $g(q_1)$ using the Lemma \ref{zw coordinates}. Using the Proposition \ref{compact to compact},  there exists a compact set $K \subset X\setminus D$ such that for any $g\in Iso_0^D(X,\omega)$, we have that $g(q_0) \in K$.  Then we can use the Lemma \ref{image goes to d} such that the assumptions of  the Lemma \ref{normal vector} hold with $\gamma$ replaced by $g \circ \gamma$. Then we can get that:
\begin{equation}\label{e v1 tilde w}
|<v_1, \partial_{\widetilde{w}^{\alpha}}>| \le 2\delta, \,\, |<Jv_1, \partial_{\widetilde{w}^{\alpha}}>| \le 2\delta
\end{equation} 
Combining (\ref{e v1 w}) and (\ref{e v1 tilde w}) and the fact that $v_1$ is a unit vector, we have that:
\begin{equation*}
dist\{span<g_* \partial_{w^1},...,g_* \partial_{w^{n-1}}>, span<\partial_{\widetilde{w}^1},...,\partial_{\widetilde{w}^{n-1}}>\}\le C \delta .
\end{equation*}
Here $dist\{\Sigma_1, \Sigma_2\}$ for two $k-$dimensional planes $\Sigma_1, \Sigma_2$ is defined as:
\begin{equation*}
dist\{\Sigma_1, \Sigma_2\} =\sup_{v\in \Sigma_2, |v|=1} d(\Sigma_1,v)+ \sup_{v\in \Sigma_1, |v|=1} d(\Sigma_2,v)
.\end{equation*}
We let $\delta$ small  such that $C\delta \le \frac{1}{2}$. Then we have that $\{p_* \circ g_* \partial_{w^{\alpha}}\}$ are linearly independent.
\end{proof}

 In the proof of the Lemma \ref{gui estimate} above, we use the following Lemma:
\begin{lem}\label{normal vector}
Let $q_1, q_2 \in \mathcal{N}_A$ with $t_i= t(q_i)$. Assume that $t_1 \ge t_2$.  Let $\gamma$ be a unit speed minimizing geodesic with $\gamma(0)=q_2$ and $\gamma(T)=q_1$. Denote $v= \nabla_s (p \circ \gamma)(q_1)$. Then for any $\delta>0$, we can let $A$ and $T$ be big depending on  $\diam (D)$ and $\delta$ such that $|v'|\le \delta$.
 \end{lem}
\begin{proof}
First, we let $\widetilde{\gamma}$ be a minimizing geodesic on $D$ with respect to $\omega_D$ such that $\widetilde{\gamma}(0)=p(q_2)$ and $\widetilde{\gamma}(T)= p(q_1)$. In particular, we have that $|\widetilde{\gamma}'|= \frac{d_D(p(q_1), p(q_2))}{T}$. Assume that $p(q_2) \in U_i$ for some $i\in \mathcal{A}$. We can let $A$ be big enough such that $q_2 \in \Phi_i( U_i \times \Delta^*)$, where $\Phi_i$ is defined before. Denote $w_n (q_2)$ as the projection of $\Phi_i^{-1}(q_2)$ to $\Delta^*$. Denote $w_n (q_2)= r_2 e^{i \theta_2}$. In the rest of the proof we use coordinate $(\widetilde{t}, \theta)$ for $\Delta^*$ such that $w_n = e^{-\frac{e^{\widetilde{t}}}{2}} e^{i \theta}$. Note that $t =\widetilde{t} +O(e^{-t})$. We let $\gamma_1$ be a minimizing geodesic connecting $q_1$ with $\Phi_i (p(q_2), t_1, \theta_2)$ with $\gamma_1(0)= q_1$ and $\gamma_1(1)= \Phi_i (p(q_2), t_1, \theta_2)$. Define $\gamma_2$ by
\begin{equation*}
\gamma_2(s)= \Phi_i(p(q_2), t_1 + \frac{t_2- t_1}{T} (s-1), \theta_2 ),
\end{equation*}
for $1\le s \le T+1$.
 Then we define $\gamma_3$ as
$$
\gamma_3(s)= \left \{
\begin{array}{rcl}
\gamma_1(s) &,\,\, &\text{ for } 0\le s \le 1 \\
\gamma_2(s) &, \,\, &\text{ for } 1\le s \le T+1
\end{array}
\right.
$$
Suppose that $|v'| >\delta$. We want to show that the length of $\gamma_3$ is shorter than the length of $\gamma$. This will contradict the fact that $\gamma$ is length-minimizing. According to the Lemma \ref{asym}, we have that 
\begin{equation*}
l(\gamma_1) \le \diam(D) +1,
\end{equation*}
if we let $A$ be big enough, and
\begin{equation*}
l(\gamma_2) \le  \int_0^T \frac{t_1- t_2}{T}(1+C e^{-\eta t_2}) ds = (1+C e^{-\eta t_2}) (t_1-t_2).
\end{equation*}
Then we have that:
\begin{equation*}
l(\gamma_3) =l(\gamma_2) +l(\gamma_1) \le (1+C e^{-\eta t_2})  (t_1-t_2) + \diam(D)+1.
\end{equation*}
Next, we estimate $t_1- t_2$ from above. Denote $\gamma_4\triangleq p(\gamma)$. Since $\gamma$ is a geodesic and $\omega$ is asymptotic to a product metric according to the Lemma \ref{asym}, we have that $\gamma_4$ is close to a geodesic in the sense that for any $\epsilon$ we can let $A$ be big enough which depends on $T$ and doesn't depend on $q_1,q_2$, such that:
\begin{equation*}
|\nabla_{\partial_s} \gamma_4'(s)| \le \epsilon.
\end{equation*}
This implies that 
\begin{equation*}
|\gamma_4'(s)| \ge \delta -C\epsilon.
\end{equation*}
Denote $s_1= \inf\{s_0\in [0,T]: t(\gamma(s))\ge t_2 \text{ for any }s\ge s_0\}$. Then we have that $t(\gamma(s_1))=t_2$. Denote $\gamma_5 \triangleq t \circ \gamma$. Then we can get that 
\begin{equation*}
t_1-t_2 \le  \int_{s_1}^T |\gamma_5'| \le \int_{s_1}^T(1+C e^{-\eta t_2})\sqrt{1-|\gamma_4'|^2} dt \le (1+C e^{- \eta t_2}) \sqrt{1- (\delta-C \epsilon)^2}T.
\end{equation*} 
Here we use the fact that for $s\ge s_1$, we have that $t(\gamma(s)) \ge t_2$, so the metric is a product metric up to an error $C e^{- \eta t_2}$, according to the Lemma \ref{asym} and the section 3.7. Then we have that:
\begin{equation*}
\begin{split}
l(\gamma_3) &\le (1+C e^{-\eta t_2})^2 \sqrt{1- (\delta-C \epsilon)^2}T  +\diam(D)+1.
\end{split}
\end{equation*}
Let $\epsilon$ be small depending on $\delta$ and let $T$ be big depending on $diam(D) $ and $\delta$ and let $A$ be big such that $t_2$ is big depending on $T$ and $\delta$. Then we get that 
\begin{equation*}
l(\gamma_3) \le (1-\frac{\delta^2}{3})T < T =l(\gamma).
\end{equation*}
This is a contraction with the fact that $\gamma$ is a minimizing geodesic.
\end{proof}

\begin{lem}\label{lem 4.18}
For any $\epsilon$, there exists $\delta_0>0$ such that for any $|w_n|\le \delta_0$ and  $g\in Iso_0^D(X,\omega)$, we have that:
\begin{equation*}
(1-\epsilon) g_D \le   g_{U_i,w_n}^* g_D \le (1+\epsilon) g_D.
\end{equation*}
Here $g_D$ is the Riemannian metric on $D$ with respect to $\omega_D$.
\end{lem}
\begin{proof}
According to the Lemma \ref{normal bundle metric}, we have that
\begin{equation*}
g_{\omega} |_{U_i\times\{w_n\}} = p^* g_D + O(e^{-\eta t}).
\end{equation*}
Since $g$ preserves $\omega$ and the complex structure $J$, it preserves $g_{\omega}$ which is the Riemannian metric with respect to $\omega$. In particular, we have that 
\begin{equation*}
g^* (g_{\omega}|_{g(U_i\times \{w_n\})})= g_{\omega} |_{U_i\times \{w_n\}}.
\end{equation*}
Using the Lemma \ref{gui estimate}, we have that $g_{U_i,w_n}$ is a local diffeomorphism. Note that $g|_{U_i\times \{w_n\}}$ is locally diffeomorphic onto its image. Combining this and the fact that  $g_{U_i, w_n}= p \circ g|_{U_i\times \{w_n\}}$, we get that  $p$ is a local diffeomorphism from $g(U_i\times \{w_n\})$ to $D$. Then for any $ q \in g(U_i\times \{w_n\})$ we can find a open neighbourhood $U$ of $q$ such that $p|_{g(U_i \times \{w_n\}) \cap U}$ is a diffeomorphism onto its image, denote the inverse of this map as $p^{-1}.$  Let $(\widetilde{w})$ and $(\widetilde{z})$ be coordinates around $q$ given by the Lemma \ref{zw coordinates}. Suppose that in $(\widetilde{w})$, $g(U_i\times \{w_n\})\cap U$ can be expressed as 
\begin{equation}
\{\widetilde{w}: \widetilde{w}_n =f(\widetilde{w}')+\sqrt{-1} g(\widetilde{w}')\}.
\end{equation}
 for some real functions $f$ and $g$. For any $\delta>0$, we can let $\delta_0$ be small enough such that we can follow the proof of the Lemma \ref{gui estimate} to get that
\begin{equation*}
dist\{span<g_* \partial_{w^1},...,g_* \partial_{w^{n-1}}>, span<\partial_{\widetilde{w}^1},...,\partial_{\widetilde{w}^{n-1}}>\}\le C \delta .
\end{equation*}
Combining the above formula with the Lemma \ref{normal bundle metric}, we can get that:
\begin{equation}
\frac{|f'(0)|}{|\widetilde{w}^n| |\log |\widetilde{w}^n||} \le C (\delta + e^{-\eta t}),\,\,\, \frac{|g'(0)|}{|\widetilde{w}^n| |\log |\widetilde{w}^n||} \le C  (\delta + e^{-\eta t}).
\end{equation}
for some constant $C$. Then we can get that:
\begin{equation*}
|(p^{-1})^* g_{\omega} -g_D| \le C(e^{-\eta t}+\delta) g_D.
\end{equation*}
Then we can let $\delta_0$ be small such that $t$ is big and $\delta$ is small to get that:
\begin{equation*}
(1-\epsilon) g_D \le (p^{-1})^* g_{\omega} \le (1+\epsilon) g_D.
\end{equation*}
This concludes the proof of this Lemma. 
\end{proof}

\begin{cor}\label{cor 4.19}
There exists a constant $C$ and $\delta_0>0$ such that for any $g\in Iso_0^D(X,\omega)$, for any $i$ and $|w_n| \le \delta_0$, we have that: 
\begin{equation*}
|\nabla g_{U_i,w_n}|_{\omega_D}\le C, \,\,\, |\nabla g_{U_i,w_n}^{-1}|_{\omega_D}\le C.
\end{equation*}
\end{cor}
\begin{proof}
This Corrolary directly follows from the Lemma \ref{lem 4.18}.
\end{proof}

\begin{lem}\label{almost holo}
For any $\epsilon>0$, there exists $\delta_0>0$ such that for any $w_n\in B_{\delta_0}(0)$  and $g \in Iso_0^D(X,\omega)$, we have that
\begin{equation*}
|\bar \partial g_{U_i,w_n}|\le \epsilon.
\end{equation*}
\end{lem}
\begin{proof}
According to the Lemma \ref{lem 4.15}, we have that
\begin{equation*}
    D g \circ J_0 (x)  - J_0 \circ D g(x)= O(e^{- \eta \min\{t_1,t_2\}}). 
\end{equation*}
Here $t_1$ is the value of the function $t$ at $x$ and $t_2$ is the value of the function $t$ at $g(x)$. Apply $p_*$ to the above formula. We can get that:
\begin{equation}\label{e 4.12}
p_* D g \circ J_0 = p_* \circ J_0 \circ D g + O(e^{-\eta\min\{t_1,t_2\}}) = J_0 \circ p_* \circ D g+ O(e^{-\eta\min\{t_1,t_2\}}).
\end{equation}
Here the second equality above uses the fact that $J_0$ is the product almost complex structure on $U_i \times \Delta^*$. As $w_{i,k}$ goes to zero, we have that $t_1$ and $t_2$ goes to $\infty$, according to the Lemma \ref{image goes to d}. This concludes the proof of this Lemma.
\end{proof}

\begin{lem}\label{gigj close}
There exists a constant $C$ independent of $g$ such that for any $g\in Iso_0^D(X,\omega)$, we have that:
\begin{equation*}
d(g_{U_i,w_n}|_{U_i \cap U_j}, g_{U_j,w_n}|_{U_i \cap U_j}) \le \frac{C}{|\log |w_n||},
\end{equation*}
for any $i,j \in \mathcal{A}$.
\end{lem}
\begin{proof}
 for any $i,j \in \mathcal{A}$ and for any $w'\in U_i \cap U_j$, there exists $\theta \in S^1$ such that  $\sigma_i(w') =e^{i\theta} \sigma_j(w')$. Using the Lemma \ref{normal bundle metric}, we have that 
 \begin{equation*}
 d(exp_{w'}(w_n \sigma_i(w')), exp_{w'}(w_n \sigma_j(w')))= d(exp_{w'}(w_ne^{i\theta} \sigma_j (w')), exp_{w'}(w_n \sigma_j (w')))  \le C \frac{1}{|\log|w_n||}.
 \end{equation*}
Here $d$ is the distance function induced by the Poincar\'e type metric $\omega$. Since $g$ preserves $\omega$, we have that 
\begin{equation*}
d(g(exp_{w'}(w_n \sigma_i)), g(exp_{w'}(w_n \sigma_j))) \le C \frac{1}{|\log|w_n||}.
\end{equation*}
Since the projection map $p$ satisfies that:
\begin{equation*}
|\nabla p|\le 2.
\end{equation*}
Then we can get that:
\begin{equation*}
d(p \circ g(exp_{w'}(w_n \sigma_i)), p \circ g(exp_{w'}(w_n \sigma_j))) \le C \frac{1}{|\log|w_n||}.
\end{equation*}
This concludes the proof of this Lemma.
\end{proof}
 
Note that we assume that $Aut_0(D)=\{Id\}$, which implies that $Aut(D)$ is discrete. However, $Aut(D)$ may not be $\{Id\}$ or even finite. As a result, when we take a sequence of $g_k \in Iso_0^D(X,\omega)$ which converges locally uniformly to some map $g$ on $X\setminus D$, even if we can prove that $g$ can be extended to $D$, we still need more work to prove that $g|_{D}=\{Id\}$. 
\begin{lem}\label{finite iso}
Assume that $Aut_0(D)=\{Id\}$. Then $Iso(D,\omega_D)$ is a finite set.
\end{lem}
\begin{proof}
Since $Aut_0(D)=\{Id\}$, we have that $Aut(D)$ is a discrete set. Note that $Iso(D,\omega_D)$ is a compact set in $Aut(D)$. Then we have that $Iso(D,\omega_D)$ is a finite set.
\end{proof}

\begin{lem}\label{lem 4.22}
For any $\epsilon>0$, there exists $\delta>0$ such that for any family of maps $\{g_i\}_{i \in \mathcal{A}}$, where $g_i$ is a map from $U_i$ to $D$ satisfying:
\begin{enumerate}
\item $(1-\delta) g_D \le g_i^* g_D \le (1+\delta)g_D$ for any $i$.
\item $ |\bar \partial g_i| \le \delta$ for any $i$.
\item $d(g_i |_{U_i \cap U_j}, g_j |_{U_i \cap U_j})\le \delta$  for any $i,j$.
\end{enumerate}
there exists some $g \in Iso(D,\omega_D)$ such that 
\begin{equation*}
|g_i - g |\le \epsilon.
\end{equation*}
\end{lem}
\begin{proof}
We prove the Lemma by contradiction. Suppose that we have a sequence of maps $\{g_{i,k}\}$, for $ i \in \mathcal{A}$ and $k \ge 1$. For each $k$, we have that:
\begin{enumerate}
\item 
\begin{equation}\label{almost isometry}
(1-\frac{1}{k}) g_D \le g_{i,k}^* g_D \le (1+\frac{1}{k})g_D
\end{equation}
 for any $i$.
\item $ |\bar \partial g_{i,k}| \le \frac{1}{k}$ for any $i$.
\item $d(g_{i,k} |_{U_i \cap U_j}, g_{j,k} |_{U_i \cap U_j})\le \frac{1}{k}$  for any $i,j$.
\end{enumerate}
and 
\begin{equation*}
\inf_{g \in Iso(D,\omega_D)} \sup_{1 \le i \le N}|g_{i,k}- g| \ge \epsilon
\end{equation*}
for some $\epsilon >0$ independent of $k$.  Using $(1-\frac{1}{k}) g_D \le g_{i,k}^* g_D \le (1+\frac{1}{k})g_D$, there exists a constant $C$ such that 
\begin{equation}\label{e 4.14}
|\nabla g_{i,k}|\le C, \,\,\, |\nabla g_{i,k}^{-1}|\le C
\end{equation}
for any $i,k$. Then we can use the Arzela-Ascoli theorem to get a subsequence of $g_{i,k}$ (still denoted as $g_{i,k}$) such that there exists $g_{i}$ such that:
\begin{equation*}
\lim_{k \rightarrow \infty} ||g_{i,k}-g_i ||_{C^{\alpha}}=0.
\end{equation*}
Combining this and (\ref{almost isometry}), for any $q\in U_i$, there exists a neighbourhood $B_{\epsilon}(q)$ such that $g_i |_{B_{\epsilon}(q)}$ preserves the distance induced by $g_D$. In fact, we can let $\epsilon>0$ be small enough. Then for any $q_1,q_2 \in B_{\epsilon}(q)$, we have that:
\begin{equation*}
\begin{split}
d(g_i(q_1),g_i(q_2)) = \lim_{k \rightarrow \infty} d(g_{i,k}(q_1), g_{i,k}(q_2)) = \lim_{k \rightarrow \infty}d_{i,k}(q_1, q_2) &\le \lim_{k \rightarrow \infty} (1+\frac{1}{k}) d (q_1,q_2) \\
& = d(q_1,q_2)
\end{split}
\end{equation*}
and
\begin{equation*}
\begin{split}
d(g_i(q_1),g_i(q_2)) = \lim_{k \rightarrow \infty} d(g_{i,k}(q_1), g_{i,k}(q_2)) = \lim_{k \rightarrow \infty}d_{i,k}(q_1, q_2)&\ge \lim_{k \rightarrow \infty} (1-\frac{1}{k}) d (q_1,q_2) \\
& = d(q_1,q_2).
\end{split}
\end{equation*}
Here $d$ is the distance function induced by $g_D$ and $d_{i,k}$ is the distance function induced by $g_{i,k}^* g_D$. Then we have shown that $g$ is distance preserving on $B_{\epsilon}(q)$. Using the fact that $ |\bar \partial g_{i,k}| \le \frac{1}{k}$, we get that $g_i$ is weakly holomorphic. This implies that $g_i$ is indeed holomorphic and smooth, using standard elliptic regularity results. Since a smooth distance preserving map is also metric preserving which can be shown using normal coordinates around $q$ and $g_i(q)$, we have that: 
\begin{equation*}
g_i^* g_D = g_D.
\end{equation*}
In conclusion, we have that:
\begin{enumerate}
\item $ g_{i}^* g_D = g_D$ for any $i$.
\item $ \bar \partial g_{i}=0$ for any $i$.
\item $g_i |_{U_i \cap U_j} = g_j |_{U_i \cap U_j}$  for any $i,j$.
\end{enumerate}
and 
\begin{equation}\label{e 4.15}
\inf_{g \in Iso(D,\omega_D)} \sup_{1 \le i \le N}|g_{i}- g| \ge \epsilon
\end{equation}
However, we can use the fact that $g_i |_{U_i \cap U_j} = g_j |_{U_i \cap U_j}$  for any $i,j$ to define $\widetilde{g}$ by $\widetilde{g}= g_i $ on $U_i$. Using the Lemma \ref{holo to auto} below, we have that $\widetilde{g}\in Iso(D,\omega_D)$. This is a contradiction with (\ref{e 4.15}).
\end{proof}

\subsubsection{holomorphic maps from $D$ to $D$ that preserve $g_D$}
In this section, we want to prove the following Lemma:
\begin{lem}\label{holo to auto}
Let $g_D$ be a metric on $D$. Suppose that $g$ is a holomorphic map from $D$ to $D$ itself and $g^* g_D= g_D$. Then we have that $g$ is an automorphism of $D$. 
\end{lem}
\begin{proof}
Since $g^* g_D= g_D$, we have that $g$ is a local diffeomorphism. As a result, the image of $g$ is open and closed and nonempty. Thus $g$ is surjective. A local diffeomorphism that is also surjective must be a covering map. In order to prove that $g$ is an automorphism of $D$, it suffices to prove that $g$ is of degree 1. Using the definition of the degree of a map, we have that:
\begin{equation*}
deg (g) \int_D dvol_D = \int_D g^* dvol_D= \int_D dvol_D.
\end{equation*}
This implies that $deg (g) =1$. This concludes the proof of this Lemma.
\end{proof}

Now we are ready to prove the main Propositions in this subsection:
\begin{proof}
(of the Proposition \ref{gk converge}). 
Using the Lemma \ref{finite iso}, there are only finite elements in $Iso(D,\omega_D)$ which we denote as $\{\widetilde{g}_i\}_{i=1}^N$. There exists a constant $\epsilon_0 >0$ such that for any $i \neq j$, we have that
\begin{equation}\label{e 4.16}
d(\widetilde{g}_i,\widetilde{g}_j) \ge \epsilon_0.
\end{equation}
Take $\epsilon = \frac{\epsilon_0}{3}$. Then, we can use the Lemma \ref{lem 4.22} to get a $\delta$ with respect to this $\epsilon$. Then we fix this $\delta$. According to the Lemma \ref{lem 4.18}, Lemma \ref{gigj close} and the Lemma \ref{almost holo}, there exists $\delta_0>0$ independent of $g$ such that  the assumptions of the Lemma \ref{lem 4.22} hold with $g_i =g_{U_i, w_n}$ for any $g\in Iso_0^D(X,\omega)$, $\delta$ that we fix before, and $|w_n|\le \delta_0$. Then the Lemma \ref{lem 4.22} implies that there exists $\widetilde{g}_{w_n} \in Iso(D,\omega_D)$ such that 
\begin{equation}\label{e 4.17}
|g_{U_i,w_n}- \widetilde{g}_{w_n}|\le \frac{\epsilon_0}{3}.
\end{equation}
This map $\widetilde{g}_{w_n}$ is uniquely determined by $w_n$ since (\ref{e 4.16}).\\
 We \textbf{claim} that $ \widetilde{g}_{w_n}= Id$ for any $|w_n| \le \delta_0$. In fact, suppose that there exists $|w_n'| \le \delta$ such that 
\begin{equation}\label{e 4.18}
\widetilde{g}_{w_n'}\neq Id.
\end{equation}
 We can consider $w_n(s)= s w_n'$. Since $g|_D =Id$, we have that $g_{U_i,0}=Id$ which implies that $\widetilde{g}_{w_n(0)}=Id$. Combining this with (\ref{e 4.17}) and (\ref{e 4.18}), we have that there exists $s$ and two different elements $\widetilde{g}_1, \widetilde{g}_2 \in Iso(D,\omega_D)$ such that 
\begin{equation*}
|g_{U_i,w_n(s)}- \widetilde{g}_i |\le \frac{\epsilon_0}{3},
\end{equation*}
for any $i=1,2$. This implies that
\begin{equation*}
|\widetilde{g}_i -\widetilde{g}_j| \le \frac{2\epsilon_0}{3}. 
\end{equation*}
This is a contradiction with (\ref{e 4.16}). This concludes the proof of the claim. Then the lemma follows immediately from this claim.
\end{proof}

\section{Characterization of isometry group}
In this section we want to prove the Theorem \ref{compact isometry}. We follow \cite{C4} to prove this proposition. The main obstacle we come across in the Poincar\'e type case compared with the smooth case in \cite{C4} is that we need to prove that $Iso_0^D(X,\omega)$ is a compact group, which we have proved in the Theorem \ref{compactness of isometry group}. The rest of the proof is similar to that in \cite{C4}. As a result, we will sketch the proof and emphasize the modifications we make in this section.

First, we decompose $Ker L$ which is seen as a $\mathbb{C}-$module of complex-valued functions.

\begin{defn}\label{defn 7.2}
\begin{enumerate}
\item Denote $E_{\lambda}$ as the eigenspace of $\overline{L}$ over  $Ker L$ for each $\lambda$ in the spectrum. 
\item Denote $E_{0,r}$ as the real functions in $E_0$. Denote $E_{0,i}$ as the purely imaginary functions in $E_0$. 
\item Define $\mathbf{h}_{\lambda}=\nabla^{1,0} E_{\lambda}$ for $\lambda > 0$.
\item Define $\mathbf{h}_0' =\nabla^{1,0} E_0$, $\mathbf{h}_0 = \mathbf{a}^D_{\parallelsum}(X) \oplus \mathbf{h}_0'$. Recall that $\mathbf{a}^D_{\parallelsum}(X)$ consists of auto parallel holomorphic vector fields in $\mathbf{h}_{\parallelsum}^D$.
\item Define $\mathfrak{l}'=\nabla^{1,0}E_{0,i}$, $\mathbf{m}=\nabla^{1,0}E_{0,r}$ and $\mathfrak{l}=\mathbf{a}^D_{\parallelsum}(X)\oplus \mathfrak{l}'$.
\end{enumerate}
\end{defn}

\begin{lem}\label{lem 7.3}
Suppose that $\omega$ is an Poincar\'e type extremal K\"ahler metric. Then the corresponding Lichnerowicz operator satisfies that:
\begin{equation*}
Ker(L |_{W^{4,2}_{0,\mathbb{C}}} )= \oplus \Sigma_{\lambda \in Spec(\bar L|_{Ker (L)})}E_{\lambda}= E_{0,r}\oplus E_{0,i}\oplus \Sigma_{\lambda \in Spec(\bar L|_{ Ker(L)}), \lambda>0}E_{\lambda}
\end{equation*}
\end{lem}
\begin{proof}
Since $\omega$ is an extremal K\"ahler metric, we have that $L$ and $\bar L$ commute with each other. Then we have that for any $v\in Ker L|_{W^{4,2}_{0,\mathbb{C}}}$, $L \bar L v= \bar L L v=0$. This implies that $\bar L$ can be seen as an operator on $ Ker (L|_{W^{4,2}_{0,\mathbb{C}}})$.  Using the first part of the proof of the Lemma \ref{ker l}, we have that $\nabla^{1,0}  L|_{W^{4,2}_{0,\mathbb{C}}} \subset h^D_{\parallelsum}$. Since $h^D_{\parallelsum}$ is of finite dimension, so is $Ker (L|_{W^{4,2}_{0,\mathbb{C}}})$. So we have that $\bar L$ is a self-adjoint operator on a finite dimensional space $Ker (L|_{W^{4,2}_{0,\mathbb{C}}})$. Then we can decompose $Ker (L|_{W^{4,2}_{0,\mathbb{C}}})$ using the eigenspaces of $\bar L$. This implies that 
\begin{equation*}
Ker (L|_{W^{4,2}_{0,\mathbb{C}}}) = \oplus \Sigma_{\lambda \in Spec(\bar L|_{Ker(L)})}E_{\lambda}.
\end{equation*}
For any $f=g+hi \in Ker L \cap Ker \bar L$, we have that $\bar L (g+hi)=0$. Take the conjugation of this formula, we have that $L g - L h i=0 $. Note that $L (g+hi)=0$. So we have that $Lg=Lh =0$. This implies that $E_0=E_{0,r} \oplus E_{0,i}$.
\end{proof}

\begin{lem}\label{lie algebra}
Let $\omega$ be a Poincar\'e type extremal K\"ahler metric on $X$. Consider the special element $X_0=\nabla^{1,0}R \in \mathbf{h}^D_{\parallelsum}$. Then we have the following relations:
\begin{enumerate}
\item For each $\lambda \in spec(\bar L |_{Ker L})$ including $\lambda=0$, and for each $Y\in \mathbf{h}_{\lambda},$
\begin{equation*}
[X_0,Y]=\lambda Y.
\end{equation*}
\item For each pair of numbers $\lambda,\mu$ in $spec(\bar L |_{Ker L})$, we have $[\mathbf{h}_{\lambda}, \mathbf{h}_{\mu}] \subset \mathbf{h}_{\lambda+\mu},$ with the usual convention that $\mathbf{h}_{\lambda+\mu}=\{0\},$ if $\lambda+\mu$ is not in the spectrum.
\item The subspace $\mathbf{a}^D_{\parallelsum}(X)$, $\mathfrak{l}'$ and $\mathbf{m}$ satisfies the relation:
       \begin{enumerate}
       \item $\mathbf{a}^D_{\parallelsum}(X)$ is in the center of $\mathbf{h}_0$.
       \item $[\mathfrak{l},\mathfrak{l}] \subset \mathfrak{l}, [\mathfrak{l},\mathbf{m}] \subset \mathbf{m}, [\mathbf{m},\mathbf{m}] \subset \mathfrak{l}$.
       \end{enumerate}
\end{enumerate}
\end{lem}
\begin{proof}
The proof of this Lemma follows the proof of Lemma 3.2 in \cite{C4} word by word.
\end{proof}

We also have the following Lemma:
\begin{lem}\label{lem 7.5}
Let $\omega$ be a Poincar\'e type K\"ahler metric. Suppose that $h \in W^{1,2}_{0,\mathbb{C}}$. Then we have that $Re (h)$ is a constant if and only if  $\nabla^{1,0}h$ is a Killing vector field.
\end{lem}
\begin{proof}
Denote $h=a +\sqrt{-1}f$.  Denote $v=\nabla^{1,0}h$. Then we have that 
\begin{equation*}
\begin{split}
&2 \mathcal{L}_{Re v} \omega = \mathcal{L}_v \omega + \mathcal{L}_{\bar v} \omega = d (\iota_v  \omega)+ d(\iota_{\bar v}  \omega) \\
&= d(\sqrt{-1}g^{i\bar j}f_{\bar j}g_{i\bar k}\sqrt{-1} d\bar z^k +\sqrt{-1}g^{\bar i j}f_{j} g_{k \bar i}\sqrt{-1} d z^k+ g^{i\bar j}a_{\bar j}g_{i\bar k}\sqrt{-1} d\bar z^k-g^{\bar i j}a_{j} g_{k \bar i}\sqrt{-1} d z^k)\\
&= -d(f_{\bar j}d \bar z^j +f_k d z^k) +dd^c a=-d^2 f+dd^c a=dd^c a.
\end{split}
\end{equation*}
So $v$ is a Killing vector field if and only if $ \mathcal{L}_{Re v} \omega=0$ if and only if $dd^c a =0$ which is equivalent to $a=C$ for some constant $C$. 
\end{proof}

In our previous paper \cite{XZ}, we proved the following decomposition of holomorphic vector fields:
\begin{prop}\label{reductivity}
Let $\omega$ be a Poincar\'e type extremal K\"ahler metric. 
 One can define in terms of $\omega$ a unique semidirect sum splitting of the Lie algebra $\mathbf{h}^D_{\parallelsum}$:
\begin{equation*}
\mathbf{h}^D_{\parallelsum}=\mathbf{a}^D_{\parallelsum}(M) \oplus \mathbf{h}^{D}_{\parallelsum,\mathbb{C}}.
\end{equation*}
\end{prop}

We also need the following lemma:
\begin{lem}\label{lem 7.6}
Let $\omega$ be a Poincar\'e type extremal K\"ahler metric. Let $\mathfrak{l}$, $\mathbf{a}^D_{\parallelsum}$ and $\mathfrak{l}'$ be defined in the Definition \ref{defn 7.2}. Then we have that:
\begin{equation*}
\{v\in \mathbf{h}^D_{\parallelsum}:v \text{ is a Killing vector field with respect to } \omega \}=\mathfrak{l}= \mathbf{a}^D_{\parallelsum} \oplus \mathfrak{l}' .
\end{equation*}
\end{lem}
\begin{proof}
First, we show that for any $v\in \mathbf{a}^D_{\parallelsum} \oplus \mathfrak{l}'$, $v$ is a Killing vector field. In any normal coordinate, for any $v=v^{k}\frac{\partial}{\partial z^k} \in \mathbf{a}^D_{\parallelsum}$ we have that $v^k_{,\bar l}= v^k_{,l}=0$. This implies that $\mathcal{L}_{Re v} \omega= 1/2 (\mathcal{L}_v \omega + \mathcal{L}_{\bar v} \omega)=1/2( d (\iota_v  \omega)+d (\iota_{\bar v}  \omega))=0$. So we have that $\mathbf{a}^D_{\parallelsum} \subset \{v\in \mathbf{h}^D_{\parallelsum}:v \text{ is a Killing vector field with respect to } \omega \}$. For any $v\in  \mathfrak{l}'$, $v=\nabla^{1,0} \sqrt{-1}f$ for some real-valued function $f$. Since we use the normal coordinate, we assume that $\omega= \Sigma_{i=1}^n \sqrt{-1} dz^i \wedge d \bar z^i$ at the given point. Then we have that 
\begin{equation*}
\begin{split}
2 \mathcal{L}_{Re v} \omega &= \mathcal{L}_v \omega + \mathcal{L}_{\bar v} \omega = d (\iota_v  \omega)+ d(\iota_{\bar v}  \omega) \\
&= d(\sqrt{-1}g^{i\bar j}f_{\bar j}g_{i\bar k}\sqrt{-1} d\bar z^k +\sqrt{-1}g^{\bar i j}f_{j} g_{k \bar i}\sqrt{-1} d z^k)\\
&= -d(f_{\bar j}d \bar z^j +f_k d z^k)=-d^2 f=0.
\end{split}
\end{equation*}
This concludes the proof of $\mathfrak{l} \subset \{v\in \mathbf{h}^D_{\parallelsum}:v \text{ is a Killing vector field with respect to } \omega \}$. 

On the other hand, for any $v\in  \{v\in \mathbf{h}^D_{\parallelsum}:v \text{ is a Killing vector field with respect to } \omega \},$ we can use the proposition \ref{reductivity}, the Lemma \ref{ker l} and the Lemma \ref{lem 7.3} to get that $v=H +\Sigma_{\lambda \ge 0} \nabla^{1,0} Y_{\lambda}$ and $Y_0=Y_0'+Y_0''$, with $H \in \mathbf{a}^D_{\parallelsum}, \nabla^{1,0}Y_{\lambda} \in \mathbf{h}_{\lambda}$ for each $\lambda \in spec(\bar L |_{ker L|_{W_{0,\mathbb{C}}^{4,2}}})$, $\nabla^{1,0}Y_0' \in \mathfrak{l}'$, $ \nabla^{1,0}Y_0''\in \mathbf{m}$. Since we have shown that 
\begin{equation*}
\mathbf{a}^D_{\parallelsum} \oplus \mathfrak{l}' \subset \{v\in \mathbf{h}^D_{\parallelsum}:v \text{ is a Killing vector field with respect to } \omega \},
\end{equation*} 
we can replace $v$ by  $v-(H+ \nabla^{1,0}Y_0')$ and assume that $v = \nabla^{1,0} Y_0'' +\Sigma_{\lambda > 0} \nabla^{1,0} Y_{\lambda}$. We first show that $Y_{\lambda}=0$ for any $\lambda>0$. We can write that $\Sigma_{\lambda > 0} Y_{\lambda} = \Sigma_i a_i f_i$ where $f_i$ is an eigenfunction with respect to the eigenvalue $\lambda_i>0$ and $f_i$ are linearly independent from each other. We have that $L(Y_0'' + \Sigma_{\lambda > 0} Y_{\lambda})=0$ and 
\begin{equation}\label{equ 7.1}
\bar L(Y_0'' + \Sigma_{\lambda > 0} Y_{\lambda}) = \Sigma_i a_i \lambda_i f_i.
\end{equation}
 Using the Lemma \ref{lem 7.5}, we have that the real part of $Y_0'' + \Sigma_{\lambda > 0} Y_{\lambda}$ is a constant. Without loss of generality, we assume that  $Y_0'' + \Sigma_{\lambda > 0} Y_{\lambda}$ is a purely imaginary function. Then we can take the conjugate of (\ref{equ 7.1}) to get that 
 \begin{equation*}
 \Sigma_i \overline{a_i} \lambda_i \bar f_i = L( \overline{Y_0'' + \Sigma_{\lambda > 0} Y_{\lambda}})=-L (Y_0'' + \Sigma_{\lambda > 0} Y_{\lambda})=0.
 \end{equation*}
 Since $f_i$ are linearly independent of each other and $\lambda_i>0$, we have that $a_i=0$ for each $i$. So we prove that $Y_{\lambda}=0$ for any $\lambda>0$ and thus $v=\nabla^{1,0} Y_0''$. Using the Lemma \ref{lem 7.5} we have that $v=0$. This concludes the proof of the Lemma. 
\end{proof}

Before we prove the Theorem \ref{compact isometry}, we need to prove the following Lemma:

\begin{lem}\label{liegroup noncompact}
For any nonzero vector field $v= \nabla^{1,0}u \in \mathbf{m}$ for a real function $u$, we have that the Lie group $\{exp (tv): t \in \bR\}$ is not contained in any compact subgroup of $Aut_0^D(X)$.
\end{lem}
\begin{proof}
We prove this Lemma by contradiction. Denote $g_t = exp(tv)$ for $t \in \bR$. Suppose that $\{g_t: t \in \bR\}$ is contained in a compact subgroup of $Aut_0^D(X)$. Then we can take a sequence $\{g_{n_k}\}_{k \in \mathbb{N}}$ with $\lim_{k \rightarrow \infty}n_k= \infty$ such that $g_{n_k}$ converge to $g$ for some $g\in Aut_0^D(X)$. For any $t\in \bR$, we have that
\begin{equation*}
g_t \circ g_{n_k}= g_{t+ n_k}= g_{n_k} \circ g_t.
\end{equation*}
Let $k \rightarrow \infty$, we get that
\begin{equation}\label{ggt commute}
g_t \circ g = g \circ g_t.
\end{equation}
First, we want to study the behavior of the orbit of $v$ starting from an arbitrary point $p\in X$, denoted as $\Sigma$, and the behaviour of $u$ along $\Sigma$. If $v(p)=0$, then the orbit consists of only one point $p$. If $p\in X \setminus D$ and $v(p) \neq 0$, then we have that $\Sigma$ doesn't intersect with $D$. In fact, $v|_D$ is a vector field parallel to $D$. So if $\Sigma$ intersects with $D$, then the whole orbit lies in $D$. Since $\omega$ is smooth on $X\setminus D$, we can use $v= \nabla^{1,0}u$ to get that $u$ is strictly increasing with respect to $t$ on $\Sigma$. If $p \in D$ and $v(p) \neq 0$, then we again use the fact that $v|_D$ is a vector field parallel to $D$ to get that $\Sigma$ lies in $D$. According to the Lemma \ref{lem 7.21}, we have that $u \in \widetilde{C}_{-\eta_0}^{5,\alpha}$ and $\nabla_{\omega_D}^{1,0}u |_D =v|_D$, where $\omega_D$ is the metric on $D$ defined by the Lemma \ref{asym}. This implies that $u$ is strictly increasing with respect to $t$ on $\Sigma$ as well. 

Now we fix a point $p\in X$ such that $v(p) \neq 0$. We consider the following two cases:

(1) $v(g(p))=0$. From the above argument, we can see that for any $t\in \bR$, $g_t \circ g (p)=g(p)$.  Combining this with (\ref{ggt commute}), we can get that
\begin{equation*}
g(g_t(p))=g_t(g(p))=g(p).
\end{equation*}
However, since $v(p)\neq 0$, we have that $g_t(p) \neq p$. This contradicts with the fact that $g\in Aut_0(X)$ which implies that $g$ is injective.

(2) $v(g(p)) \neq 0$. Then from the above argument we have that for $t>0$, 
\begin{equation}\label{ugtug}
u(g_t \circ g (p)) > u(g(p)).
\end{equation}
 However, we have that
\begin{equation*}
g_t \circ g(p)= \lim_{k\rightarrow \infty} g_{t+ n_k} (p),\,\,\, g(p)= \lim_{k \rightarrow \infty} g_{n_k} (p).
\end{equation*}
Since $n_k \rightarrow \infty$, we can find a subsequence of $\{t+n_k\}$, denoted as $\{a_k\}$, and a subsequence of $n_k$, denoted as $\{b_k\}$ such that $a_k < b_k < a_{k+1}$ for any $k$. Since $u$ is strictly increasing along $\Sigma$, we have that
\begin{equation*}
u(g_{a_k}(p) )< u( g_{b_k} (p) )< u(g_{a_{k+1}}(p)).
\end{equation*} 
Here we use the fact that $u\in C^0(X)$ because we have $u \in \widetilde{C}_{-\eta_0}^{5,\alpha}$ according to the Lemma \ref{lem 7.21}. Let $k \rightarrow \infty$ in the above formula. We get that 
\begin{equation*}
u(g_t \circ g(p)) =u (g(p)).
\end{equation*}
This contradicts with (\ref{ugtug}). This concludes the proof of this Lemma.
\end{proof}

Now we can prove the Theorem \ref{compact isometry}.
\begin{proof}
(of the Theorem \ref{compact isometry}.)
We prove the proposition by contradiction. Suppose that $Iso_0^D(X,\omega)$ is not a maximal compact, connected subgroup of $Aut_0^D(X,\omega)$. Then there exists a compact, connected subgroup $G \subset Aut_0^D(X,\omega)$ that properly contains $Iso_0^D(X,\omega)$. Let $Y$ be an element of the Lie algebra $\widetilde{\mathfrak{l}}$ of $G$ that is not in $\mathfrak{l}$. Denote:
\begin{equation*}
Y=H +\Sigma_{\lambda \ge 0} \nabla^{1,0} Y_{\lambda} \text{ and } Y_0=Y_0' +Y_0'',
\end{equation*}
with $Y_0' \in E_{0,i}$, $Y_0'' \in E_{0,r}$, $H \in \mathbf{a}^D_{\parallelsum}, \nabla^{1,0}Y_{\lambda}\in \mathbf{h}_{\lambda}$ for each $\lambda \in spec(ad X_0)$. Denote $Z_0 =\nabla^{1,0} (R \sqrt{-1}) \in \mathfrak{l}'$, we can consider the adjoint action of the one-parameter group of isometries generated by $Z_0$ on $Y$. We then have 
\begin{equation*}
ad\, exp(t Z_0) (Y) = H + \nabla^{1,0} Y_0' + \nabla^{1,0} Y_0'' + \Sigma_{\lambda>0} e^{\lambda t \sqrt{-1}} \nabla^{1,0}Y_{\lambda} \in \widetilde{\mathfrak{l}}.
\end{equation*}
Then, we can take appropriate linear combinations of the resulting elements for sufficiently many values of $t$ to get that:
\begin{equation*}
H + \nabla^{1,0} Y_0' + \nabla^{1,0}Y_0'' \in \widetilde{\mathfrak{l}},
\end{equation*}
and
\begin{equation*}
\nabla^{1,0}Y_{\lambda} \in  \widetilde{\mathfrak{l}} 
\end{equation*}
for each $\lambda>0$. If, for some $\lambda>0$, $Y_{\lambda} \neq 0$, then we have that $L \neq \bar L$. Using the Lemma \ref{lie algebra}, we have that $Z_0$ and $\nabla^{1,0} Y_{\lambda}$ generate a solvable, non-abelian Lie subalgebra of $\widetilde{\mathfrak{l}}$. This is impossible since $\widetilde{\mathfrak{l}}$ generates a compact group. As a result, we have that $\Sigma_{\lambda>0}\nabla^{1,0} Y_{\lambda}=0$ and
\begin{equation*}
Y= H + \nabla^{1,0}Y_0' +\nabla^{1,0}Y_0'' \in \widetilde{\mathfrak{l}}.
\end{equation*}
Since $Y \notin \mathfrak{l}$, we have that $\nabla^{1,0}Y_0'' \neq 0$. Note that we have that $\nabla^{1,0}Y_0'' \in \mathfrak{l} \subset \widetilde{\mathfrak{l}}$. By definition, $Y_0'' \in Ker L \cap Ker \bar L$ is a real-valued function. Then we can use the Lemma \ref{liegroup noncompact} to get that $Y=0$. This concludes the proof of the Theorem  \ref{compact isometry}. The above proof essentially follows \cite{C4} and \cite{L}.
\end{proof}

\section{Extremal K\"ahler vector field}
In this section, we want to prove the following proposition:
\begin{prop}\label{unique extremal vector field}
Let $\omega_i \in [\omega]$ be two Poincar\'e type extremal K\"ahler metrics, such that 
\begin{equation*}
Iso_0^D(M,\omega_1)= Iso_0^D(M,\omega_2).
\end{equation*}
Then we have $\nabla_{\omega_1}^{1,0}(R_{\omega_1})= \nabla_{\omega_2}^{1,0}(R_{\omega_2})$.
\end{prop}
The above proposition is an adaptation of a result due to Futaki-Mabuchi \cite{FM} to the Poincar\'e type case. One can see Berman-Berndtsson\cite{BB}. for a detailed formulation. Note that Auvray \cite{A4} defined the Poincar\'e type Futaki character. We will use Berman-Berndtsson's formulation to sketch the proof of the Proposition \ref{unique extremal vector field} for the convenience of readers. 

For any $V\in \mathbf{h}^D_{\parallelsum, \mathbb{C}}$, and any Poincar\'e type K\"ahler metric $\omega$, there exists a function $h$ with $\int h \omega^n=0$ such that $V=\nabla^{1,0}h$. We can define $h^V_{\omega}$ to be $h$. 

\begin{lem}\label{V different metric}
If $\omega_u= \omega_X + i \partial \bar \partial u$ is a Poincar\'e type K\"ahler metric and $\omega_X$ is a smooth K\"ahler metric on $X$, then for any $V \in \mathbf{h}^D_{\parallelsum, \mathbb{C}}$, we have that:
\begin{equation*}
h^V_{\omega_u}= h^V_{\omega_X} + V(u).
\end{equation*}
\end{lem}
\begin{proof}
Since we have that 
\begin{equation*}
i \bar \partial (h^V_{\omega_X}+V(u))= \iota_V \omega_u. 
\end{equation*}
As a result, we have that 
\begin{equation*}
h^V_{\omega_u}= h^V_{\omega_X} + V(u)+ c(u),
\end{equation*}
where $c(u)$ is a constant on $X$. We can calculate 
\begin{equation}\label{cdot}
0= (\frac{d}{dt}) \int_X h^V_{\omega_{tu}}\omega_{tu}^n = \int_X (V(u)+ \dot{c}(tu))\omega_{tu}^n + n \int_X h^V_{\omega_{tu}} i \partial \bar \partial u \wedge \omega_{tu}^{n-1}.
\end{equation}
Since $u$ is a Poincar\'e type K\"ahler potential, it satisfies that:
\begin{equation*}
|\nabla_{\omega_{tu}} u|_{\omega_{tu}}, |\nabla^2_{\omega_{tu}} u|_{\omega_{tu}} < \infty.
\end{equation*}
Then we can use the Lemma \ref{G-S} to do integration by part in (\ref{cdot}) to get that $\dot{c}(tu)=0$, so we have that $c(u)=0$ since $c(0)=0$.
\end{proof}

Then, we can define a bilinear form on $\mathbf{h}^D_{\parallelsum, \mathbb{C}}$ by 
\begin{equation*}
<V,W>_{\omega}= \int_X h^V_{\omega} h^{W}_{\omega} \omega^n.
\end{equation*}
We can prove the following proposition:
\begin{prop}\label{prop 6.3}
$<,>_{\omega}$ only depends on the cohomology class $[\omega]$.
\end{prop}
\begin{proof}
We take a curve of metrics $\omega_t= \omega+ i \partial \bar \partial u_t$ in  $\mathcal{PM}_{[\omega]}$ and differentiate,
\begin{equation*}
(\frac{d}{dt})\int_X h^V_{\omega_t} h^W_{\omega_t} \omega_t^n = \int_X (V(\dot{u})h^W_{\omega_t} +W(\dot{u})h^V_{\omega_t})\omega_t^n +n \int_X h_{\omega_t}^V h_{\omega_t}^W i \partial \bar \partial \dot{u} \wedge \omega_t^{n-1}.
\end{equation*}
Then, we can use the Lemma \ref{G-S} to do integration by part to get that the above integral is equal to zero.
\end{proof}

For any compact subgroup $K$ of $Aut_0^D(X)$, we define 
\begin{equation*}
\mathbf{h}^D_{\parallelsum, \mathbb{C},K} \triangleq \{v\in \mathbf{h}^D_{\parallelsum}: \text{ the flow induced by }Im V \text{ lie in }K\}.
\end{equation*}
We can also prove the following proposition:
\begin{prop}\label{prop 6.4}
For any compact subgroup $K$ of $Aut_0^D(X)$ the restriction of $<,>$ to $\mathbf{h}^D_{\parallelsum, \mathbb{C},K}$ is real valued and positive definite, in particular non-degenerate.
\end{prop}
\begin{proof}
Taking averages of an arbitrary K\"ahler form, we can represent our form by a K-invariant K\"ahler form $\omega$ using the Proposition \ref{prop 6.3}. Then we can use the proof of the Lemma \ref{lem 7.6} to get that $h^V_{\omega}$ is real-valued if $V\in \mathbf{h}^D_{\parallelsum, \mathbb{C},K}$. Then this proposition follows immediately.
\end{proof}
Recall the Poincar\'e type Futaki character defined by Auvray: For any $Z \in \mathbf{h}^D_{\parallelsum, \mathbb{C}}$. we define
\begin{equation*}
\mathcal{F}^D_{[\omega_X]}(Z) = \int_{X\setminus D} R_{\omega_1} \frac{1}{2} Re h_{\omega_1}^Z \frac{\omega_1^n}{n !}.
\end{equation*}
Auvray \cite{A3} proved that the Poincar\'e type Futaki character does not depend on $\omega$ of class $[\omega_X]$, provided it is of Poincar\'e type:
\begin{lem}
Let $\widetilde{\omega}$ be any Poincar\'e type metric in $\mathcal{PM}^D_{[\omega_X]}$, and $Z \in \mathbf{h}^D_{\parallelsum, \mathbb{C}}$. Then $\mathcal{F}^D_{[\omega_X]}(Z)= \int_{M} R_{\widetilde{\omega}} \frac{1}{2} Re h^Z_{\widetilde{\omega}} \frac{\widetilde{\omega}^n}{n!}.$
\end{lem}

Now, we are ready to prove the main Proposition in this section.
\begin{proof}
(of the Proposition \ref{unique extremal vector field})
Denote $V_i= \nabla_{\omega_i}^{1,0}R_{\omega_i}$ for $i=1,2$. Denote $K = Iso_0^D(M,\omega_1)=Iso_0^D(M,\omega_2)$. Using the Lemma \ref{lem 7.6}, we have that $Im V_1$ and $Im V_2$ both generate transformations lying in $K$. Then for any $Z\in \mathbf{h}^D_{\parallelsum, \mathbb{C},K}$, we have that:
\begin{equation*}
\mathcal{F}^D_{[\omega_X]}(Z) = \int_{M} R_{\omega_1} \frac{1}{2} Re h^Z_{\omega_1} \frac{\omega_1^n}{n!} = <V_1,Z>
\end{equation*}
and
\begin{equation*}
\mathcal{F}^D_{[\omega_X]}(Z) =\int_{M} R_{\omega_2} \frac{1}{2}Re h^Z_{\omega_2} \frac{\omega_2^n}{n!} = <V_2,Z >.
\end{equation*}
Then we get that
\begin{equation*}
<V_1,Z>= <V_2,Z >.
\end{equation*}
Using the Proposition \ref{prop 6.4}, we can get that $V_1=V_2$.
\end{proof}

\section{Uniqueness of Poincar\'e type extremal K\"ahler metric}

In this section, we want to prove the Theorem \ref{main theorem} and the Theorem \ref{unique extk}. Let $\omega_{i}= \omega + dd^c \varphi_i$, $i=1,2$ be two Poincar\'e type extremal K\"ahler metrics. We can use the Theorem \ref{compact isometry} and the fact that any two maximal compact subgroups of a Lie group are conjugate to each other to get that there exists $g\in Aut_0^D(X)$ such that $Iso_0^D(X, g^* \omega_1)= Iso_0^D(X,\omega_2)$. Therefore, in the rest of the section, we can assume that $Iso_0^D(X, \omega_1)= Iso_0^D(X,\omega_2)$ by replacing $\omega_1$ with $g^* \omega_1$. Denote $K=Iso_0^D(X, \omega_1)= Iso_0^D(X,\omega_2)$. Denote
\begin{equation*}
C_{K,\delta}^{k,\alpha}= \{ \varphi \in C_{\delta}^{k,\alpha}: \varphi = \varphi \circ \sigma \text{ for any } \sigma \in K \}.
\end{equation*}
and
\begin{equation*}
\widetilde{C}^{k,\alpha}_{K,\delta} \triangleq \{\varphi \in \widetilde{C}^{k,\alpha}_{\delta}: \varphi = \varphi \circ \sigma \text{ for any } \sigma \in K\}.
\end{equation*}
Let $X_1 = \nabla_{\omega_1}^{1,0} R_{\omega_1}$ be the holomorphic vector field corresponding to the metric $\omega_{1}$. From now on, we use $\omega_1$ as the background metric. We denote $$\omega_{\varphi} \triangleq \omega_1 + dd^c \varphi.$$

\subsection{Difference between Poincar\'e type extremal K\"ahler metrics}
Using the Lemma \ref{asym}, we can find extremal K\"ahler metrics $\widetilde{\omega}$ such that
\begin{equation}\label{omegaj dec}
\omega=\frac{a_j \sqrt{-1} dz^1 \wedge d\bar z^n}{2 |z^n|^2 log^2(|z^n|)}+p^* \widetilde{\omega} +O(|log(|z^n|)|^{-\eta})
\end{equation}
near $D_j$. Then  we can calculate that
\begin{equation*}
Ric_{\omega} = \frac{- \sqrt{-1} dz^n \wedge d \bar z^n}{2|z^n|^2 \log^2 (|z^n|)} + p^* Ric_{\widetilde{\omega}} +O(e^{-\eta t}).
\end{equation*}
Then we have that:
\begin{equation}\label{R dec}
\begin{split}
R_{\omega}&=2n \frac{Ric_{\omega}\wedge \omega^{n-1}}{\omega^n}\\
&=2n\frac{(n-1)p^* Ric_{\widetilde{\omega}}\wedge p^* \widetilde{\omega}^{n-2}\wedge (\frac{a_j \sqrt{-1} dz^n d \bar z^n}{2 |z^n|^2 \log^2 (|z^n|)})+(\frac{- \sqrt{-1} dz^n d \bar z^n}{2|z^n|^2 \log^2 (|z^n|)})\wedge p^* \widetilde{\omega}^{n-1}+O(e^{-\eta t})}{n p^* \widetilde{\omega}^{n-1}\wedge(\frac{a_j \sqrt{-1} dz^n d \bar z^n}{2|z^n|^2 \log^2 (|z^n|)})+O(e^{-\eta t})}\\
& =p^* R_{\widetilde{\omega}}-\frac{2}{a_j} +O(e^{-\eta t}),
\end{split}
\end{equation}
If $\omega$ is a cscK metric, then $\widetilde{\omega}$ is also a cscK metric according to the Lemma \ref{asym}. Thus we have that
\begin{equation*}
R_{\omega}= \underline{R}, \,\,\, R_{\widetilde{\omega}} = \underline{R}_{D_j}.
\end{equation*}
Then (\ref{R dec}) implies that $a_j = \frac{2}{\underline{R}_{D_j} - \underline{R}}$. In particular, $a_j$ depends only on  $[\omega]$, $X$ and $D_j$. We can see that the proof above can't be directly applied to the extremal K\"ahler metrics. Denote the constant $a_j$ corresponding to $\omega$ in (\ref{omegaj dec}) as $a_j(\omega)$. Our observation is that $a_j(\omega)$ depends on $[\omega]$, $X$ and $D_j$ even if $\omega$ is only an extremal K\"ahler metric instead of a cscK metric. Now we are ready to proof the Theorem \ref{a for extremal metric}:
\begin{proof}
(of the Theorem \ref{a for extremal metric}) First, we perform gauge fixing. Since we assume that $Aut_0(D)=\{Id\}$ and $D$ is a smooth divisor, we can follow the argument at the beginning of this section to find $g\in Aut_0^D(X)$ such that $Iso_0^D(X,g^* \omega_1)= Iso_0^D(X,\omega_2)$. Since $a_j(g^* \omega_1)= a_j(\omega_1)$, we can replace $\omega_1$ by $g^* \omega_1$ and assume that $Iso_0^D(X,\omega_1)= Iso_0^D(X,\omega_2)$. Then, we can get that the extremal K\"ahler vector field of $\omega_1$ and $\omega_2$ are the same according to the Proposition \ref{unique extremal vector field}. Denote their extremal K\"ahler vector field as $V$. 
Using the Lemma \ref{asym}, near $D_j$ we can write $\omega_i$ as
\begin{equation*}
\omega_i=\frac{a_j(\omega_i) \sqrt{-1} dz^1 \wedge d\bar z^n}{2|z^n|^2 log^2(|z^n|)}+p^* \widetilde{\omega_i} +O(|log(|z^n|)|^{-\eta}).
\end{equation*}
Here $\widetilde{\omega}_i$ is an extremal K\"ahler metric on $D$. Since we assume that $Aut_0(D)=\{Id\}$, we have that $\widetilde{\omega}_i$ is in fact a cscK metric and $\widetilde{\omega}_1= \widetilde{\omega}_2$, using the uniqueness of cscK metric on $D$.  Thus we have that 
\begin{equation*}
R_{\widetilde{\omega}_1}= R_{\widetilde{\omega}_2}.
\end{equation*}
According to (\ref{R dec}), in order to prove that $a_j(\omega_1)= a_j(\omega_2)$, it suffices to prove that $R_{\omega_1}=R_{\omega_2}$ on $D$. Then we can use the Lemma \ref{V different metric} to get that 
\begin{equation}\label{R1 R2 compare}
R_{\omega_1}- \underline{R}= R_{\omega_2}- \underline{R} + V(\varphi_1 - \varphi_2).
\end{equation}
Note that $V|_D$ is a holomorphic vector field on $D$. Since $Aut_0(D)=\{Id\}$, there is no nontrivial holomorphic vector field on $D$. As a result, $V|_D=0$. This implies that the norm of $V$ with respect to a Poincar\'e type metric converges to zero when we go to $D$.  Since $\varphi_1$ and $\varphi_2$ are two Poincar\'e type K\"ahler potentials, their derivatives with respect to a Poincar\'e metric is bounded. As a result $V(\varphi_1 - \varphi_2) |_D=0$. Then (\ref{R1 R2 compare}) implies that
\begin{equation*}
R_{\omega_1}= R_{\omega_2}
\end{equation*}
holds on $D$. According to the above argument, this concludes the proof of this Theorem.
\end{proof}
Define $E_{\beta}^{k,\alpha}= C_{\beta}^{k,\alpha} \oplus \Sigma_{i=1}^N \chi_i$. Here $\chi_i$ is a cut-off function supported in a small neighborhood of $D_i$ and it is equal to $1$ in a smaller neighborhood of $D_i$. For any $u\in u_1+ \Sigma_i \lambda_i \chi_i$, we define its norm as:
\begin{equation*}
||u||_{E_{\beta}^{k,\alpha}} \triangleq ||u_1||_{C_{\beta}^{k,\alpha}} + \Sigma_i |\lambda_i|.
\end{equation*}
We want to use the following Lemma proved by Auvray (See the Proposition 3.5 in \cite{A}):
\begin{lem}\label{delta estimate}
Let $(k,\alpha) \in \bN \times(0,1)$, $\eta \in C_{-\beta}^{k,\alpha}(\Lambda^{1,1})$ an exact $2-$form, $\beta>0$, and $\varphi$ the $\partial \bar \partial -$ potential of $\eta$ with zero mean with respect to some Poincar\'e type  K\"ahler metric $\omega$. Then $\varphi$ is in fact in $E_{\beta}^{k+2,\alpha}(\omega)$ and there exists a constant $C=C(\beta,k,\alpha,\omega)$ such that $||\varphi||_{E_{\beta}^{k+2,\alpha}}\le C||\eta||_{C_{\beta}^{k,\alpha}}.$ 
\end{lem}
Note that the definition of $C_{-\beta}^{k,\alpha}$ is the same as the definition of $C_{\beta}^{k,\alpha}$ in \cite{A}.
 The difference between two Poincar\'e type extremal K\"ahler metrics can be characterized as follows:

\begin{lem}\label{difference two metric}
Suppose that $D$ is a smooth  divisor and $Aut_0(D)=\{Id\}$. Let $\omega_3= \omega + dd^c \varphi_3, \omega_4= \omega+ dd^c \varphi_4$ be two Poincar\'e type extremal K\"ahler metrics in the same cohomology class. Then we have that
\begin{equation*}
\varphi_3 - \varphi_4  \in \widetilde{C}_{-\eta}^{\infty}.
\end{equation*}
\end{lem}
\begin{proof}
First we prove that
\begin{equation*}
\varphi_3 - \varphi_4 + \Sigma_{i=1}^N (a_j(\omega_3)- a_j(\omega_4))t \chi_i \in \widetilde{C}_{-\eta}^{\infty}.
\end{equation*}
In fact, using the Lemma \ref{asym}, near $D_j$ we can write $\omega_i$ as
\begin{equation*}
\omega_i=\frac{a_j(\omega_i) \sqrt{-1} dz^1 \wedge d\bar z^n}{2|z^n|^2 log^2(|z^n|)}+p^* \widetilde{\omega_{i,j}} +O(|log(|z^n|)|^{-\eta}),
\end{equation*}
for some metric $\widetilde{\omega_{i,j}}$ on $D_j$. Combining the above formula with (\ref{t coordinate 2}), we have that
\begin{equation*}
\omega_i = a_j(\omega_i) dd^c(-t) + p^*\widetilde{\omega_{i,j}} +O(e^{-\eta t}). 
\end{equation*}
Let $\varphi_{5,j}$ be a smooth function on $D_j$ such that
\begin{equation*}
\widetilde{\omega_{3,j}}= \widetilde{\omega_{4,j}}+ dd^c \varphi_{5,j}
\end{equation*}
holds on $D_j$.
Then we have that
\begin{equation*}
\begin{split}
dd^c (\varphi_{3,j} -\varphi_{4,j})=\omega_3 - \omega_4 & = (a_j(\omega_3)- a_j(\omega_4)) dd^c (-t) + p^* (\widetilde{\omega_{3,j}}- \widetilde{\omega_{4,j}}) +O(e^{-\eta t}) \\
 & = dd^c ( -(a_j(\omega_3)- a_j(\omega_4))t \chi_j  + \chi_j p^*\varphi_{5,j} ) +O(e^{-\eta t}).
\end{split}
\end{equation*}
Then we can use the Lemma \ref{delta estimate} to get that:
\begin{equation*}
\varphi_3 -\varphi_4 +\Sigma_{j=1}^N [(a_j(\omega_3)- a_j(\omega_4))t \chi_j - \chi_j p^*\varphi_{5,j}- \lambda_j \chi_j ] \in C_{-\eta}^{\infty}
\end{equation*}
for some constants $\lambda_j$. Using the Theorem \ref{a for extremal metric}, we have that $a_j(\omega_3)= a_j(\omega_4)$.  Then this concludes the proof of this Lemma.
\end{proof}

\subsection{Gauge fixing}

First, we want to fix the gauge for $\omega_1$.
Denote  $N_K$ as the normalizer of $K$
 in $Aut_0^D(X)$ consisting of $g\in Aut_0^D (X)$ such that $g K g^{-1}= K$. For any $g\in N_K$, we have that $g^* \omega_1 \in [\omega_1]$. So there exists a real-valued function $\varphi $ such that 
 \begin{equation}\label{e 8.1}
 g^* \omega_1 = \omega_{\varphi}= \omega_1+ \sqrt{-1} \partial \bar \partial \varphi,\,\,\, \int \varphi \omega_1^n =0.
 \end{equation}
Since $g^*\omega_1$ and $\omega_1$ are both Poincar\'e type extremal K\"ahler metrics, we can use the Lemma \ref{difference two metric} to get that $\varphi \in \widetilde{C}^{\infty}_{-\eta}$. For any $g\in N_K$, we have that $g^* \omega_1$ is $K-$invariant. As a result, we can get that the function $\varphi$ in (\ref{e 8.1})  is also $K-$invariant. Then we can define a map $\Psi^{\omega_1}$ from $N_K$ to $\widetilde{C}_{K,-\eta}^{\infty}$ by 
\begin{equation*}
\Psi^{\omega_1} (g)= \varphi.
\end{equation*} 
Define $S_{K,\omega_1}=\{g^* \omega_1: g\in N_K\}$. Then we can prove the following Lemma:
\begin{lem}\label{J proper}
Let $\omega$ be a Poincar\'e type  metric. Let $\omega_0= \omega+ dd^c \varphi_0$ with $||\varphi_0||_{L^{\infty}}< + \infty$ be another Poincar\'e type metric. Then $J_{\omega_0}$ is a proper functional over $S_{\omega}.$
\end{lem}
\begin{proof}
We can compute that for $\varphi \in \Psi^{\omega}(Aut_0^D(X))$:
\begin{equation*}
\begin{split}
J_{\omega_{0}} (\varphi)- J_{\omega}(\varphi) & = \frac{1}{n!} \Sigma_{p=1}^{n-1} \int_X \varphi \sqrt{-1} \partial \bar \partial \varphi_0 \wedge \omega^{n-p-1} \wedge \omega_{0}^p \\
& = \frac{1}{n!} \Sigma_{p=0}^{n-1} \int_X \varphi_0 \sqrt{-1} \partial \bar \partial \varphi \wedge \omega^{n-p-1} \wedge \omega_{0}^p \\
& = \frac{1}{n!} \Sigma_{p=0}^{n-1} \int_X \varphi_0 \omega_{\varphi} \wedge \omega^{n-p-1} \wedge \omega_{0}^p - \frac{1}{n!} \int_X \varphi_0 \omega^{n-p} \wedge \omega_{0}^p.
 \end{split}
\end{equation*}
Then we have that 
\begin{equation}\label{2j bdd}
|J_{\omega_{0}} (\varphi)- J_{\omega} (\varphi)| \le C_0 \sup_M |\varphi_0|,
\end{equation}
where $C_0$ is a constant independent of $\varphi$.

According to the section 3.6,  the functionals $J_{\omega}$ and $J_{\omega_0}$ are strictly convex along smooth Poincar\'e type geodesics. Note that $J_{\omega}$ has a critical point $\omega$ and $S_{\omega}$ is a finite dimensional space. we have that $J_{\omega}$ is proper on $S_{\omega}$. Combining this with (\ref{2j bdd}) and the fact that $S_{\omega}$ is a finite-dimensional space, we have that $J_{\omega_0}$ is also proper on $S_{\omega}$.
\end{proof}
 Next we prove the following proposition:
\begin{prop}\label{prop 8.4}
Let $\omega_1$ be a Poincar\'e type extremal K\"ahler metric. Then the image of the tangent space $(\Psi^{\omega_1})_*(T_{Id}{N_{K}})$ coincides with the space generated by the real-valued functions $f\in \widetilde{C}_{K,-\eta}^{\infty}$, such that $\nabla_{\omega_1}^{1,0} f \in \mathbf{h}^D_{\parallelsum}$ and $\int f \omega_1^n =0$.
\end{prop}
\begin{proof}
 Let $g_t$ be a smooth path in $N_K$, such that $g_0 =Id$ and such that the derivative $\frac{d g_t}{dt}|_{t=0}$ identifies with a holomorphic vector field which we denote by $X$. Since $g_t \in N_K$ and $\omega_1$ is a Poincar\'e type extremal K\"ahler metric and $K$-invariant, we have that $g_t^* \omega_1$ is also a Poincar\'e type extremal K\"ahler metric and $K$-invariant. Denote $\psi_t = \Psi^{\omega_1}(g_t)$. Using the Proposition \ref{reductivity}, we have that
\begin{equation*}
\mathbf{h}^D_{\parallelsum}=\mathbf{a}^D_{\parallelsum}(M) \oplus \nabla^{1,0} E_{0,r}\oplus \nabla^{1,0} E_{0,i}\oplus \Sigma_{\lambda \in Spec(\bar L|_{ N(L)}), \lambda>0} \nabla^{1,0}E_{\lambda}\end{equation*}
Here $L$ is defined using $g^* \omega_1$.
Then we can write
\begin{equation*}
X= X_a + \nabla^{1,0}(f_0 + \Sigma_{\lambda>0 } f_{\lambda}),
\end{equation*}
where $X_a \in \mathbf{a}^D_{\parallelsum}(M)$ and $f_{\lambda} \in E_{\lambda}$ for $\lambda \ge 0$. Since $g_t \in N_K$, we have that for any $\sigma \in K$, 
\begin{equation*}
g_t^* \sigma^* (g_t^{-1})^* \omega_1= \omega_1.
\end{equation*}
Differentiate with respect to $t$, we have that:
\begin{equation*}
\begin{split}
0 = (\frac{d}{dt} g_t^* \sigma^* (g_t^{-1})^* \omega_1)|_{t=0} &= \sigma^* (\frac{d}{dt}(g_t^{-1})^* \omega_1)|_{t=0} + (\frac{d}{dt} g_t^* \sigma^* \omega_1)|_{t=0}\\
& = \sqrt{-1} \partial \bar \partial [(f+ \bar f)- (f +\bar f)\circ \sigma],
\end{split}
\end{equation*}
where $f= f_0 + \Sigma_{\lambda >0} f_{\lambda}. $ As a result, for any $\sigma \in K$,
\begin{equation}\label{e 8.2}
 f + \bar f = ( f + \bar f ) \circ \sigma. 
\end{equation}
Apply $\bar L$ on both sides of $(\ref{e 8.2})$ for k times, we get that:
\begin{equation*}
\Sigma_{\lambda >0} \lambda^k( f_{\lambda} - f_{\lambda} \circ \sigma) =0.
\end{equation*}
Thus we infer that $f_{\lambda}= f_{\lambda} \circ \sigma$ for any $\sigma \in K$ and $\lambda>0$. Consider
\begin{equation*}
X \triangleq  Im (\nabla^{1,0} R_{\omega_1}) = \frac{-\sqrt{-1}}{2}[g^{\alpha \bar \beta} R_{\omega_1, \bar \beta}\frac{\partial }{\partial z_{\alpha}} - g^{\alpha \bar \beta}R_{\omega_1, \alpha} \frac{\partial}{\partial \bar z_{\beta}}].
\end{equation*}
Then $exp(t X)$ is a one parameter subgroup of $K$. Using the fact that $f_{\lambda}$ is $K-invariant$, we have that:
\begin{equation*}
0 =\frac{d}{dt}exp(t X)^* f_{\lambda} = X (f_{\lambda})=\frac{-\sqrt{-1}}{2}[R_{\omega_1, \bar \delta}f_{\lambda, \delta}- R_{\omega_1, \delta}f_{\lambda, \bar \delta}].
\end{equation*}
Then we have that:
\begin{equation*}
\lambda f_{\lambda} = \bar L f_{\lambda} = -(L - \bar L)f_{\lambda} = R_{\omega_1, \bar \delta} f_{\lambda, \delta} - R_{\omega_1,\delta} f_{\lambda, \bar \delta}=0.
\end{equation*}
Thus we have that $f_{\lambda}=0$ for any $\lambda>0$.  Thus
\begin{equation*}
X = X_a + \nabla^{1,0} f_0,
\end{equation*}
with $f_0 \in Ker L \cap Ker \bar L$ is a $K-$invariant complex-valued function. Therefore, $Re(f_0)$ and $Im (f_0)$ are both $K-$invariant and belong to $Ker L \cap Ker \bar L$. Differentiate $g_t^* \omega_1 = \omega_{\varphi_t}$ at $t=0$, we get that:
\begin{equation*}
\sqrt{-1} \partial \bar \partial \dot{\varphi}_0 = \sqrt{-1} \partial \bar \partial Re (f_0).
\end{equation*}
Using \ref{e 8.1}, we have that $\int \dot{\varphi}_0 \omega_1^n =0$. Without loss of generality, we can assume that $\int Re (f_0) \omega_1^n =0$. Thus, we can get that $\dot{\varphi}_0= Re (f_0)$. Then this Lemma follows immediately.
\end{proof}

\begin{lem}\label{gauge extk}
Suppose that $D$ is a smooth divisor. Suppose that $Aut_0(D)=\{Id\}$. Let $\omega_1, \omega_2$ be two Poincar\'e type extremal K\"ahler metrics. Then $J_{\omega_2}$ has a unique minimum and hence a critical point, $g^*\omega_1$, on $S_{K,\omega_1}.$ This implies that $d J_{\omega_2}|_{g^* \omega_1}$ annihilates all
real-valued functions $f\in \widetilde{C}_{K,-\eta}^{\infty}$, such that $\nabla_{\omega_1}^{1,0} f \in \mathbf{h}^D_{\parallelsum}$ and $\int f \omega_1^n =0$.
\end{lem}
\begin{proof}
Using the Lemma \ref{difference two metric}, we have that 
\begin{equation*}
||\varphi_1 -\varphi_2||_{L^{\infty}}< + \infty.
\end{equation*}
Then, we can use the Lemma \ref{J proper} to get that $J_{\omega_2}$ is proper on $S_{\omega_1}$. Since $S_{\omega_1}$ has a finite dimension, there is a critical point $g^* \omega_1 \in S_{\omega_1}$ which is a minimum point of $J_{\omega_2}$ on $S_{\omega_1}$. Since $J_{\omega_2}$ is strictly convex according to section 3.6, the critical point of $J_{\omega_2}$  on $S_{\omega_1}$ is unique.  The second part of this Lemma follows from the Proposition \ref{prop 8.4}
\end{proof}

\subsection{K-invairiant functions}

If $\omega_1$ is an extremal K\"ahler metric other than a cscK metric, it is possible that the Lichnerowicz operator with respect to $\omega_1$ may not be real-valued. However, this difficulty can be addressed by considering $K-$invariant functions. 

\begin{lem}\label{real rho}
Suppose that $\varphi$ is $K-$invariant. Then $h_{\omega_{\varphi}}^{X_1}$ is real-valued.
\end{lem}
\begin{proof}
We can compute that:
\begin{equation*}
\sqrt{-1}\bar \partial h_{\omega_{\varphi}}^{X_1} = \iota_{X_1} \omega_{\varphi} = \iota_{X_1} \omega_1 + \iota_{X_1}(\sqrt{-1} \partial \bar \partial (\varphi - \varphi_1))= \sqrt{-1} \bar \partial (R_{\varphi_1}+ X_1 (\varphi - \varphi_1)).
\end{equation*}
Thus 
\begin{equation*}
h_{\omega_{\varphi}}^{X_1} = R_{\varphi_1} + X_1(\varphi - \varphi_1)- \int_M (R_{\varphi_1}+ X_1 (\varphi - \varphi_1))\omega_{\varphi}^n.
\end{equation*}
The imaginary part of $h_{\omega_{\varphi}}^{X_1}$ is given by 
\begin{equation*}
Im (h_{\omega_{\varphi}}^{X_1}) = Im (X_1)(\varphi- \varphi_1)- \int_M Im (X_1) (\varphi- \varphi_1)\omega_{\varphi}^n.
\end{equation*}
Using the Lemma \ref{lem 7.6}, we know that $Im (X_1)$ is in the Lie algebra of $K$. Since $\omega_1$ is K-invariant and we assume that $\omega$ is K-invariant without loss of generality, $\varphi_1$ is also K-invariant. So
 $(\varphi- \varphi_1)$ is K-invariant, which implies that $Im (h_{\omega_{\varphi}}^{X_1})=0$.
\end{proof}

\begin{lem}\label{r real}
Suppose that $\varphi$ is $K-$invariant. Then we have that $R_{\bar \alpha} \varphi^{\bar \alpha}$ is real-valued.
\end{lem}
\begin{proof}
Denote $X= \nabla^{1,0}R$. Then we have that
\begin{equation*}
R_{\bar \alpha} \varphi^{\bar \alpha} = X(\varphi) = (Re X (\varphi) + \sqrt{-1} Im X (\varphi)).
\end{equation*}
By the Lemma \ref{lem 7.6}, $Im X$ lies in the Lie algebra of $K$. Thus we have that $ Im X (\varphi) =0$. This concludes the proof of this Lemma.
\end{proof}

\subsection{Proof of the Theorem \ref{main theorem}}

We define the functional $\mathcal{F}_K$ by the formula:
\begin{equation*}
\begin{split}
\mathcal{F}_K : \bR^n  \times \widetilde{C}^{5,\alpha}_{K,\delta} \times \bR  &\rightarrow \widetilde{C}^{1,\alpha}_{K,\delta} \times \bR \\
(\lambda, u, t_1) & \rightarrow R_{t \Sigma_{i=1}^N \chi_i  \lambda_i  + u}- \underline{R} -(1-t_1)(tr_{ t \Sigma_{i=1}^N \chi_i  \lambda_i   + u} \omega_2 -n ) - h_{\omega_{ t \Sigma_{i=1}^N \chi_i  \lambda_i  + u}}^{X_1}.
\end{split}
\end{equation*}
Here $\lambda=(\lambda_1,...,\lambda_N) \in \bR^n$. $\chi_i$ is a cut-off function which is supported in a small neighborhood of $D_i$ and is equal to $1$ in a smaller neighborhood of $D_i$.  Here $R_{ t \Sigma_{i=1}^N \chi_i  \lambda_i + u}$ means the scalar curvature of $\omega_1 + dd^c( t \Sigma_{i=1}^N \chi_i  \lambda_i )$. $tr_{ t \Sigma_{i=1}^N \chi_i  \lambda_i  + u} \omega_2$ means $tr_{\omega_1 + dd^c( t \Sigma_{i=1}^N \chi_i  \lambda_i +u)}\omega_2.$ We can replace $t$ by $t_K = \int_K t(gx) dg$ which is $K-$invariant. Since $K$ is a compact subgroup of $Iso^D(X,\omega_D)$, we can get that 
\begin{equation*}
t_K -t =O(e^{-t}).
\end{equation*}
Thus we can assume that $t$ and $\chi$ are $K-$invariant without loss of generality.

Denote 
\begin{equation*}
\mathcal{H}_{K, \delta, \varphi_1}= \{u \in \widetilde{C}_{K, \delta}^{5,\alpha}: Re L_{\omega_1}u=0\},
\end{equation*}
and
\begin{equation*}
\mathcal{H}_{K,\delta,\varphi_1,l}^{\bot}=\{u \in \widetilde{C}^{l,\alpha}_{K,\delta}: u \bot v, \text{ for any }v\in \mathcal{H}_{K, \delta, \varphi_1}\}.
\end{equation*}

Define a bilinear operator $B_{\varphi}(\cdot, \cdot)$:
\begin{equation*}
\begin{split}
B_{\varphi}(u,v) & \triangleq <\partial \bar \partial v, \partial \bar \partial \Delta_{\varphi} u>_{\varphi} + \Delta_{\varphi}<\partial \bar \partial v, \partial \bar \partial u >_{\varphi} + <\partial \bar \partial \Delta_{\varphi}v, \partial \bar \partial u>_{\varphi} \\
& + u_{,\bar \alpha p}v_{, \beta \bar p} (Ric_{\varphi})_{\alpha \bar \beta} + u_{, \bar p \beta} v_{, p \bar \alpha} (Ric_{\varphi})_{\alpha \bar \beta}.
\end{split}
\end{equation*}
In order to simplify the writing, we will use the following notation:
\begin{equation*}
<\partial v_1, \bar \partial v_2>_{\varphi} =\Sigma_{\alpha, \beta} \frac{\partial v_1}{\partial z_{\alpha}}\frac{\partial v_2}{ \partial z_{\bar \beta}} g^{\alpha \bar \beta}.
\end{equation*}
The following Lemma is due to Chen-Paun-Zeng \cite{CPZ}.  Since the proof of the Lemma is purely local, the proof in the Poincar\'e type case is the same as the proof in the smooth case.
\begin{lem}\label{lem 8.6}
Let $\omega_{\varphi}\in [\omega]$ be an extremal metric, and let $v, \xi$ be real-valued two smooth functions such that $L_{\varphi}v = \bar L_{\varphi} v =0$. Then we have the next identity,
\begin{equation*}
L_{\varphi}<\partial v, \bar \partial \xi>_{\varphi}= <\partial v, \bar \partial L_{\varphi} \xi>_{\varphi} + B_{\varphi}(v,\xi).
\end{equation*}
\end{lem}

Then we have the following Lemma:
\begin{lem}\label{lem deco}
Suppose that $D$ is smooth and connected. Suppose that $Aut_0(D)=\{Id\}$. Suppose that $\omega$ is a Poincar\'e type K\"ahler extremal K\"ahler metric. Denote $K=Iso_0^D(X,\omega)$. Then there exists a constant $\delta_1>0$ such that for any $\eta_0 \in (0, \delta_1)$,
\begin{equation*}
\widetilde{C}_{K,-\eta_0}^{1,\alpha} =  Ker L  |_{\widetilde{C}_{K,-\eta_0}^{5,\alpha}} \oplus  L(t\chi) \oplus L ( \widetilde{C}_{K,-\eta_0}^{5,\alpha})
\end{equation*}
\end{lem}
\begin{proof}
Note that $\omega$ is invariant under the holomorphic transformations in $K$. $t$ and $\chi$ can be assumed to be $K-invariant$ as well. Note that according to the Lemma \ref{r real}, $L=Re L$ when they act on $K-$invariant functions. Since $Aut_0(D)=\{Id\}$, we have that $Ker Re L_D=\{0\}$. Then this Lemma follows directly from the Proposition \ref{holder dec}.
\end{proof}

Now we are ready to prove the Theorem \ref{main theorem}:
\begin{proof}
(of the Theorem \ref{main theorem}). By differentiating the first term of $\iota_{X_1} \omega_{\varphi}= \sqrt{-1} \bar \partial h_{\omega_{\varphi}}^{X_1}$, we get that
\begin{equation*}
\bar \partial \dot{h}_{\omega_\varphi}^{X_1}= \bar \partial X_1(\dot{\varphi}).
\end{equation*}
This implies that $\dot{h}_{\omega_\varphi}^{X_1}- X_1(\dot{\varphi})$ is constant. On the other hand, by differentiating $\int_M h_{\omega_{\varphi}}^{X_1} \omega_{\varphi}^n =0$, we infer that we have
\begin{equation*}
\int_M (\dot{h}_{\omega_\varphi}^{X_1}+ h_{\omega_{\varphi}}^{X_1} \Delta_{\varphi}(\dot{\varphi})) \omega_{\varphi}^n =0.
\end{equation*}
Using Integration by parts and $\iota_{X_1} \omega_{\varphi}= \sqrt{-1} \bar \partial h_{\omega_{\varphi}}^{X_1}$, we get that:
\begin{equation*}
\int_M \dot{h}_{\omega_\varphi}^{X_1}- X_1(\dot{\varphi})\omega_{\varphi}^n =0.
\end{equation*}
Thus, we have $\dot{h}_{\omega_\varphi}^{X_1}= X_1(\dot{\varphi})$. Plugging in the definition of $X_1$, this is equivalent to 
\begin{equation}\label{e 8.4}
\dot{h}_{\omega_\varphi}^{X_1} =<\partial \dot{\varphi}, \bar \partial R_{\varphi}>_{\omega_{\varphi}}.
\end{equation}
Then we can calculate the derivative of $\mathcal{F}_K $ at $(\varphi_1,1)$:
\begin{equation}\label{e 8.5}
\begin{split}
d \mathcal{F}_K |_{(\varphi_1,1)}: \bR^n \times \widetilde{C}_{K,-\eta}^{5,\alpha} \times \mathbb{R}  & \rightarrow \widetilde{C}_{K,-\eta}^{1,\alpha}(M) \times \bR \\
(\lambda,u,s) & \rightarrow - L( u +t \Sigma_{i=1}^N \chi_i  \lambda_i )+s(tr_{\omega_1}\omega_2 -n).
\end{split}
\end{equation} 
Here we use the Lemma \ref{r real} to get that $<\partial \dot{\varphi}, \bar \partial R_{\varphi}>_{\omega_{\varphi}}$ is real-valued. Thus, $L u$ in (\ref{e 8.5}) is real-valued. 
According to the Lemma \ref{difference two metric}, the K\"ahler potentials of $\omega_1$ and $\omega_2$ are bounded from each other, we can use the Lemma \ref{gauge extk} to find $g\in N_K$ such that $g^* \omega_1$ is the minimum point of $J_{\omega_2}$ on $S_{K,\omega_1}$. From now on we replace $\omega_1$ by $g^* \omega_1$ and assume that $\omega_1$ is the minimum point of $J_{\omega_2}$ on $S_{K,\omega_1}$.  As a result,
\begin{equation*}
tr_{\omega_1}\omega_2 - n \in \mathcal{H}_{K,-\eta,\varphi_1,5}^{\bot}.
\end{equation*}
Then we can define the following map:
\begin{equation*}
\begin{split}
\Pi : \bR^n \times (\mathcal{H}_{K, -\eta, \varphi_1} \oplus \mathcal{H}_{K,-\eta,\varphi_1,5}^{\bot}) \times \bR & \rightarrow \mathcal{H}_{K, -\eta, \varphi_1} \oplus \mathcal{H}_{K,-\eta,\varphi_1,1}^{\bot} \times \bR \\
(\lambda, u, w, t_1) &\rightarrow (u+ \pi_2 \circ \mathcal{F}_K( \lambda  , u + w,t_1), t_1)
\end{split}
\end{equation*}
Here $\pi_2$ is the projection to $\mathcal{H}_{K,-\eta,\varphi_1,1}^{\bot}$.
The derivative of $\Pi$ at $(0,0,0,1)$ is:
\begin{equation*}
d \Pi (\lambda,u,w,s)=(u- L_{\omega_1}(\Sigma_{i=1}^N \lambda_i \chi_i t )- L_{\omega_1}(w)+s(tr_{\omega_1}\omega_2 -n),s).
\end{equation*}
Using the Lemma \ref{holder dec} and the assumption $Aut_0(D)=\{Id\}$, we have that $d \Pi$ is a bijection at $(0,0,0,1)$. Then we can use the implicit function theory to get that there exists $\epsilon_0>0$ such that for any  $||u||\le \epsilon_0$ and $|t_1 -1|\le \epsilon_0$, there exists $\widetilde{t}(u,t_1)$, $\widetilde{u}(u,t_1)$, $w(u, t_1)$  and $\lambda (u,t_1)$ such that
\begin{equation*}
\begin{split}
&\Pi(\lambda(u,t_1), \widetilde{u}(u,t_1), w(u,t_1), \widetilde{t}(u,t_1))\\
&=(\widetilde{u}(u,t_1) +\pi_2 \circ \mathcal{F}_K(\lambda(u,t_1), \widetilde{u}(u,t_1)+ w(u,t_1), \widetilde{t}(u,t_1)),\widetilde{t}(u,t_1)) \\
&= (u+0, t_1).
\end{split}
\end{equation*}
This implies that $\widetilde{u}=u$, $\pi_2 \circ \mathcal{F}_K(\lambda(u,t_1), \widetilde{u}(u,t_1)+ w(u,t_1), \widetilde{t}(u,t_1))=0$  and $\widetilde{t}=t_1$. Then we can get that:
\begin{equation}\label{pi2f01}
\pi_2 \circ \mathcal{F}_K(\lambda (u,t_1), u + w(u, t_1),t_1) =0.
\end{equation}
Consider the functional 
\begin{equation*}
P(u,t_1) \triangleq \pi_1 \circ \mathcal{F}_K(\lambda(u,t_1), u+ w(u,t_1),t_1),
\end{equation*}
where $\pi_1$ is the projection onto the factor $\mathcal{H}_{K, -\eta, \varphi_1,0}$. Here we define 
\begin{equation*}
\mathcal{H}_{K,-\eta, \varphi_1,0}= \{u \in \mathcal{H}_{-\eta, \varphi}: \int u \omega_{\varphi_1}^n =0\}.
\end{equation*}
We want to solve the equation
\begin{equation*}
P(u_{t_1},t_1)=0
\end{equation*}
for each $1-\epsilon_1 < t_1 \le 1$ with $u_{t_1} \in \mathcal{H}_{K,-\eta, \varphi_1,0}$.
Denote $\psi(u,t_1)\triangleq \Sigma_{i=1}^N \lambda_i(u,t_1)t\chi_i + w(u,t_1)$. Take the derivative of (\ref{pi2f01}) with respect to $t_1$, we get that
\begin{equation*}
- L_{\omega_1} \frac{\partial \psi}{\partial t_1}|_{(0,1)} + tr_{\omega_1} \omega_2 -n=0.
\end{equation*}
Take the derivative of (\ref{pi2f01}) with respect to $u$, we get that
\begin{equation*}
0=L_{\omega_1}(\frac{\partial \psi}{\partial u}|_{(0,1)}(v)) = L_{\omega_1}(\Sigma_{i=1}^N \frac{\partial \lambda_i}{\partial u}|_{(0,1)}(v)t \chi_i + \frac{\partial w}{\partial u}|_{(0,1)}(v))
\end{equation*}
for any $v\in \mathcal{H}_{K, -\eta, \varphi_1,0}$. Using the Lemma \ref{lem 7.21}, we have that $\Sigma_{i=1}^N \frac{\partial \lambda_i}{\partial u}|_{(0,1)}(v)t \chi_i  + \frac{\partial w}{\partial u}|_{(0,1)}(v)$ must be bounded which implies that $\frac{\partial \lambda_i}{\partial u}|_{(0,1)}(v)=0$. Thus we have that 
\begin{equation*}
L_{\omega_1}(\frac{\partial w}{\partial u}|_{(0,1)}(v))=0,
\end{equation*}
which implies that 
\begin{equation*}
\frac{\partial w}{\partial u}|_{(0,1)}(v)=0,
\end{equation*}
since $\frac{\partial w}{\partial u}|_{(0,1)}(v) \in \mathcal{H}_{K,-\eta,\varphi_1,l}^{\bot}$.  Thus we have that
\begin{equation*}
\frac{\partial \psi}{\partial u}|_{(0,1)}(v)=0.
\end{equation*}

We claim that $P(u,1)=0$ for any $u\in \{v\in \mathcal{H}_{K, -\eta, \varphi_1,0}: ||v||\le \epsilon_1 \}$ with $\epsilon_1$ to be a small constant. In fact, consider the corresponding holomorphic transformation $g_u$ of $u$ according to the Proposition \ref{prop 8.4}. Then we have that $g_u^* \omega$ is also a Poincar\'e type cscK metric. This implies that $P(u,1)=0$.
Then, we can define
\begin{equation*}
\widetilde{P}(u,t_1) \triangleq \frac{P(u,t_1)}{t_1-1}
\end{equation*}
and it can extended as a continuous function on $\mathcal{H}_{K,-\eta, \varphi_1,0} \times [0,1]$, because of the equality 
\begin{equation*}
\widetilde{P} (u,1) =\lim_{t_1 \rightarrow 1-} \frac{P(u,t_1)}{t_1-1} =\frac{\partial P}{\partial t_1}|_{(u,1)}.
\end{equation*}
It suffices to solve the equation $\widetilde{P}(u_{t_1},t_1)=0$ by showing that $\frac{\partial \widetilde{P}}{\partial u}|_{(0,1)}$ is invertible. First we write 
\begin{equation*}
\begin{split}
\widetilde{P}(u,1) & = \frac{\partial}{\partial t_1} P|_{(u,1)} = \pi_1 [- L_{ u + \psi_{u,1}} \frac{\partial \psi}{\partial t_1}|_{(u,1)} + tr_{u +\psi_{u,1}}\omega -n] \\
& = \pi_1 [-\Delta_{u + \psi_{u,1}}^2 \frac{\partial \psi}{\partial t_1}|_{(u,1)} - (\frac{\partial \psi}{\partial t_1}|_{(u,1)})_{,\bar \alpha \beta} (Ric_{u + \psi_{u,1}})_{\alpha \bar \beta}\\
& + tr_{ u + \psi_{u,1}}\omega -n].
\end{split}
\end{equation*}
We compute 
\begin{equation*}
\begin{split}
\frac{\partial}{\partial u} \widetilde{P}|_{(0,1)}(v) & =\pi_1 \{ <\partial \bar \partial v, \partial \bar \partial \Delta_{\varphi_1} \xi>_{\varphi_1}  + \Delta_{\varphi_1}<\partial \bar \partial v, \partial \bar \partial \xi>_{\varphi_1} + <\partial \bar \partial \Delta_{\varphi_1}v , \partial \bar \partial \xi>_{\varphi_1} + \xi_{,\bar \alpha p} v_{, \bar p \beta}(Ric_{\varphi_1})_{\alpha \bar \beta}\}\\
& + \xi_{\bar p \beta} v_{,p \bar \alpha} (Ric_{\varphi_1})_{\alpha \bar \beta} -<\partial \bar \partial v, \chi>_{\varphi_1} - L_{\omega_1} \frac{\partial^2 \psi}{\partial u \partial t_1}|_{(0,1)} (v)\} \\
& = \pi_1 [B_{\varphi_1}(v, \xi)- <\partial \bar \partial v, \chi>_{\varphi_1}],
\end{split}
\end{equation*}
where $\xi= \frac{\partial \psi}{\partial t_1}|_{(0,1)}$ and $B_{\varphi_1}(v, \xi)$ is the operator in Lemma \ref{lem 8.6}. Then we can use the above formula and the Lemma \ref{lem 8.6} to get that:
\begin{equation*}
\begin{split}
\frac{\partial}{\partial u} \widetilde{P}|_{(0,1)}(v) &= \pi_1 [L_{\omega_1}(<\partial v, \bar \partial \xi>_{\varphi_1})-<\partial v,\bar \partial L_{\omega_1}\xi>-<\partial \bar \partial v, \omega>_{\varphi_1}]\\
&= \pi_1 (-<\partial v, \bar \partial (tr_{\omega_1}\omega-n)>_{\varphi_1}- <\partial \bar \partial v, \omega>_{\varphi_1}).
\end{split}
\end{equation*}
Then we can see that
\begin{equation*}
\begin{split}
\int \frac{\partial \widetilde{P}}{\partial u}|_{(0,1)}(v)v \omega_{\varphi}^n &= \int (- <\partial v, \bar \partial (tr_{\omega_1}\omega -n)>_{\varphi_1} v - <\partial \bar \partial v, \omega>_{\varphi_1} v)\omega_{\varphi_1}^n \\
& = \int v_{,\bar \alpha} v_{, \beta} \omega_{\alpha \bar \beta} \omega_{\varphi_1}^n \ge 0,
\end{split}
\end{equation*}
Since $\int v \omega_{\varphi_1}^n =0$, the integral above  is positive and is equal to zero if and only if $v=0$. Therefore, $\frac{\partial \widetilde{P}}{ \partial u}|_{(0,1)}$ is injective and therefore bijective. Then we can use the implicit function theorem to get $u_t$ such that $P(u_{t_1},t_1)=0$ for $t_1$ sufficiently close to 1. Combining this with (\ref{pi2f01}), we have that 
\begin{equation*}
\mathcal{F}_K(\lambda(u_{t_1},t_1),  u_{t_1} + w(u_{t_1},t_1),t_1) =C_{t_1}
\end{equation*}
for some constant $C_{t_1}.$ Since the integral of $\mathcal{F}_K(\lambda(u_{t_1},t_1),  u_{t_1} + w(u_{t_1},t_1),t_1)$ with respect to $\omega_{\Sigma_i \lambda_i(u_{t_1},t_1) t \chi_i + u_{t_1}+w(u_{t_1},t_1)}^n$ is $0$, we have that $C_{t_1}=0$. This concludes the proof of this theorem.
\end{proof}

\subsection{Energy functional $\mathcal{E}_V$}

Before we prove the Theorem \ref{unique extk}, we want to study the following functional: For any $V\in \mathbf{h}^D_{\parallelsum, \mathbb{C}}$, we define an associated energy functional $\mathcal{E}_V$ by letting:
\begin{equation*}
d \mathcal{E}_V |_{\omega} (\dot u) = \int_{X \setminus D} \dot{u} h^V_{\omega} \omega^n.
\end{equation*}
Then, we can prove the following proposition:
\begin{prop}\label{mathcal ev}
Let $\omega_u = \omega+ dd^c u \in \mathcal{PM}_{\Omega,V}$ depends smoothly on two real parameters $s$ and $t$. Assume that $\omega \in \mathcal{PM}_{\Omega,V}$ and $u$ is invariant under $ImV$. Then we have that:
\begin{equation*}
(\frac{d}{ds}) \int_X \dot{u_t} h_{\omega_u}^V \omega_u^n = \int_X (\ddot{u}_{st}- (\partial \dot{u}_t, \partial \dot{u}_s)_{\omega_u})h^V_{\omega_u} \omega_u^n,
\end{equation*}
where $(,)_{\omega_u}= Re<,>_{\omega_u}$ is the real scalar product defined by $\omega_u$. 
 \end{prop}
 \begin{proof}
This proposition follows from the Proposition 4.14 in \cite{BB}. Since both $\omega$ and $u$ are invariant under $ImV$, we have that $h^V_{\omega_u}$ is real valued. Then the proposition follows using integration by parts which is due to the Lemma \ref{G-S}.
 \end{proof}

\begin{lem}\label{mathcal ev linear}
$\mathcal{E}_V$ is linear along Poincar\'e type $C^{1,1}$ geodesics in $\mathcal{PM}_{\Omega,V}$.
\end{lem}
\begin{proof}
Using the Proposition \ref{mathcal ev}, we have that $\mathcal{E}_V$ is a well defined function and
\begin{equation*}
\frac{d^2}{dt^2} \mathcal{E}_V (u) = \int_X (\ddot{u}_{tt}- |\bar \partial \dot{u}_t|^2_{\omega_u}) h^V_{\omega_u} \omega_u^n. 
\end{equation*}
This implies that $\mathcal{E}_V$ is linear on smooth Poincar\'e type geodesic in $\mathcal{PM}_{\Omega,V}$. By approximation (c.f. the section 3 in \cite{XZ}), we have that $\mathcal{E}_V$ is linear along Poincar\'e type $C^{1,1}$ geodesics in $\mathcal{PM}_{\Omega,V}$.
\end{proof}

\subsection{Proof of the Theorem  \ref{unique extk}}

Now we are ready to prove the Theorem \ref{unique extk}.
\begin{proof}
(of the Theorem \ref{unique extk}). First, we want to prove that there exists $g\in Aut_0^D(X)$ such that $g^* \omega_1 = \omega_2$, under the assumption that the K\"ahler potentials of $\omega_1$ and $\omega_2$ are bounded from each other. We first fix the gauge. Using the Theorem \ref{compact isometry} and the fact that the maximal compact connected subgroups of $Aut_0(X)$ are conjugate (using a result by Matsuchima), we can assume that $Iso_0^D(X,\omega_1)=Iso_0^D(X,\omega_2)=K$ by replacing $\omega_1$ with $g^* \omega_1$ for an appropriate map $g\in Aut_0(X)$. According to the Lemma \ref{difference two metric}, the K\"ahler potentials of $\omega_1$ and $\omega_2$ are bounded from each other. Then we can use the Lemma \ref{gauge extk} to replace $\omega_1$ by $g_1^* \omega_1$ and assume that $\omega_1$ is the minimum point of the functional $J_{\omega_2}|_{S_{K, \omega_1}}$ by gauge fixing. Since $J_{\omega_2}$ is strictly convex on $S_{K, \omega_2}$ and $\omega_2$ is a critical point of $J_{\omega_2}$. We have that $\omega_2$ is the minimum point of $J_{\omega_2}$ on $S_{K, \omega_2}$. Since $g_1\in N_K$, we still have 
\begin{equation*}
Iso_0^D(X,\omega_1)=Iso_0^D(X,\omega_2).
\end{equation*}
Using the Proposition \ref{unique extremal vector field}, we have that $\nabla_{\omega_1}^{1,0}(R_{\omega_1})= \nabla_{\omega_2}^{1,0}(R_{\omega_2})$ which we denote as $X$. Then we can use the proof of the Theorem \ref{main theorem} to get two paths of twisted extremal metrics, $\varphi_{k,t_1}$ with $\varphi_{k,1}= \varphi_k$ for $k=1,2$ satisfying
\begin{equation}\label{e 8.8}
R_{\varphi_{k,t_1}} - \underline{R}- \rho_{\varphi_{k,t_1}}(X) - (1-t_1)(tr_{\varphi_{k,t_1}}\omega_2 -n)=0.
\end{equation}
Define the modified K-energy on $\mathcal{PM}_{[\omega]}$ by:
\begin{equation*}
\frac{d \mathcal{E}_K}{dt} = \int_M (-(R_{\varphi}- \underline{R})+h^X_{\omega_{\varphi}}(X))\frac{d\varphi}{dt} \omega_{\varphi}^n.
\end{equation*}
$\mathcal{E}_K$ can be written as 
\begin{equation*}
\mathcal{E}_K = \mathcal{M} + \mathcal{E}_X
\end{equation*}
According to \cite{XZ}, we have that $ \mathcal{M}$ is weakly convex along any K-invariant Poincar\'e type $C^{1,1}$ geodesic. Combining this with the Lemma \ref{mathcal ev linear}, we have that $\mathcal{E}_K$  is weakly convex along any K-invariant Poincar\'e type $C^{1,1}$ geodesic. Note that   $J_{\omega_2}$ is strictly convex along K-invariant Poincar\'e type $C^{1,1}$ geodesic. As a result, for any $t_1\in (0,1)$,
\begin{equation*}
\mathcal{E}_K +(1-t_1)J_{\omega_2}
\end{equation*}
is strictly convex along any K-invariant Poincar\'e type $C^{1,1}$ geodesic. Note that any two K-invariant K\"ahler metrics can be connected by a K-invariant Poincar\'e type $C^{1,1}$ geodesic. As a result, the critical point of $\mathcal{E}_K +(1-t_1)J_{\omega_2}$ is unique. Since a solution to (\ref{e 8.8}) is a critical point of $\mathcal{E}_K +(1-t_1)J_{\omega_2}$, we have that the solution to (\ref{e 8.8}) is unique. As a result, we have that $\varphi_{1,t_1}= \varphi_{2,t_1}$. As $t_1 \rightarrow 1$, we get that $\varphi_1 = \varphi_2$. This concludes the proof of this Theorem. 
\end{proof}

\end{document}